%% file: Paper_for_arXiv_Recycling_Subspaces_for_coupledAdjoint_v2.tex
\documentclass[preprint]{elsarticle}

\usepackage[utf8]{inputenc}
\usepackage[T1]{fontenc}
\usepackage{array}
\usepackage{multicol}
\usepackage{multirow}
\usepackage{tabularray}
\usepackage{ragged2e}
\usepackage{makecell}
\usepackage{multirow}
\usepackage{diagbox}
\usepackage{longtable}
\usepackage{booktabs}
\usepackage{graphicx}
\usepackage{float}
\usepackage{siunitx}
\usepackage{amsmath,amssymb,amsthm,mathtools}
\usepackage{empheq}
\usepackage{hyperref}
\usepackage[textheight = 21cm, textwidth = 16cm]{geometry}
\renewcommand{\arraystretch}{0.9}

\usepackage{subcaption}
\usepackage{caption}
\usepackage[export]{adjustbox}

\newtheorem{thm}{Theorem}[section]

\newdefinition{rmk}{Remark}[section]
\newdefinition{definition}{Definition}[section]
\newproof{pf}{Proof}

\usepackage[nodisplayskipstretch]{setspace}   

\usepackage[svgnames,dvipsnames,table]{xcolor}

\usepackage[ruled, linesnumbered]{algorithm2e}
\SetKwComment{Comment}{$\triangleright$\ }{}
\SetAlgoHangIndent{0em}    

\SetCommentSty{footnotesize}
\SetNlSty{texttt}{}
\SetAlgoLined
\DontPrintSemicolon

\input{macros.tex} 

\DeclareMathAlphabet\mathbfcal{OMS}{cmsy}{b}{n}  

\newcommand{\ols}[1]{\mskip.5\thinmuskip\overline{\mskip-.5\thinmuskip {#1} \mskip-.5\thinmuskip}\mskip.5\thinmuskip} 

\begin{document}
	
	\begin{frontmatter}
		
		\title{Recycling Krylov Subspaces for Efficient Partitioned Solution of Aerostructural Adjoint Systems}
		
		\author[1]{Christophe Blondeau\corref{cor1}\fnref{fn1}}
		\ead{Christophe.Blondeau@onera.fr}
		\author[2]{Mehdi Jadoui\fnref{fn2}}

        \affiliation[1]{organization={ONERA}, addressline={Université Paris-Saclay, F-92322 Châtillon}, country={France}}
		\cortext[cor1]{Corresponding author}
        \fntext[fn1]{{Aerodynamics, Aeroelasticity, Acoustics Dept., DAAA.}}
        \fntext[fn2]{{At the time of writing, Mehdi Jadoui was a PhD student at ONERA.}}        
\begin{abstract}
	
Robust and efficient solvers for coupled-adjoint linear systems are crucial to successful aerostructural optimization. Monolithic and partitioned strategies can be applied. The monolithic approach is expected to offer better robustness and efficiency for strong fluid-structure interactions. However, it requires a high implementation cost and convergence may depend on appropriate scaling and initialization strategies. On the other hand, the modularity of the partitioned method enables a straightforward implementation while its convergence may require relaxation. In addition, a partitioned solver leads to a higher number of iterations to get the same level of convergence as the monolithic one. 

The objective of this paper is to accelerate the fluid-structure coupled-adjoint partitioned solver by considering techniques borrowed from approximate invariant subspace recycling strategies adapted to sequences of linear systems with varying right-hand sides. Indeed, in a partitioned framework, the structural source term attached to the fluid block of equations affects the right-hand side with the nice property of quickly converging to a constant value. We also consider deflation of approximate eigenvectors in conjunction with advanced inner-outer Krylov solvers for the fluid block equations.
We demonstrate the benefit of these techniques by computing the coupled derivatives of an aeroelastic configuration of the ONERA-M6 fixed wing in transonic flow. For this exercise the fluid grid was coupled to a structural model specifically designed to exhibit a high flexibility. All computations are performed using RANS flow modeling and a fully linearized one-equation Spalart-Allmaras turbulence model. Numerical simulations show up to 39\% reduction in matrix-vector products for GCRO-DR and up to 19\% for the nested FGCRO-DR solver.
			
\end{abstract}
		

\end{frontmatter}

\section{Introduction}

We are interested in robust and efficient solvers for the solution of the coupled-adjoint linear system which is crucial to successful aerostructural optimization.
Fine control of an aerodynamic shape or of a structural layout leads to a high-dimensional parameter space. For high-fidelity simulations, gradient-based optimizers in conjunction with the adjoint approach are the methods of choice. However, the coupled adjoint linear system is inherently ill-conditioned as it embeds matrix blocks of different scales and structures. In addition, the fluid block coming from the exact linearization of the RANS equations associated to a turbulence model is often very stiff. Besides, a strong level of fluid-structure interaction is known to be detrimental to the robustness and efficiency of existing solution techniques.   
For such linear systems, monolithic and partitioned (or segregated) strategies can be applied. 

The former approach is expected to offer better robustness and efficiency for strong fluid-structure interactions \cite{heil2004, heil2008}. However, it requires a high implementation cost and convergence may depend on appropriate scaling and initialization strategies. On the other hand, the modularity of the partitioned method enables a straightforward implementation while its convergence may require relaxation. In addition, a partitioned solver leads to a higher number of iterations to get the same level of convergence as the monolithic one. A review of partitioned simulations of fluid-structure interactions involving black-box solvers is proposed in \cite{degroote2013partitioned}.

The partitioned approach simply consists in solving, in an alternating way, the aerodynamic and the structural sub-problems by applying the Linear Block Gauss-Seidel algorithm (LBGS). It accounts for the interdisciplinary coupling by adding a source term to the right-hand side of each set of disciplinary adjoint equations. The modularity of this approach makes it rather interesting since it takes advantage of the specific routines designed for each sub-problem and does not require a high implementation cost. Nevertheless, this approach becomes rapidly inefficient and could even diverge for strong fluid-structure coupling even though the addition of some level of relaxation helps to mitigate this issue. On the other hand, the monolithic approach consists in solving the fluid and structural equations simultaneously making it more robust in the sense that it is less sensitive to the strength of the fluid-structure interaction. The coupled adjoint system is generally solved by using Krylov subspace methods. The challenging aspect of such an approach is then to develop advanced preconditioning strategies combined with numerical ingredients so that Krylov methods reach the best performances in terms of robustness and efficiency.

The objective of this work is to improve the efficiency of the existing partitioned solver \cite{Achard2018, Jadoui_AIAA2022} by considering techniques borrowed from Krylov subspace recycling strategies adapted to sequences of linear systems with varying right-hand sides \cite{parks2006recycling}.
We will demonstrate the benefit of these advanced techniques by computing the coupled derivatives for the ONERA-M6 fixed wing in transonic flow \cite{AileM6}. For this exercise the fluid grid is coupled to a structural model specifically designed to exhibit a high flexibility. All computations are performed using RANS flow modeling and a fully linearized one-equation Spalart-Allmaras turbulence model.

\begin{figure}[H]
	\begin{subfigure}{.5\textwidth}
		\raggedright
		\includegraphics[trim={0.5cm 0.5cm 1.5cm 2.0cm},clip,width=0.95\textwidth]{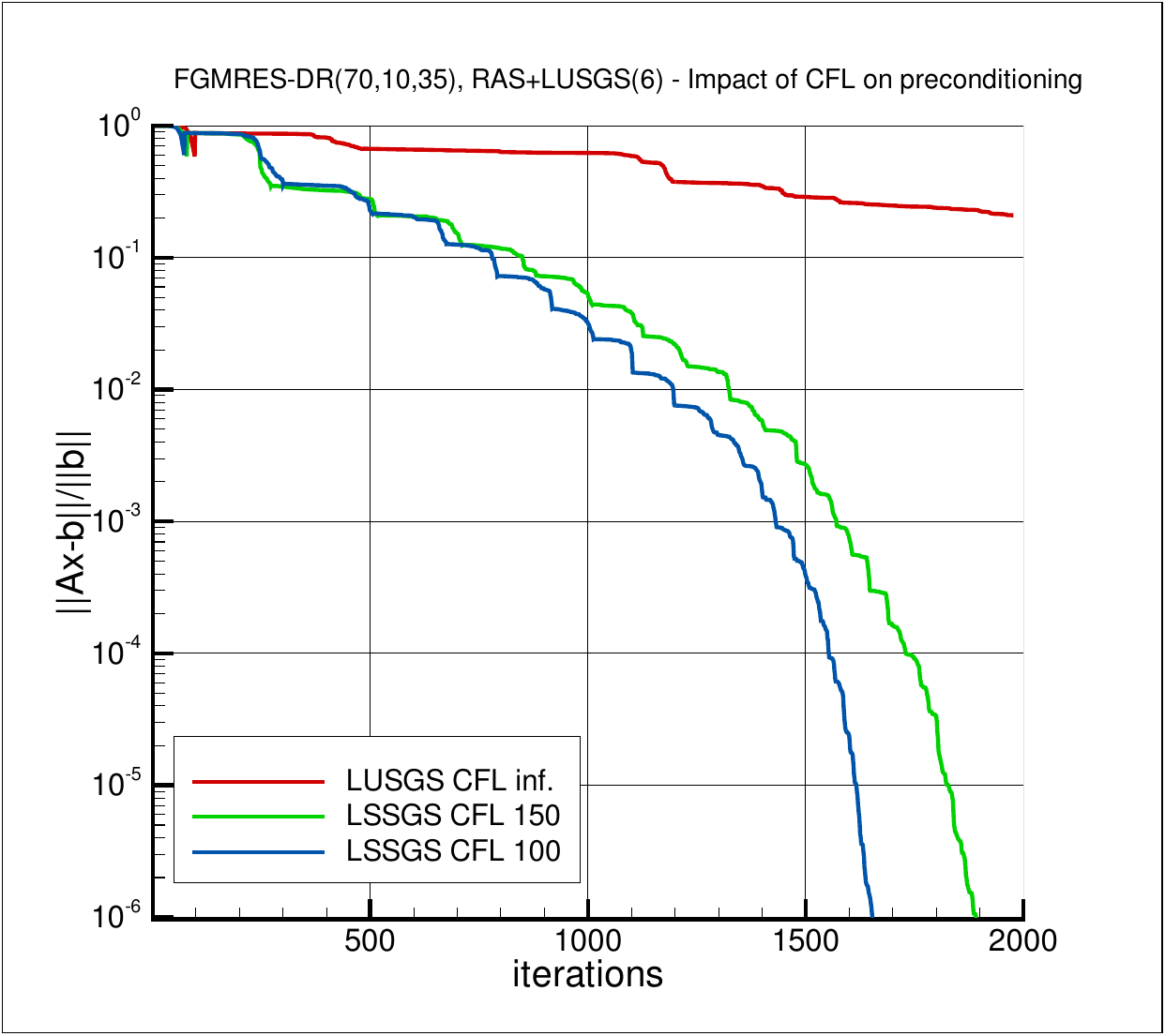}
	\end{subfigure}%
	\begin{subfigure}{.5\textwidth}
		\centering
		\includegraphics[trim={0.5cm 0.5cm 1.5cm 2.0cm},clip,width=0.95\textwidth]{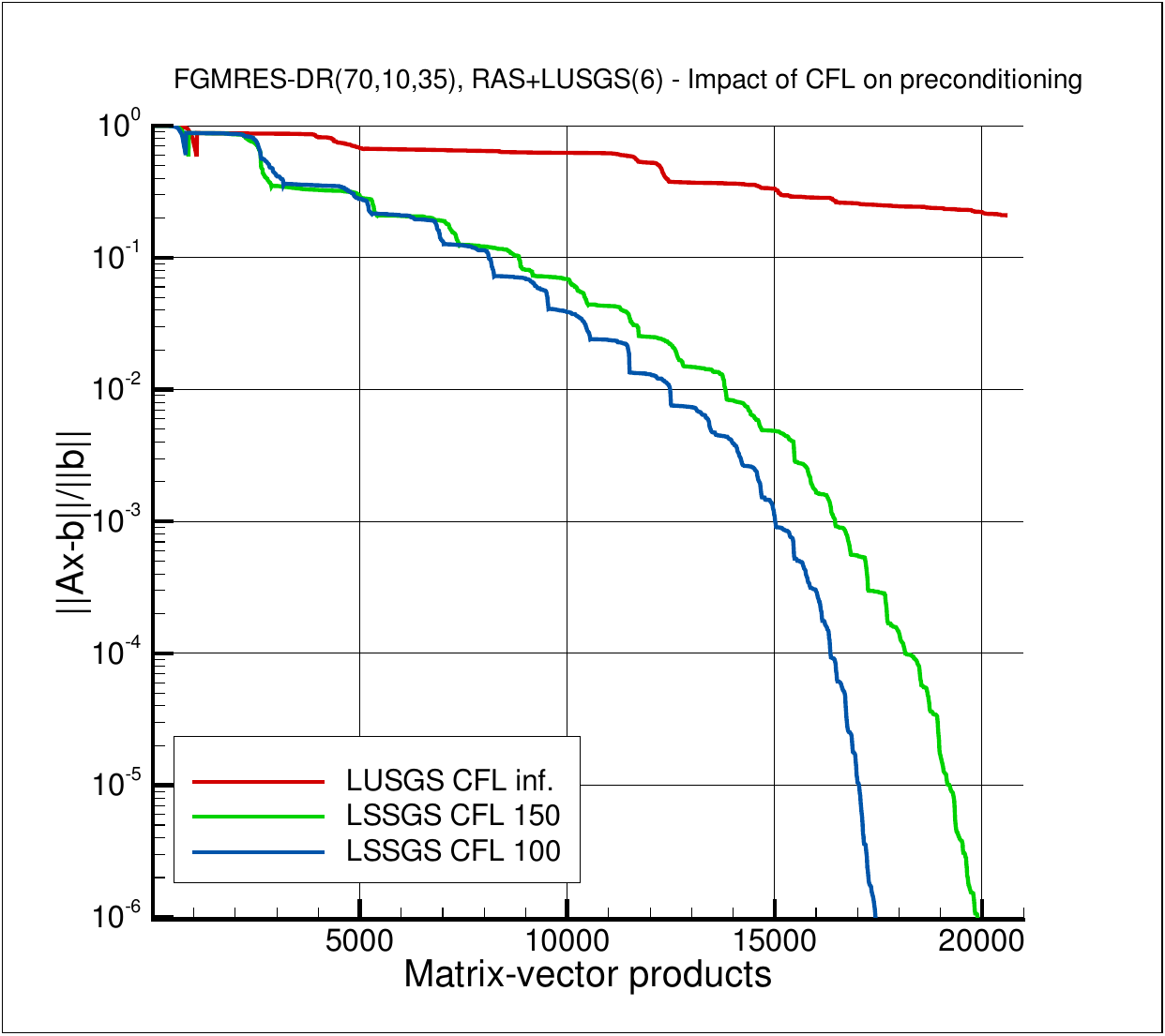}
	\end{subfigure}
	\caption{Coupled-adjoint relative residual norm convergence history of FGMRES-DR(70,10,35). Impact of CFL on the performance of the LU-SGS preconditioner.}
	\label{fig:FGMRESDR_70_10_35_SAlin_LUSGS_6_CFL_effect_iterations}
\end{figure}

As an example, Figure \ref{fig:FGMRESDR_70_10_35_SAlin_LUSGS_6_CFL_effect_iterations} illustrates the performance of the existing embedded coupled-adjoint segregated solver for structured grids in a former version of the elsA CFD code \footnote{elsA and is the joint property of ONERA and Safran.}, applied to the M6 wing test case.

These convergence curves will serve as a reference for comparison with the improved solver efficiency allowed by the developments performed in this work.
A flexible GMRES Krylov solver with deflated restarting is used with an outer Krylov basis of 70 vectors, an inner Krylov basis of size 10 and a deflation subspace of 35 vectors (between fluid cycles). At the beginning of this work, the only available preconditioner for structured grids in elsA was a combination of a Restrictive Additive Schwarz (RAS) domain decomposition method coupled with a Lower-Upper Symmetric Gauss Seidel (LU-SGS) relaxation. This preconditioner is controlled by two parameters: the number of relaxations and the CFL coefficient. The number of relaxations is fixed at 6 and the CFL is varied.
In \cite{jadoui2022comparative} we point out that a proper tuning of the CFL coefficient is crucial for this type of preconditioners and how new block ILU type preconditioners accelerate the rate of convergence.

Also, after each fluid-structure coupling, a cold restart is performed and the spectral information from the previous fluid-structure cycle is discarded. As a consequence, a plateau is observed after each restart which dramatically hampers convergence in the first cycles. We will demonstrate that well-chosen subspace recycling strategies will eliminate these convergence stagnations.

\bigbreak
In a high-fidelity aerostructural optimization context, Zhang and Zingg \cite{zhang2018efficient} implemented a robust monolithic solution method for both aerostructural analysis and coupled adjoint problem. A three-field formulation was adopted involving the mesh, the flow and the structural states. The performance of the monolithic method as well as the partitioned one was investigated through a comparative study by varying the level of fluid-structure coupling. For the coupled adjoint problem solution, a GCROT \cite{de1999truncation} Krylov solver has been used in conjunction with a block Gauss-Seidel preconditioner. The monolithic adjoint solution has been 60$\%$ more efficient than the partitioned one for strong coupling. For weak fluid-structure coupling, the monolithic solution still outperformed the partitioned one with a better efficiency of 40$\%$. In terms of computational time, the monolithic method showed 50 $\%$ to over 60 $\%$ faster than the partitioned method. A similar comparative study of both monolithic and partitioned approaches was performed by Kenway et al. \citep{kenway2014scalable} except that a Block Jacobi preconditioner was applied to the coupled adjoint system. The aerodynamic and structural block preconditioners were solved by a preconditioned Krylov method (restarted Generalized Minimal RESidual - GMRES \citep{saad1986}) and a direct factorization method respectively. A Flexible Krylov method (e.g FGMRES \citep{saad1993flexible}) has been used for the coupled adjoint system solution. The Common Research Model (CRM) wing-body-tail configuration was sized by considering two critical load cases: 1$g$ cruise condition with moderate elastic deformation and a 2.5$g$ pull-up with significantly more deflection. For the same memory footprint, the best monolithic solution seems to outperform the best partitioned one by reducing the time by 19$\%$ for the 1$g$ load and by 29$\%$ for the 2.5$g$ case. These numerical experiments demonstrate the great benefit of using monolithic approach in the strong coupling case but at the price of a robust preconditioner for the Krylov solver. We note however that both studies only considered inviscid flow modeling. These conclusions about the monolithic solver efficiency might be mitigated by the added stiffness of adjoint system matrices produced by a RANS fluid model associated to a linearized turbulence model.
\bigbreak
Although the satisfactory performance of monolithic solvers, advanced strategies that could accelerate the partitioned algorithm using black-box solvers have received less attention. We recall that a partitioned algorithm consists in approximately solving the aerodynamic adjoint block at each fluid-structure iteration resulting in a sequence of adjoint linear systems with varying right-hand sides. As already mentioned, the structural source term that affects the right-hand side of the fluid block has the nice property of rapidly converging to a constant value. The corollary of this property is that after several fluid-structure couplings, the subsequent fluid systems should greatly benefit from recycling spectral information from the previous fluid-structure cycles. At the start of this work, the current partitioned adjoint solver did not take advantage of recycling and at each update of the structural source term the Krylov solver did a cold start from the previous solution.

The principle of deflation is to remove the influence of a system's subspace on the iterative process. This is usually beneficial when directions of certain subspaces hamper convergence. Deflation of an eigenspace can be performed in two ways: the linear system (matrix and right-hand side) is left-multiplied by a projector $\mathbf{P}$, i.e. deflation by projection \cite{vuik1999efficient, tang2007deflation, Gaul2014}, or some eigenvectors are added to the Krylov subspace, i.e. deflation by augmentation \citep{morgan2002gmres}. A survey of deflation and augmentations techniques can be found in \citep{coulaud:hal-00803225}. A specific type of deflation preconditioning aims at solving a rank-deficient projected system, using a Krylov solver, in a certain subspace outside of the problematic subspace. The solution is then complemented with the solution in the latter subspace. In \citep{Gaul2014} the author makes the post-correction superfluous by using a projection as right-preconditioner instead. For deflation by augmentation, adding eigenvectors to the Krylov subspace can effectively deflate corresponding eigenvalues from the spectrum because when these directions are included in the solution approximation, the convergence of the Krylov solver continues according to the modified spectrum. The deflation by augmentation led to the well-known FGMRES-DR solver \citep{giraud2010} and its extension to inner-outer Krylov solver \citep{simoncini2002flexible,jadoui2022comparative}.
\bigbreak
Unfortunately, the deflated restarted GMRES framework based on subspace augmentation is not adapted for solving sequences of linear systems \citep{parks2006recycling}. Fortunately, some authors have proposed new strategies in order to reuse information accumulated in previous fluid-structure cycles and use it to accelerate the solution of the next linear system. Krylov subspaces recycling methods seem to be the suitable choice. Historically, De Sturler suggested the Generalized Conjugate Residuals with inner Orthogonalization (GCRO) method \citep{de1996nested}, an improved version of the recursive GMRES (GMRESR) solver \citep{van1994gmresr} by maintaining an orthogonality condition between the outer and the inner spaces generated by GMRESR. This way, it provides the optimal correction to the solution in a global search space. Later, Parks et al. formulated the GCRO-DR algorithm \citep{parks2006recycling} that combines GCRO and deflation techniques by augmentation introduced by Morgan \citep{morgan2002gmres}. They demonstrated better performances of GCRO-DR compared to the GMRES-DR in a  long sequence of linear systems from a fracture mechanics problem. Carvhalo et al. extended GCRO-DR to the flexible case (FGCRO-DR) \citep{carvalho2011flexible} and they conducted in-depth analysis of both flexible methods. In particular, they showed that both methods can be algebraically equivalent if a certain colinearity condition is satisfied at each cycle. In 2013, Niu et al. introduced Loose GCRO-DR (LGCRO-DR) \citep{niu2013accelerated} for improving the convergence of GCRO-DR by recycling both spectral information and approximate error information. The error is defined as the distance between the current iterate and the exact solution of the system. It is not known by definition but a fair approximation to it can be computed. This idea was initially proposed by Baker et al. \citep{baker2005technique} and mimics the idea behind GMRESR of including approximations to the error in the current approximation space. This error information is interesting since it represents in some sense the previous Krylov space generated in the previous cycle and subsequently discarded. In addition to that, LGCRO-DR is straightforward and economic to implement. 
\bigbreak
In this work, we investigate advanced Krylov subspace methods using subspace recycling strategies for accelerating the partitioned solver applied to the linear coupled-adjoint system. More specifically, we compare GMRES-DR and FGMRES-DR to GCRO-DR and FGCRO-DR with and without subspace recycling. This work benefits from the recent achievements to improve efficiency of the fluid adjoint solution by applying nested Krylov subspace methods \citep{jadoui2022comparative}. The numerical experiments are performed on an aeroelastic configuration of the ONERA M6 fixed wing in transonic viscous flow.

This paper is organized as follows. In section  \ref{Aerostructuraladjointsystem} we briefly recall the theoretical background of aeroelastic and coupled-adjoint equations. The partitioned algorithm is also outlined. The aeroelastic numerical test case is presented in section \ref{ONERA_M6}.
To support later comparison with GCRO-DR anf FGCRO-DR, we then review the fundamentals of FGMRES-DR in section \ref{FGMRESDR_Review} with a focus on numerical implementation and application to the solving of the fluid and coupled adjoint systems. Section \ref{GCRODR_theory} is devoted to the description of the GCRO algorithm with adaptations related to variable preconditioning and subspace recycling. We also take the opportunity to give some insights related to an efficient implementation of the deflation and recycling strategy. The numerical experiments are then repeated with the GCRO-DR solver and show very promising reduction in terms of matrix-vector products compared to the standard implementation. Finally, the flexible case is presented in section \ref{FGCRODR_theory}.


\section{Aerostructural adjoint system} \label{Aerostructuraladjointsystem}


\subsection{Aeroelastic equilibrium}
Let us denote the state variables of the coupled system $\mathbf{W}$ and $\mathbf{U}$, representing the fluid conservative variables and the structural displacements respectively. At the aeroelastic equilibrium, the state variables and the meshes satisfy the discretized equations of fluid and structural mechanics simultaneously:


\begin{equation}
	\left\{ 
	\begin{aligned} 
	\mathbf{R}_a(\mathbf{X}_a,\mathbf{W},\mathbf{U}) &= \mathbf{0} \\
	\mathbf{R}_s(\mathbf{X}_s,\mathbf{W},\mathbf{U}) &= \mathbf{0}
	\end{aligned} 
	\right.
	\label{eq1}
\end{equation}
where $\mathbf{R}_a$ is the discrete aerodynamic residual and $\mathbf{R}_{s}$ the discrete structural residual. These two blocks of equations are coupled through aerodynamic forces $\mathbf{Q}_{a}$ loading the skin of the structure and the structural displacements $\mathbf{U}$ deforming the fluid mesh. The structural mesh is noted $\mathbf{X}_s$. In the following we introduce two aerodynamic grids $\mathbf{X}_a$ and $\mathbf{X}_{a0}$. $\mathbf{X}_a$ denotes the deformed aerodynamic grid at the aeroelastic equilibrium at the outcome of the aeroelastic analysis. $\mathbf{X}_{a0}$ is called the reference mesh which supports the aerodynamic shape parametrization. Typically for an aircraft design study the reference mesh is chosen as the jig shape or the flight shape in reference nominal flight conditions. The load, displacement and mesh deformation operators then merely depend on $\mathbf{X}_{a0}$ for an aeroelastic or coupled-adjoint analysis. The structural loads $\mathbf{Q}_s$ are obtained with a suitable load transfer technique applied to $\mathbf{Q}_a$ such that

\begin{equation}
\mathbf{Q}_s(\mathbf{Q}_a(\mathbf{W},\mathbf{X}_a),\mathbf{X}_{a0},\mathbf{X}_s) = \mathbf{T}^{Q}_{surf}(\mathbf{X}_{a0},\mathbf{X}_s) \mathbf{Q}_a(\mathbf{W},\mathbf{X}_a)
\end{equation}

\noindent where $\mathbf{T}^{Q}_{surf}$ represents a linear load transfer operator. The subscript $surf$ stipulates that the associated linear operators or data relate to the fluid-structure interface. The structural displacements alter the fluid grid locations through the relation:

\begin{equation}
\mathbf{X}_{a} = \mathbf{X}_{a0} + \delta\mathbf{X}_{a}(\delta\mathbf{X}_{a,surf},\mathbf{X}_{a0}) = \mathbf{X}_{a0} + \mathbf{T}_{vol}(\mathbf{X}_{a0})\delta\mathbf{X}_{a,surf}
\end{equation}

\noindent with $\mathbf{T}_{vol}(\mathbf{X}_{a0})$ the volume operator performing the deformation of the fluid domain. The vector $\delta\mathbf{X}_{a,surf}$ corresponds to the displacements of the fluid nodes at the fluid-structure interface.

\begin{equation}
\delta\mathbf{X}_{a,surf} = \delta\mathbf{X}_{a,surf}(\mathbf{X}_{a0},\mathbf{X}_{s},\mathbf{U}) = \mathbf{T}^{U}_{surf}(\mathbf{X}_{a0},\mathbf{X}_{s})\mathbf{U}
\end{equation}

\noindent where $\mathbf{T}^{U}_{surf}(\mathbf{X}_{a0},\mathbf{X}_{s})$ represents a linear displacement transfer operator.


\subsection{Partitioned strategy for the coupled adjoint system}

Let us consider a scalar aeroelastic objective function $J(\mathbf{W},\mathbf{U},\mathbf{X}_{a},\mathbf{X}_{s})$ and a design parameter $p$. One way to obtain the coupled adjoint equations is to formulate an augmented objective function by adding the total variation of the residuals $\mathbf{R}_{s}$ and $\mathbf{R}_{a}$ to the total derivative  $dJ/dp$. More specifically, we define $d\tilde{J}/dp$ as
\begin{align}
	\frac{d\tilde{J}}{dp} &= \frac{dJ}{dp} + \lambda^{T}_{a} \frac{d \mathbf{R}_{a}}{dp} + \lambda^{T}_{s} \frac{d\mathbf{R}_{s}}{dp},
	\label{eq2}
\end{align}
\noindent where
\begin{equation}
	\frac{dJ}{dp} = \frac{\partial J}{\partial \mathbf{W}} \frac{d \mathbf{W}}{dp} + \frac{\partial J}{\partial \mathbf{X}_{a}} \frac{d \mathbf{X}_{a}}{dp} + \frac{\partial J}{\partial \mathbf{X}_{s}} \frac{d \mathbf{X}_{s}}{dp} + \frac{\partial J}{\partial \mathbf{U}} \frac{\partial \mathbf{U}}{\partial p}.
	\label{eq3}
\end{equation}

In Eq. (\ref{eq2}) the total variations of residuals are exactly zero since they represent constraints related to the satisfaction of the equilibrium equations at the outcome of the aeroelastic analysis. For simplicity, we restrict here to the specific case of a shape design parameter not affecting the structural geometry nor the structural stiffness. In addition, the explicit dependency of the objective function with respect to the structural states is dropped, i.e.,  we consider only derivatives of aerodynamic coefficients. The full derivation for the general case can be found in \citep{Achard2017}. It is worth to mention that these assumptions do not lead to any loss of generality of the work presented in this paper since we focus on solution techniques for the adjoint system. As $\mathbf{X}_{s}$ does not depend on the design parameter $p$, we have $d \mathbf{X}_{s}/dp = \mathbf{0}$. Under the same assumption we also have $d\mathbf{K}/dp=\mathbf{0}$. After some algebra manipulation we end up with the following expression for $d\tilde{J}/dp$ in which the total derivatives $d\mathbf{W}/dp$ and $d\mathbf{U}/dp$ have been factored out:
\begin{equation}
\begin{split}
\frac{d \tilde{J}}{dp} & = \left( \frac{\partial J}{\partial \mathbf{W}} + \lambda^{T}_{a} \frac{\partial \mathbf{R}_{a}}{\partial \mathbf{W}} - \lambda^{T}_s\mathbf{C} \right) \frac{d \mathbf{W}}{dp} + \left( \frac{\partial J}{\partial \mathbf{X}_{a}} \mathbf{A} + \lambda^{T}_{a} \frac{\partial \mathbf{R}_{a}}{\partial \mathbf{X}_{a}}\mathbf{A} + \lambda^{T}_{s}(\mathbf{K} - \mathbf{D})\right) \frac{d \mathbf{U}}{dp} \\
& + \left( \frac{\partial J}{\partial \mathbf{X}_{a}}  + \lambda^{T}_{a} \frac{\partial \mathbf{R}_{a}}{\partial \mathbf{X}_{a}} \right) \mathbf{B} \frac{d \mathbf{X}_{a0}}{dp}-\lambda^{T}_{s}\mathbf{E}\frac{d \mathbf{X}_{a0}}{dp}
\end{split}
\label{eq4}
\end{equation}

Constant matrices $\mathbf{A}$ to $\mathbf{E}$ are defined analytically with the following formulas (see \citep{Achard2017,Achard2018}):
\begin{align}
\mathbf{A} &= \mathbf{T}_{vol} \mathbf{T}^{U}_{surf} \\
\mathbf{B} &= \frac{\partial \mathbf{X}_a}{\partial \mathbf{X}_{a0}} = \mathbf{I} + \frac{\partial \mathbf{A}}{\partial \mathbf{X}_{a0}} \mathbf{U}\\
\mathbf{C} &= \mathbf{T}^{Q}_{surf} \frac{\partial \mathbf{Q}_a}{\partial \mathbf{W}} \\
\mathbf{D} &= \mathbf{T}^{Q}_{surf} \frac{\partial \mathbf{Q}_a}{\partial \mathbf{X}_a} \mathbf{T}^{U}_{surf}  \\
\mathbf{E} &= \mathbf{T}^{Q}_{surf} \frac{\partial \mathbf{Q}_a}{\partial \mathbf{X}_a}\mathbf{B} + \frac{\partial \mathbf{Q}_s}{\partial \mathbf{X}_{a0}}
\label{eq5}
\end{align}

The coupled adjoint linear system is obtained by canceling factors related to $d\mathbf{W}/dp$ and $d\mathbf{U}/dp$ in Eq. (\ref{eq4}) to give
\begin{equation}
\begin{bmatrix}
\left[ {\frac{\partial \mathbf{R}_{a}}{\partial \mathbf{W}}}\right]^{T} &
-\mathbf{C}^{T} \\
             &                  \\
\mathbf{A}^{T}\left[\frac{\partial \mathbf{R}_{a}}{\partial \mathbf{X}_{a}}\right]^{T} & \mathbf{K}^{T} - \mathbf{D}^{T}
\end{bmatrix}
\begin{bmatrix}
\lambda_{a} \\
            \\
\lambda_{s}
\end{bmatrix} =
\begin{bmatrix}
- \left[ \frac{\partial J}{\partial \mathbf{W}} \right]^{T} \\
        				\\
-\mathbf{A}^{T} \left[ \frac{\partial J}{\partial \mathbf{X}_{a}} \right]^{T}
\end{bmatrix}
\label{Coupled_Adjoint_system}
\end{equation}

The process for solving the adjoint system follows an iterative block scheme. Algorithm~\ref{alg:LBGS-algorithm} details the Linear Block Gauss-Seidel (LBGS) scheme applied for the solution of system (\ref{Coupled_Adjoint_system}). In this derivation, we use the structural flexibility $\mathbf{S}$ which is a small reduced matrix relating the set of structural forces to the set of structural displacments pertaining to the fluid-structure coupling. The relaxation factor $\theta_{s}$ has been introduced on the adjoint vector $\lambda_{s}$. Assuming that the coupled system is solved to machine accuracy, the total derivative reconstruction is given by

\begin{equation}
\frac{d J}{dp} = \left( \frac{\partial J}{\partial \mathbf{X}_{a}}  + \lambda^{T}_{a} \frac{\partial \mathbf{R}_{a}}{\partial \mathbf{X}_{a}} \right) \mathbf{B} \frac{d \mathbf{X}_{a0}}{dp}-\lambda^{T}_{s}\mathbf{E}\frac{d \mathbf{X}_{a0}}{dp}
\label{dJdp}
\end{equation}

In the expression above the computation of the product of the geometrical sensitivities with the matrix $\mathbf{B}$ is not trivial. If one already has at hand a linearized version of the operator $\mathbf{A}$, i.e. of $\mathbf{T}_{surf}^U$ and $\mathbf{T}_{vol}$, it can be applied to $d \mathbf{X}_{a0}/{dp}$ as many times as the number of design variables. This is the most straightforward manner but the benefit of the adjoint formulation is then mitigated by the cost of the gradient assembly. The other way is to transpose the first term in the right-hand side of Eq. (\ref{dJdp}) and compute products like $[ \partial\mathbf{T}_{surf}^U / \partial\mathbf{X}_{a0} ]^T \mathbf{v}$ and $[ \partial\mathbf{T}_{vol} / \partial\mathbf{X}_{a0} ]^T \mathbf{v}$, where $\mathbf{v}$ has fluid grid size. We call this mode the geometrical adjoint of the mesh deformation and displacement transfer operators. These two modes of gradient assembly have been implemented in the coupled-adjoint module of the elsA software. The iterations stop after a maximum number $n_{cpl}$ of fluid-structure couplings or when the relative fluid and structural residuals get lower than a prescribed tolerance: $r_A \le \epsilon_A$ and $r_S \le \epsilon_S$. Typically we choose $\epsilon_A = \epsilon_S = 10^{-6}$.
\vspace{\baselineskip}

\begin{algorithm}[H]
\setstretch{1.1}
\fontsize{10}{14}\selectfont
\caption{Partitioned LBGS strategy for coupled-adjoint solution}
\label{alg:LBGS-algorithm}
{
\KwData{ $\mathbf{U, W}$ $\mathbf{X}_{a}, \mathbf{X}_{s}, \mathbf{X}_{a0}, \lambda^{0}_{a}, \lambda^{0}_{s}, \theta_{s}, \epsilon_A, \epsilon_S, n_{cpl}$ }
$\mathbf{RHS}_{stru} \gets \mathbf{0}$ \;
\lIf{$\lambda^{0}_{s} \neq \mathbf{0} $}
{
$\mathbf{RHS}_{stru} \gets \left( \mathbf{T}^{Q}_{surf} \frac{\partial \mathbf{Q}_{a}}{\partial \mathbf{W}}  \right)^T \mathbf{S}^{T}\lambda^{0}_{s}$   \Comment*[f]{Restart from a previous structural adjoint solution}
}
\lIf{$\lambda^{0}_{a} = \mathbf{0}$}
{
$\left[ \frac{\partial \mathbf{R}_a}{\partial \mathbf{W}} \right]^T \lambda^{0}_{a} = - \left[ \frac{\partial J}{\partial \mathbf{W}} \right]^{T} + \mathbf{RHS}_{stru}$ \Comment*[f]{Approximate solution of the fluid adjoint problem}
}

\For{$k \leftarrow 1, n_{cpl}$}
{
$A_{Xs,surf} \leftarrow \left( \frac{\partial \mathbf{Q}_a}{\partial \mathbf{X}_{a,surf}} \right)^{T}\left( \mathbf{T}^{Q}_{surf} \right)^T \mathbf{S}^{T}\lambda^{k-1}_{s}$ \Comment*[f]{Structural geometric adjoint}

$ A_{Xa} \leftarrow - \left( [\lambda^{k-1}_{a}]^{T} \frac{\partial \mathbf{R}_{a}}{\partial \mathbf{X}_{a}} + \frac{\partial J}{\partial \mathbf{X}_{a}}  \right)^{T} $ \Comment*[f]{Aerodynamic geometric adjoint}

$A_{Xa,surf} \leftarrow (\mathbf{T}_{vol})^{T} \mathbf{A}_{Xa}$ \Comment*[f]{Mesh deformation adjoint}

$\lambda^{k}_{s} \leftarrow (\mathbf{T}^{U}_{surf})^{T}(\mathbf{A}_{Xs,surf} + \mathbf{A}_{Xa,surf})$ \Comment*[f]{Structural adjoint vector}

$\lambda^{k}_{s} \leftarrow (1-\theta_{s}) \lambda^{k-1}_{s} + \theta_{s} \lambda^{k}_{s}$ \Comment*[f]{Relaxation (optional)} \label{alg:Relaxation_step}

$\mathbf{RHS}_{stru} \leftarrow \left( \mathbf{T}^{Q}_{surf} \frac{\partial \mathbf{Q}_{a}}{\partial \mathbf{W}} \right)^{T}\mathbf{S}^{T} \lambda^{k}_{s}$ \Comment*[f]{Update structural rhs}

$\left[ \frac{\partial \mathbf{R}_a}{\partial \mathbf{W}} \right]^T \lambda^{k}_{a} = - \left[ \frac{\partial J}{\partial \mathbf{W}} \right]^{T} + \mathbf{RHS}_{stru}$ \Comment*[f]{Approximate solution of the fluid adjoint problem}

\lIf{($r_A \le \epsilon_A$ and $r_S \le \epsilon_S)$} 
{
	goto 14  \Comment*[f]{check norm of fluid and structure relative residuals}
}
}
$\frac{dJ}{dp} \leftarrow \left( \frac{\partial J}{\partial \mathbf{X}_{a}}  + \lambda^{T}_{a} \frac{\partial \mathbf{R}_{a}}{\partial \mathbf{X}_{a}} \right) \mathbf{B} \frac{d \mathbf{X}_{a0}}{dp}-\lambda^{T}_{s}\mathbf{E}\frac{d \mathbf{X}_{a0}}{dp}$  \label{assembly} \Comment*[f]{Objective function gradient assembly}

}
\end{algorithm}

\section{ONERA-M6 wing aeroelastic analysis} \label{ONERA_M6}

\begin{figure}[H]
	\centering
	\parbox[h][8.0cm][c]{.5\textwidth}{
		\includegraphics[trim={0.3cm 0.cm 1.cm 0.5cm},clip,width=0.49\textwidth]{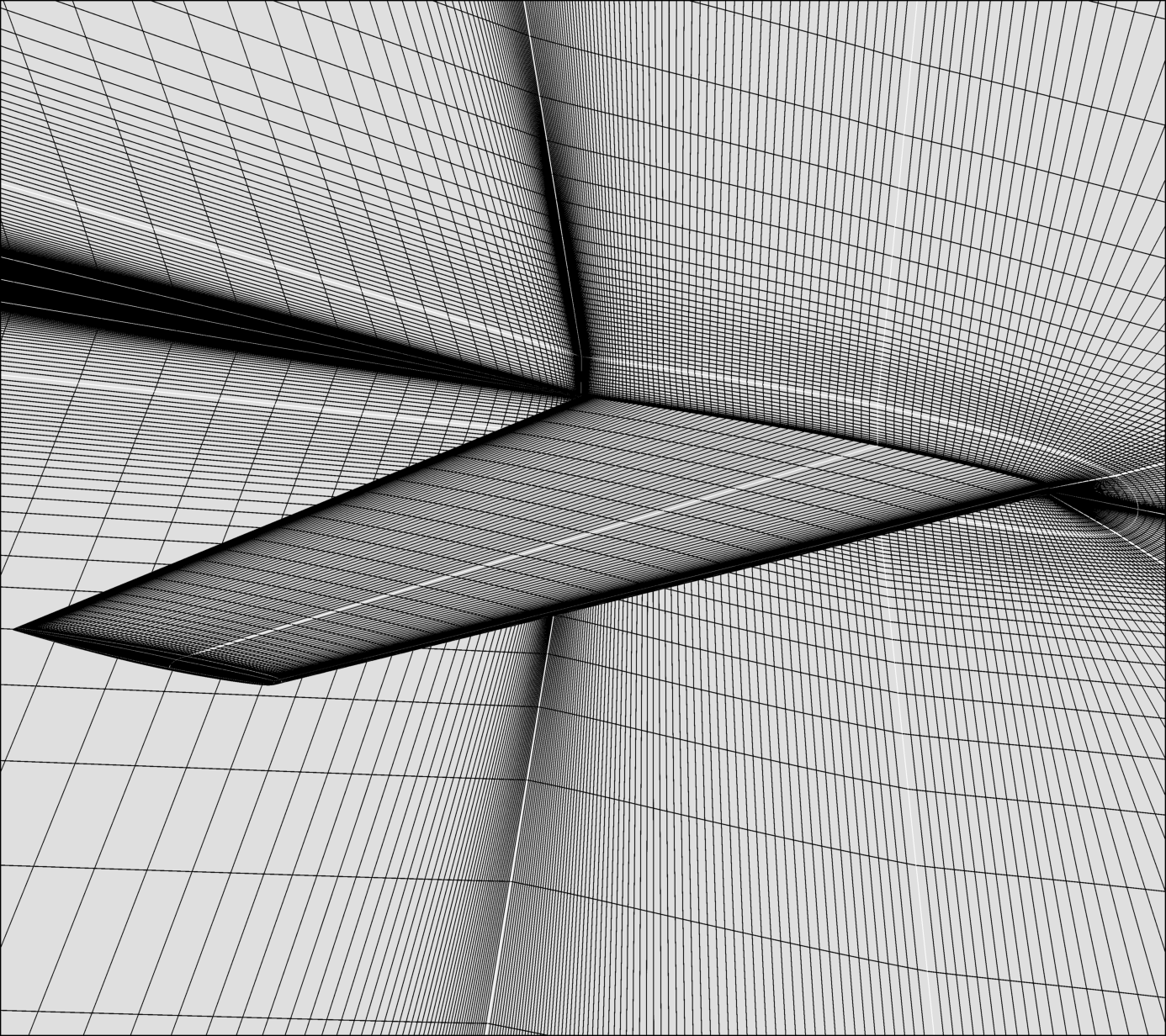}
	}\hfill
	\parbox[h][8.0cm][c]{.5\textwidth}{
		\includegraphics[trim={1.cm 7.5cm 0.2cm 7.5cm},clip,width=0.49\textwidth]{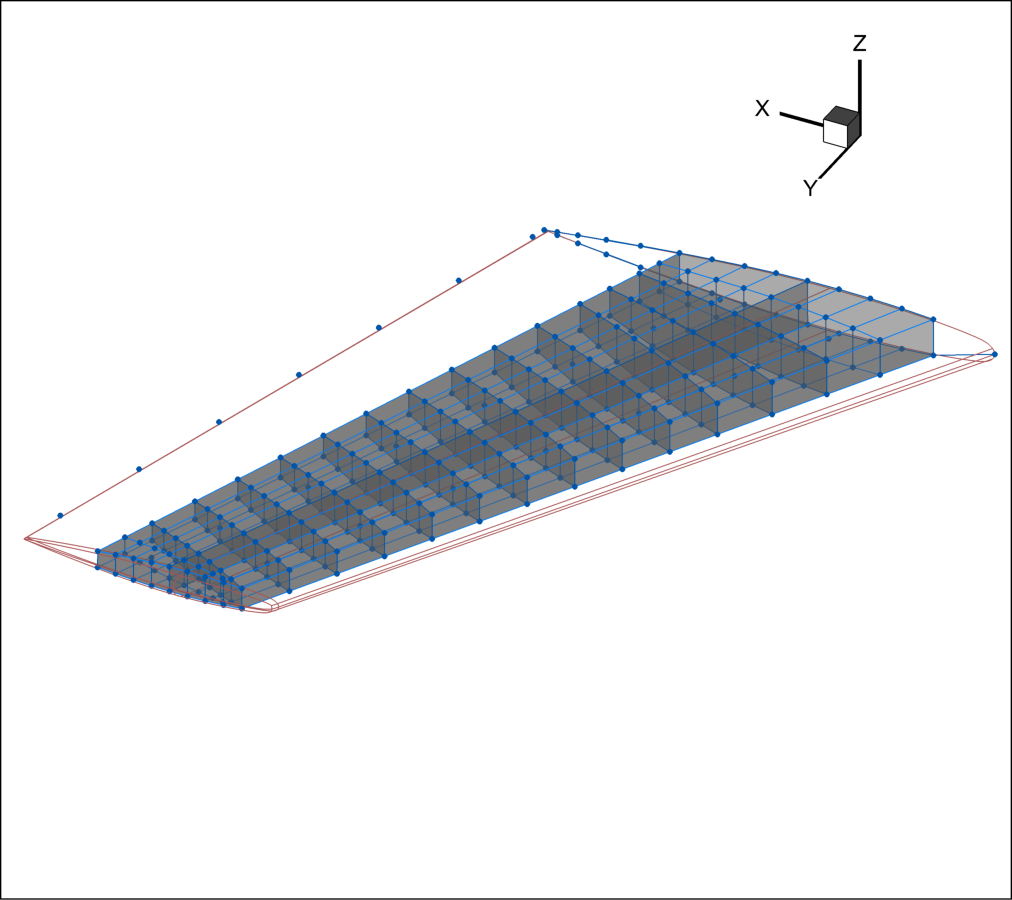}
	}
	\caption{M6 wing aeroelastic model: 42 block-structured RANS CFD mesh and FEM internal layout.}
	\label{fig:Meshes}
\end{figure}

\begin{figure}[H]
	\centering
	\parbox[h][6.5cm][c]{.5\textwidth}{
		\includegraphics[trim={1.5cm 0.5cm 1.5cm 0.5cm},clip,width=0.49\textwidth]{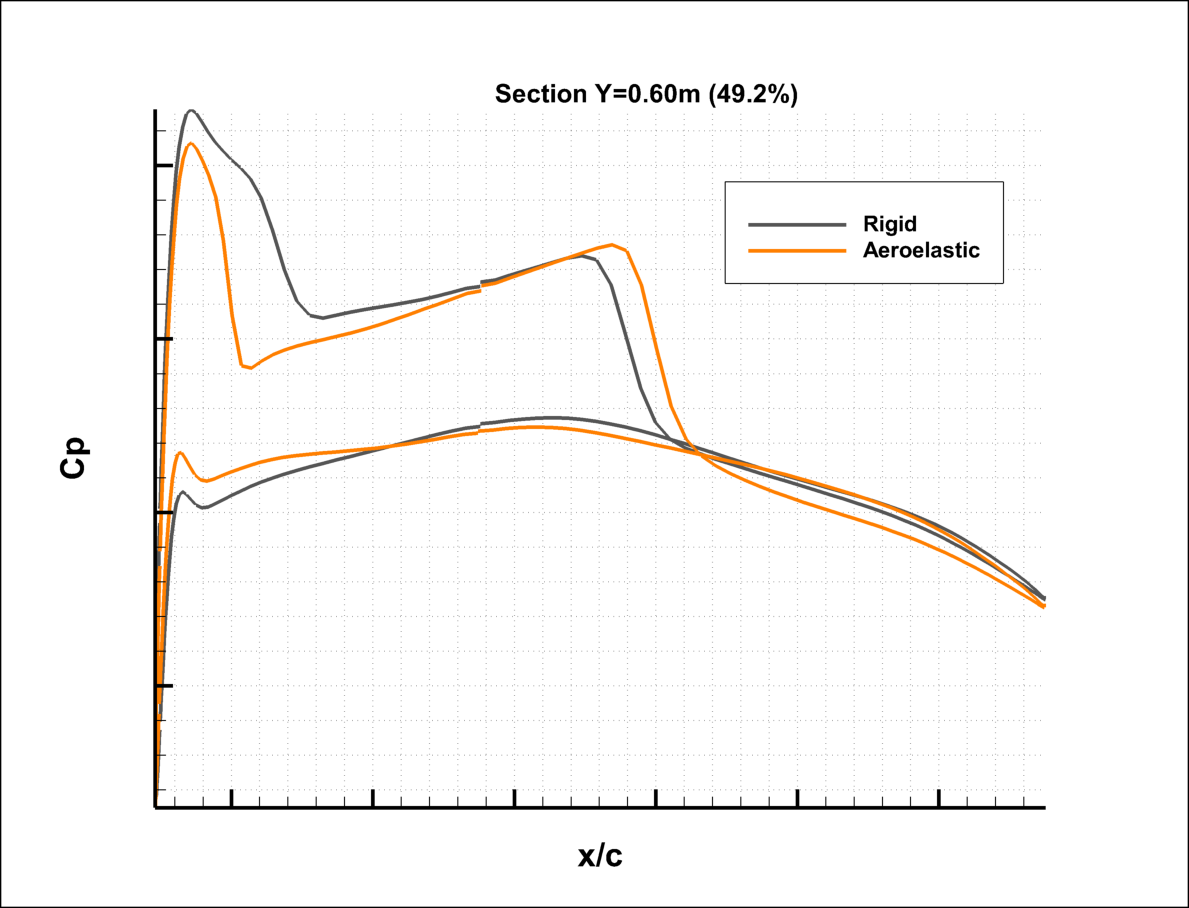}
	}\hfill
	\parbox[h][6.5cm][c]{.5\textwidth}{
		\includegraphics[trim={1.5cm 0.5cm 1.5cm 0.5cm},clip,width=0.49\textwidth]{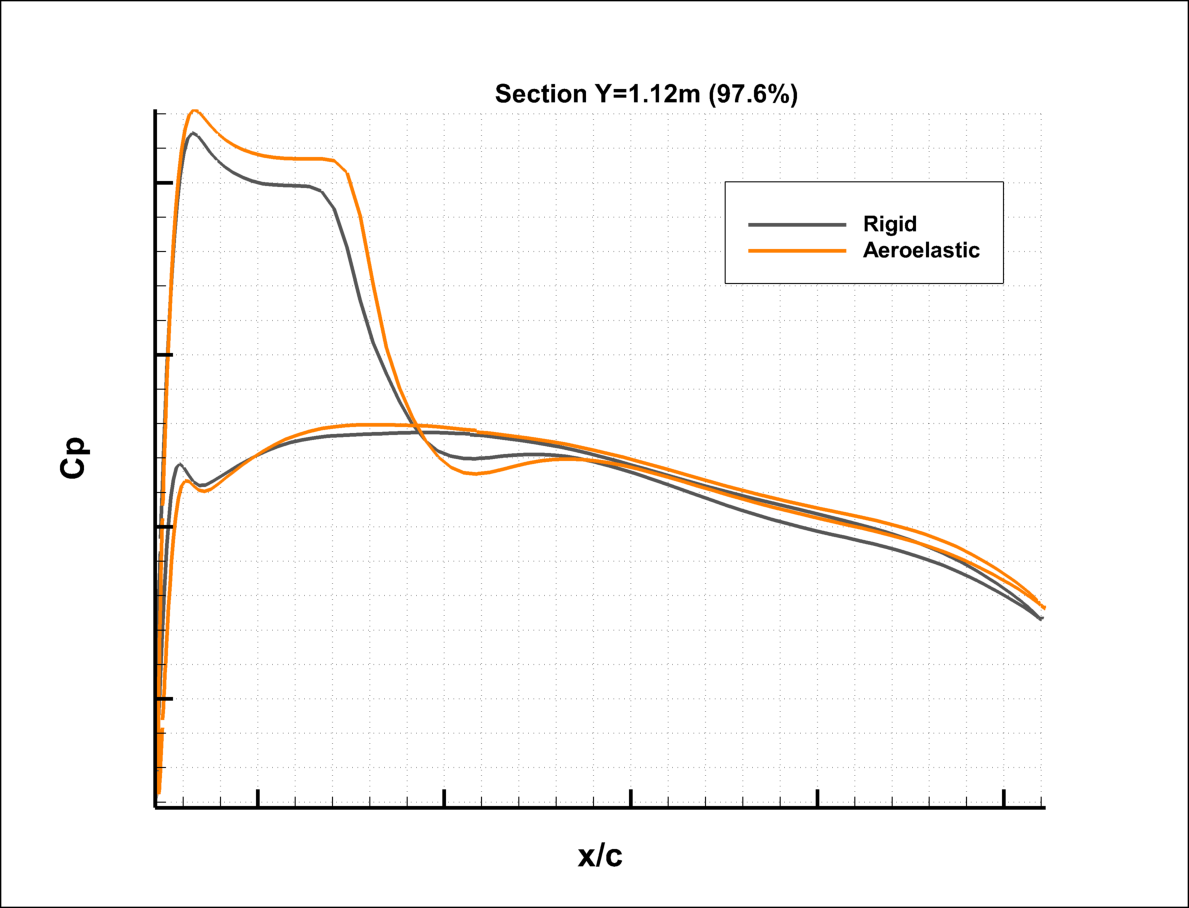}
	}
	\caption{Comparison of rigid and aeroelastic pressure coefficient section plots at y=0.60m and y=1.12m.}
	\label{fig:Cp_sections}
\end{figure}

\begin{figure}[H]
	\centering
	\includegraphics[trim={0.cm 0.cm 0.cm 0.cm},clip,width=0.60\textwidth]{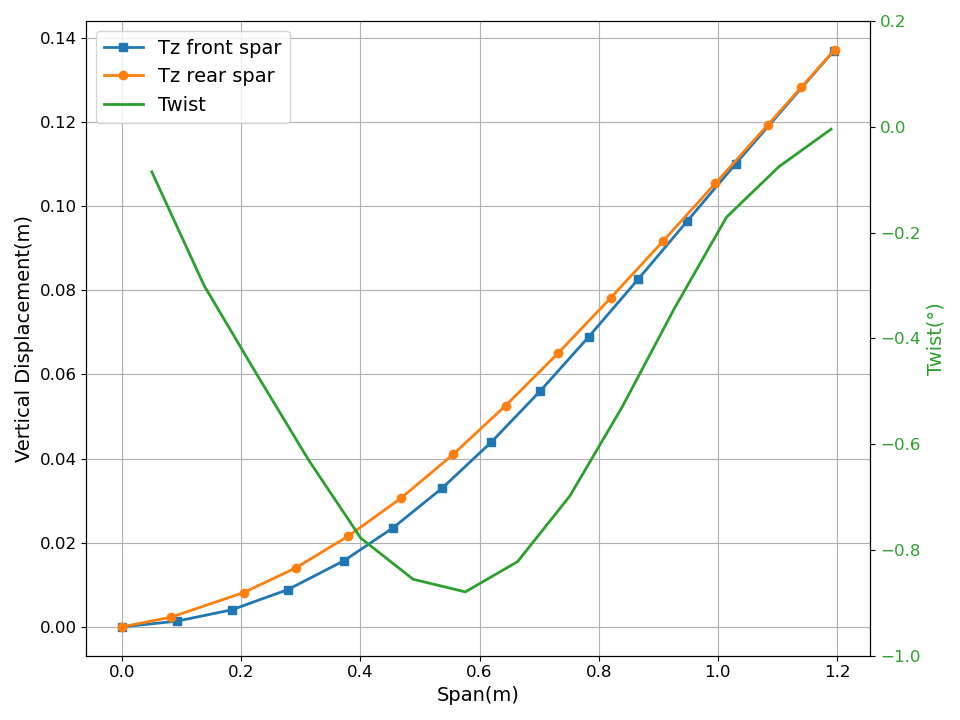}
	\caption{Vertical displacement and twist increment distribution at aeroelastic equilibrium.}
	\label{fig:DispTwistAEL}
\end{figure}

In this study the numerical experiments have been performed with the well known ONERA-M6 fixed wing configuration which has been extensively used for CFD solvers validation in transonic flow conditions. In this work we use the RANS solver provided by elsA for the steady rigid and aeroelastic analyses \citep{cambier:hal-01293795, cambier_elsA_2014}. The elsA adjoint and coupled-adjoint solvers have also been used for the computation of rigid and flexible derivatives. The latest improvements to the Krylov solvers for the solution of the adjoint linear system are described in \citep{jadoui2022comparative}.
A multi-block structured mesh featuring a C-H topology is used (Fig.~\ref{fig:Meshes}). It consists of 3.8 million grids divided into 42 blocks.
The flight conditions are a free-stream Mach number of 0.84 at an angle of attack of 3.06 degrees. The convective fluxes are discretized by an upwind Roe scheme associated to a Monotonic Upstream-centered Scheme for Conservation Laws (MUSCL) reconstruction and a Van Albada flux limiter. The one-equation Spalart-Allmaras turbulence model is selected. The surface contours in the bottom plot of Fig.~\ref{fig:Cp_contours} below show typical results for the ONERA-M6 wing. The pressure coefficient contours identify a lambda-shock along the mid-chord of the wing. For the aeroelastic analysis a simple but realistic finite element model has been designed (Fig.~\ref{fig:Meshes}). The stiffness of this model can be easily tuned to get stronger or weaker fluid-structure interaction. The pressure coefficient contours at the aeroelastic equilibrium are plotted in the upper part of Fig.~\ref{fig:Cp_contours} and can be compared to the rigid contours. The maximum vertical displacement is 0.14 meters corresponding to 11.7 \% of the wing span. To get a better insight of the effect of flexibility on the pressure distribution, we report in Fig.~\ref{fig:Cp_sections} the Cp distributions for two sections at $y=0.60$ m and $y=1.12$ m. The vertical displacement distributions associated to the front and rear spars as well as the twist increment distribution are plotted in Fig.~\ref{fig:DispTwistAEL}. The rigid analysis results in a lift coefficient $C_L = 0.27$ whereas the aeroelastic analysis, at the same angle of attack, results in a lower lift coefficient $C_L = 0.23$.
\vspace{0.5cm}

\begin{figure}[H]
	\centering
	\includegraphics[trim={4.5cm 19.7cm 0.5cm 1.0cm},clip,width=0.70\textwidth]{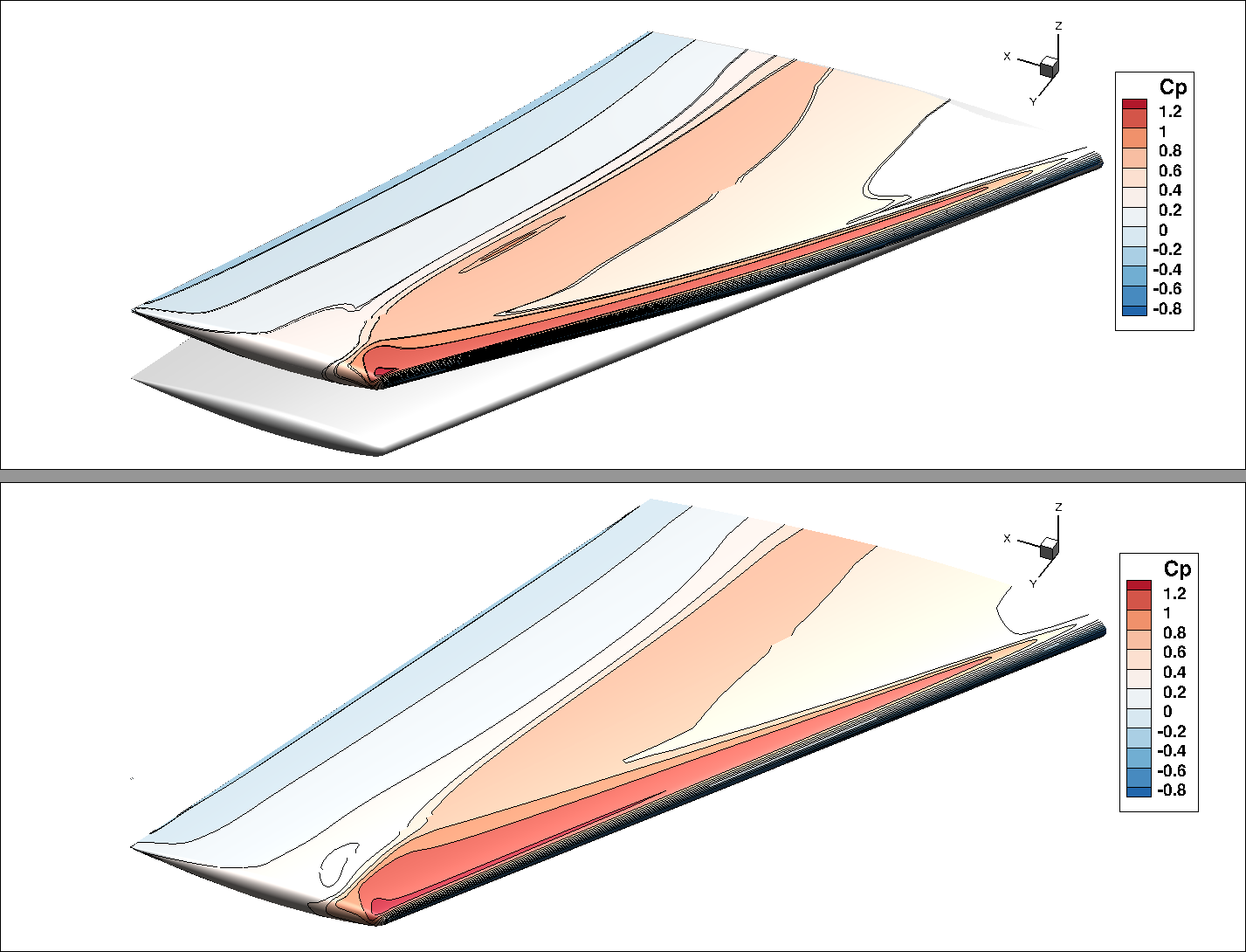}
	\includegraphics[trim={4.5cm 0.1cm 0.5cm 19.7cm},clip,width=0.70\textwidth]{M6_RANS_CP.png}
	\caption{Pressure coefficient contour plots for the rigid and aeroelastic steady flows.}
	\label{fig:Cp_contours}
\end{figure}

\section{Minimal residual Krylov subspace methods combined with spectral deflation} \label{FGMRESDR_Review}

In this section we focus on a particular minimal residual norm Krylov subspace method for the solution of linear systems with a non-symmetric real coefficient matrix of type
\begin{equation}
	\label{eq:Linear_System}
	\mathbf{Ax} = \mathbf{b}, \qquad \mathbf{A} \in \mathbb{R}^{N \times N}; \quad \mathbf{b}, \hspace{0.1cm} \mathbf{x} \hspace{0.1cm} \in \mathbb{R}^N,
\end{equation}

\noindent with the initial guess $\mathbf{x}_0$ and the associated residual $\mathbf{r}_0 = \mathbf{b} - \mathbf{Ax}_0$. The GMRES method \cite{saad1986} computes the correction $\mathbf{z}_i$ in the $i$th Krylov subspace $\mathcal{K}_i(\mathbf{A},\mathbf{r}_0) \equiv \operatorname{span}\{\mathbf{r}_0,\mathbf{Ar}_0,\mathbf{A}^2\mathbf{r}_0,\cdots,\mathbf{A}^{i-1}\mathbf{r}_0\}$ that minimizes the norm of the residual $\mathbf{r}_i = \mathbf{b} - \mathbf{A}(\mathbf{x}_0+\mathbf{z}_i) = \mathbf{r}_0 - \mathbf{Az}_i$. The relation between the minimal residual correction $\mathbf{z}_i$ and the orthogonality of the new residual $\mathbf{r}_i$ to the \emph{shifted Krylov space} $\mathcal{AK}_i(\mathbf{A},\mathbf{r}_0) \equiv \operatorname{span}\{\mathbf{Ar}_0,\mathbf{A}^2\mathbf{r}_0,\cdots,\mathbf{A}^{i}\mathbf{r}_0\}$ is given by the following theorem \cite{saad2003}:
\vspace{6pt}
\begin{thm}
    The vector $\mathbf{z}_i \in \mathcal{K}_i(\mathbf{A},\mathbf{r}_0)$ satisfies $\mathbf{z}_i = \underset{\mathbf{z} \in \mathcal{K}_i(\mathbf{A},\mathbf{r}_0)}{\operatorname{argmin}} \| \mathbf{r}_0 - \mathbf{Az} \|_2 \Leftrightarrow \mathbf{r}_i \; \bot \; \mathcal{AK}_i(\mathbf{A},\mathbf{r}_0)$
\end{thm}
\vspace{6pt}
This is known as the optimality property of the residual. In this work, a right-preconditioned system is considered so that (\ref{eq:Linear_System}) becomes
\begin{align}
	\label{eq:Variable_Precond_System}
	\mathbf{A}\mathcal{M}(\mathbf{t}) = \mathbf{b}, \\
	\mathbf{x} = \mathcal{M}(\mathbf{t})
\end{align}
with $\mathbf{t} \in \mathbb{R}^{N}$ and $\mathcal{M}:\mathbb{R}^{N} \rightarrow \mathbb{R}^{N}$ the preconditioning operator which may be a nonlinear function.

\subsection{Flexible GMRES algorithm with right preconditioning}

Saad proposed a minimal residual norm subspace method based on the standard GMRES approach \cite{saad1986} that allows a variable nonlinear preconditioning function $\mathcal{M}_{j}:\mathbb{R}^{N} \rightarrow \mathbb{R}^{N}$ at each iteration $j$ \cite{saad1993flexible}.
Starting from an initial guess $\mathbf{x}_{0}$, the flexible Arnoldi relation is written as:
\begin{equation}
	\label{eq:Flexible_Arnoldi}
	\mathbf{AZ}_{m} = \mathbf{V}_{m+1}\Bar{\mathbf{H}}_{m},
\end{equation}

\noindent where the matrices $\mathbf{V}_{m+1} \in \mathbb{R}^{N \times (m+1)}, \mathbf{Z}_{m} \in \mathbb{R}^{N \times m}$ and $\bar{\mathbf{H}}_{m} \in \mathbb{R}^{(m+1) \times m}$ stand for the orthonormal basis of the Krylov space, the solution space and the upper Hessenberg matrix respectively. The approximate solution is written as $\mathbf{x}_{m} = \mathbf{x}_{0} + \mathbf{Z}_{m}\mathbf{y}_{m}$ where $\mathbf{y}_{m}$ minimizes $|| \mathbf{r}_{0} -\mathbf{AZ}_{m}\mathbf{y} ||_{2}$ over $\mathbf{x}_{0} + \operatorname{span}\lbrace \mathbf{Z}_{m} \rbrace$, with $\mathbf{Z}_{m}=\mathbfcal{M}(\mathbf{V}_{m})=\left[\mathcal{M}_{1}(\mathbf{v}_{1}),\cdots,\mathcal{M}_{m}(\mathbf{v}_{m}) \right]$ and both $\mathbf{Z}_{m}$ and $\mathbf{V}_{m}$ need to be stored. We point out that the operator $\mathbfcal{M}$ represents the action of the nonlinear operators $\mathcal{M}_{j}$ on the set of basis vectors $\mathbf{v}_{j}$.

The restarted FGMRES($m$, $m_{i}$) pseudocode is presented in Algorithm \ref{alg:FGMRES}. We denote by $m_{i}$ the size of the Krylov subspace associated to the GMRES solver devoted to the inner linear system. We point out that the stopping criterion is essentially based on the least-squares relative residual $||\mathbf{c} - \bar{\mathbf{H}}_{m}\mathbf{y}_{m}||/||\mathbf{b}||$, with $\mathbf{c}=\|\mathbf{r}_0\|\mathbf{e}_1$ (see \cite{fraysse2008algorithm} for several definitions of stopping criteria and some practical considerations for the implementation of the key points of the algorithm). The latter is a cheap approximation of the true residual $\| \mathbf{Ax}-\mathbf{b} \|/\|\mathbf{b}\|$, assuming $\mathbf{V}_{m+1}$ is orthonormal.

In an inner-outer FGMRES the preconditioning operation $\mathbf{z}_{j} = \mathcal{M}_{j}(\mathbf{v}_{j})$ of step \ref{alg:innerGMRES} can be thought of as a means of approximately solving $\mathbf{Az}_j=\mathbf{v}_j$ where $\mathbf{M}^{-1}_j \approx \mathbf{A}^{-1}$ is the inner preconditioner. 
To prevent unnecessary propagation of rounding errors, the relative true residual is only computed at the end of each cycle and is used to construct the first vector of the next Krylov subspace basis.
In the case of large and ill-conditioned linear systems, least-squares and true residuals may differ due to loss of orthogonality during the construction of the Krylov basis.
A standard way to tackle such a phenomenon is to ask for a second iteration of the Modified Gram-Schmidt algorithm (loop from line 6 to 9) in order to strengthen the orthogonality of the Krylov basis.

In this work, we consider the standard GMRES with stationary right-preconditioning on one hand, and the specific class of flexible nested GMRES strategy on the other hand where an inner GMRES acts as an iterative preconditioner. Two stationary preconditioning strategies are considered. The first one consists in a block version of a standard LU-SGS iterative solver. LU-SGS is applied to a first order diagonally dominant upwind approximation of the flux Jacobian matrix inspired by~\cite{yoon1988}. This operator is based on a first order spatial discretization of the convective and of the viscous fluxes using a simplifying thin layer assumption~\cite{peter2007large}. This strategy leads to a very compact stencil for the preconditioning matrix which will be denoted by $\textbf{J}^{APP}_{O1}$ in the sequel of this document. 
\newline
\newline
\begin{algorithm}[H]
    \fontsize{10}{14}\selectfont    
    \label{alg:FGMRES}
    \caption{Right-preconditioned FGMRES($m$,$m_{i}$)}
    {
        Choose an initial guess $\mathbf{x}_{0}$, a convergence threshold $\epsilon$ and krylov size $m$\;
        Compute $\mathbf{r}_{0} = \mathbf{b}-\mathbf{Ax}_{0}$, $\beta =||\mathbf{r}_{0}||$, $\mathbf{c}=\beta\mathbf{e}_1$ and  $\mathbf{v}_{1} = \mathbf{r}_{0}/\beta$ \;
        \For{$j \leftarrow 1$ \KwTo $m$} {
            $\mathbf{z}_{j} = \mathcal{M}_{j}(\mathbf{v}_{j})$ \label{alg:innerGMRES}  \Comment*[f]{Inner iteration solved with GMRES($m_i$)} \; 
            $\mathbf{w} = \mathbf{Az}_{j}$ \label{alg:matvecprod} \;
            \For{$i \leftarrow 1$ \KwTo $j$} {
                $h_{i,j} = \mathbf{v}_{i}^{T}\mathbf{w}$ \label{alg:scalprod} \;
                $\mathbf{w} = \mathbf{w} -h_{i,j}\mathbf{v}_{i}$ \;
            }
            $h_{j+1,j} =||\mathbf{w}||_2$ and $\mathbf{v}_{j+1} = \mathbf{w}/h_{j+1,j}$ \;
            Solve the least-squares problem $min_{\mathbf{y}} ||\mathbf{c} - \bar{\mathbf{H}}_{j}\mathbf{y}||_2$ for $\mathbf{y}_{j}$ \;
            Exit if $||\mathbf{c} - \bar{\mathbf{H}}_{j}\mathbf{y}_{j}||/||\mathbf{b}|| \le \epsilon $ \;
        }
        Compute $\mathbf{x}_{m} = \mathbf{x}_{0} + \mathbf{Z}_{m}\mathbf{y}_{m}$ where $\mathbf{Z}_{m} = [\mathbf{z}_{1},...,\mathbf{z}_{m}]$ \;
        Set $\mathbf{x}_{0} = \mathbf{x}_{m}$ and go to 2 \;
    }
\end{algorithm}  \hfill \break

The second one is a Block Incomplete LU (BILU($k$)) factorization applied to either an approximate or exact first-order flux Jacobian matrix. For the so-called first-order exact Jacobian matrix $\textbf{J}^{EXA}_{O1}$ a first-order spatial Roe scheme is used for the discretization of the mean-flow convective fluxes and a 5-point corrected centered discretization scheme is used for the diffusive fluxes.  More specifically, the spatial gradients at the cell interfaces are modified to avoid high frequency oscillations (see~\cite{resmini2015} or~\cite{puigt2014}). The BILU($k$) preconditioner will be applied either to the first-order approximate Jacobian matrix $\textbf{J}^{APP}_{O1}$, or to the first-order exact Jacobian matrix $\textbf{J}^{EXA}_{O1}$. $\textbf{J}^{EXA}_{O1}$ is different from $\textbf{J}^{APP}_{O1}$ when it comes to memory footprint. Indeed, $\textbf{J}^{EXA}_{O1}$ has a 9-point stencil in 2D whereas a 5-point stencil is associated with $\textbf{J}^{APP}_{O1}$. In 3D, we have a 7-point stencil for $\textbf{J}^{APP}_{O1}$ and a stencil of 19 points for $\textbf{J}^{EXA}_{O1}$. Consequently, a better robustness is achieved but at the price of about twice the storage for $\textbf{J}^{EXA}_{O1}$ compared to $\textbf{J}^{APP}_{O1}$.
For a better understanding, we reproduce the 3D stencils in Fig.~\ref{fig:stencils} below:
\vspace{6pt}
\begin{figure}[H]
    \centering
    \begin{minipage}[b]{0.42\textwidth}
        \includegraphics[trim=-0.5cm -0.2cm 0.cm 0.cm, clip, width=0.9\linewidth, left]{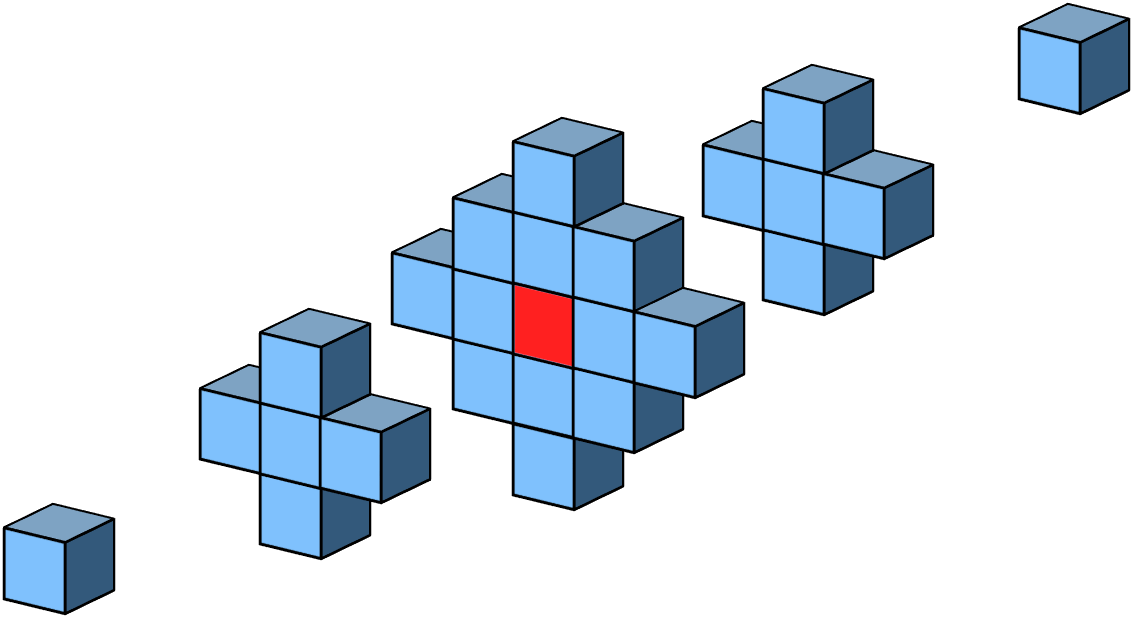}
        \centering
        \small{a)~Exact $2^{nd}$ order (size~25)}
        \label{fig:stencil_JO2}
    \end{minipage}\hfill
    \begin{minipage}[b]{0.28\textwidth}
        \includegraphics[trim=0.cm -1.2cm 0.cm 0.cm, clip, width=0.85\linewidth, center]{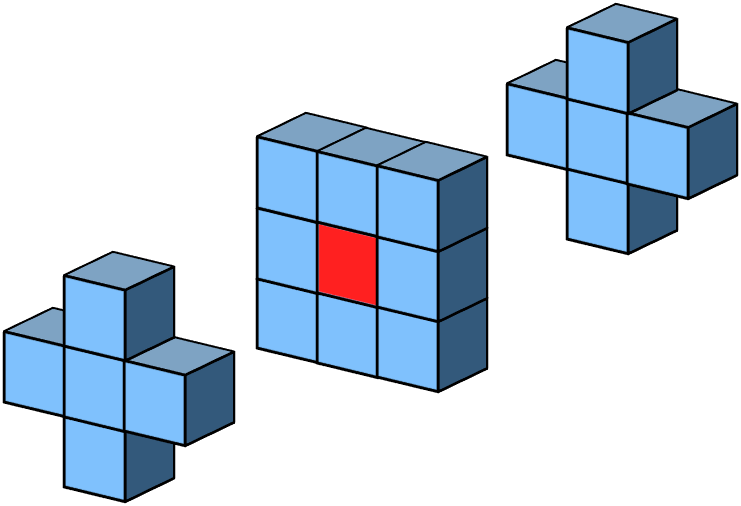}
        \centering
        \small{b)~Exact $1^{st}$ order (size~19)}
        \label{fig:Stencil_JO1_exa}
    \end{minipage}\hfill
    \begin{minipage}[b]{0.30\textwidth}
        \includegraphics[trim=0.cm -2.5cm -1.0cm 0cm, clip, width=0.70\linewidth, right]{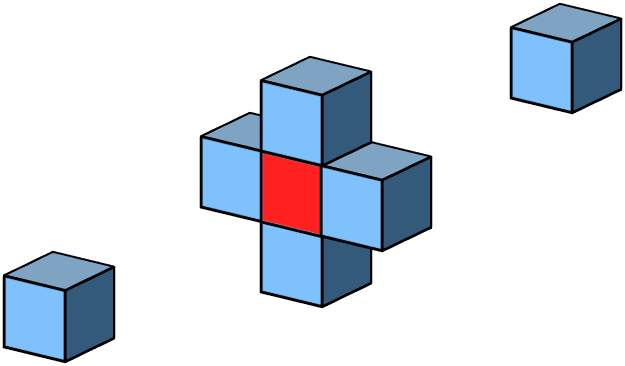}
        \centering
        \small{c)~Approx. $1^{st}$ order (size~7)}
        \label{fig:Stencil_JO1_app}
    \end{minipage}
	\vspace{0pt}
    \caption{Stencils of the various Jacobian matrices. Stencil (b) is used for BILU(k) preconditioners and stencil (c) for BILU(k) or iterative LU-SGS preconditioners.}
    \label{fig:stencils}
\end{figure}
We point out that the relevant numerical ingredients that characterize the GMRES algorithm are the matrix product (step \ref{alg:matvecprod} in algorithm \ref{alg:FGMRES}), the preconditioning strategy (step \ref{alg:innerGMRES}) and the scalar product (step \ref{alg:scalprod}). These algebraic operations are global in conjunction with a domain decomposition method. More specifically, the globalization of the preconditioner (step \ref{alg:innerGMRES}) is achieved with a Restricted Additive Schwarz method~\cite{cai1999}. In addition, the product by the operator $\mathbf{A}$ (step \ref{alg:matvecprod}) is exact. We thus get a global and parallel FGMRES($m$,$m_{i}$).

\subsection{Deflated restarting} \label{DeflatedRestarting_FGMRESDR}

The main drawback of the restarted GMRES($m$) and FGMRES($m$, $m_i$) is the loss of spectral information contained in the current Krylov subspace during the restarting procedure.
Let us recall the definition of a Ritz pair \cite{morgan1998} as it plays an important role in the strategy of deflated restarting.
\vspace{6pt}
\begin{definition}[Ritz pair] \label{Ritzpair}
    Consider a subspace $\mathcal{U}$ of $\mathbb{R}^N$. Given an operator $\mathbf{B} \in \mathbb{R}^{N \times N}$, a scalar $\lambda \in \mathbb{R}$ and $\mathbf{y} \in \mathbb{R}^N$, ($\lambda$, $\mathbf{y}$) is a Ritz pair of  $\mathbf{B}$ with respect to $\mathcal{U}$ if and only if the residual of the eigenvalue problem $\mathbf{By} = \lambda \mathbf{y}$ satisfies the following orthogonality condition:
    \begin{equation}
        \label{eq:Eigenvalue_equation}
        (\mathbf{By} - \lambda \mathbf{y}) \perp \mathcal{U} \quad \forall \mathbf{y} \in \mathcal{U}.
    \end{equation}    
\end{definition}

Choosing $\mathcal{U} \equiv \mathcal{K}_{m}(\mathbf{A},\mathbf{r}_{0})$ and $\mathbf{B} \equiv \mathbf{A}$, we have
\begin{equation}
    (\mathbf{Ay} -\lambda \mathbf{y}) \perp \mathcal{K}_{m}(\mathbf{A},\mathbf{r}_{0}) \quad \forall \mathbf{y} \in \mathcal{K}_{m}(\mathbf{A},\mathbf{r}_{0}).
    \label{eq:Ritz_problem}
\end{equation}

Recalling $\mathbf{V}_{m}$ is an orthonormal basis of $\mathcal{K}_{m}(\mathbf{A},\mathbf{r}_{0})$, we can write $\mathbf{y}=\mathbf{V}_{m}\mathbf{g}$, $\mathbf{g} \in \mathbb{R}^m$, and using the standard Arnoldi relation $\mathbf{AV}_{m} = \mathbf{V}_{m+1}\bar{\mathbf{H}}_{m}$ in (\ref{eq:Ritz_problem}), we obtain the standard eigenvalue problem
\begin{equation}
    \mathbf{H}_{m}\mathbf{g} = \lambda \mathbf{g},
    \label{eq:Ritz_spectral}
\end{equation}

\noindent where $\mathbf{H}_{m} = \left[ \mathbf{I}_m \ \mathbf{0}_{m \times 1} \right] \bar{\mathbf{H}}_{m}$. Thus, the spectral residual norm of the Ritz pair $\left\{\lambda,\mathbf{y}=\mathbf{V}_{m}\mathbf{g}\right\}$ satisfies:
\begin{equation}
    \label{eq:Residual_eigenvalue_norm}
    \| \mathbf{A}(\mathbf{V}_{m}\mathbf{g}) - \lambda (\mathbf{V}_{m}\mathbf{g}) \| = \| \mathbf{A}\mathbf{V}_{m}\mathbf{g}-\mathbf{V}_{m}\mathbf{H}_{m}\mathbf{g}\| = \| \mathbf{V}_{m+1}\Bar{\mathbf{H}}_{m}\mathbf{g} - \mathbf{V}_{m}\mathbf{H}_{m}\mathbf{g} \| = | h_{m+1,m} | |\mathbf{e}_{m}^{T}\mathbf{g}|.
\end{equation}

From a small value of $|h_{m+1,m}| |\mathbf{e}_{m}^{T}\mathbf{g}|$, the Arnoldi method takes the Ritz pair as a good approximation of the eigenvalue-eigenvector pair of the operator $\mathbf{A}$ \cite{arnoldi1951,saad1980}.
Indeed, neglecting the last row of the rectangular upper Hessenberg matrix $\Bar{\mathbf{H}}_{m}$ leads to:
\begin{equation}
    \label{eq:Hessenberg_Restriction_operator}
    \mathbf{H}_{m} = \mathbf{V}_{m}^{T}\mathbf{AV}_{m}.
\end{equation}

Therefore, the spectrum of $\mathbf{H}_{m}$ naturally approximates a part of the spectrum of $\mathbf{A}$.
The idea of deflation techniques is to keep relevant spectral information from $\mathbf{H}_{m}$ in the search space of the next cycle to expect a better convergence of the Krylov iterative methods.
In~\cite{giraud2010}, Giraud et al. take into account both smallest and largest eigenvalues to maximize the deflation effect.
In contrast, Morgan~\cite{morgan2002gmres} only deflates the smallest ones. This last strategy will be adopted for our numerical experiments.

Actually, Ritz values of the operator $\mathbf{A}$ give a good approximation of its exterior eigenvalues, i.e., of largest magnitude.
Unfortunately, interior eigenvalues are of greater interest because they are generally responsible for the convergence stagnation. An alternative which does better at finding eigenvalues nearest zero was proposed by Morgan \cite{morgan1991computing, morgan1995restarted}
who introduced the harmonic Ritz values of $\mathbf{A}$ with respect to $\mathcal{K}_{m}(\mathbf{A},\mathbf{r}_{0})$, which are equivalent to the Ritz values of $\mathbf{A}^{-1}$ with respect to $\mathcal{AK}_{m}(\mathbf{A},\mathbf{r}_{0})$. In practice we will use the following definition  

\vspace{6pt}
\begin{definition}[harmonic Ritz pair] \label{HarmonicRitzpair}
    Consider a subspace $\mathcal{U}$ of $\mathbb{R}^N$. Given an operator $\mathbf{B} \in \mathbb{R}^{N \times N}$, a scalar $\lambda \in \mathbb{R}$ and $\mathbf{y} \in \mathbb{R}^N$, ($\lambda$, $\mathbf{y}$) is a harmonic Ritz pair of  $\mathbf{B}$ with respect to $\mathcal{U}$ if and only if the residual of the eigenvalue problem $\mathbf{By} = \lambda \mathbf{y}$ satisfies the following orthogonality condition:
    \begin{equation}
        \label{eq:Eigenvalue_equation}
        (\mathbf{By} - \lambda \mathbf{y}) \perp \mathbf{B}\mathcal{U} \quad \forall \mathbf{y} \in \mathcal{U}.
    \end{equation}    
\end{definition}

Choosing $\mathcal{U} \equiv \mathcal{K}_{m}(\mathbf{A},\mathbf{r}_{0})$ and $\mathbf{B} \equiv \mathbf{A}$, we have,

\begin{equation}
    (\mathbf{Ay} -\lambda \mathbf{y}) \perp \mathcal{AK}_{m}(\mathbf{A},\mathbf{r}_{0}) \quad \forall \mathbf{y} \in \mathcal{K}_{m}(\mathbf{A},\mathbf{r}_{0}).
    \label{eq:Harmonic_Ritz_problem}
\end{equation}

Recalling $\mathbf{V}_{m}$ is an orthonormal basis of $\mathcal{K}_{m}(\mathbf{A},\mathbf{r}_{0})$, we can write $\mathbf{y}=\mathbf{V}_{m}\mathbf{g}$, $\mathbf{g} \in \mathbb{R}^m$, and using the flexible Arnoldi relation \eqref{eq:Flexible_Arnoldi} in (\ref{eq:Harmonic_Ritz_problem}), we obtain the following generalized eigenvalue problem
\begin{align}
    & \label{eq:Generalized_eigenvalue_problem_a}
	(\mathbf{AZ}_{m})^T (\mathbf{AZ}_{m} \mathbf{g} - \lambda \mathbf{V}_{m} \mathbf{g}) = \mathbf{0}  \\    
    \Leftrightarrow \quad & \Bar{\mathbf{H}}_{m}^{T} \Bar{\mathbf{H}}_{m}\mathbf{g} = \lambda  \bar{\mathbf{H}}^T_{m} \mathbf{V}^T_{m+1} \mathbf{V}_m \mathbf{g} \nonumber \\
	\label{eq:Generalized_eigenvalue_problem_b}
    \Leftrightarrow \quad &  \boxed{\Bar{\mathbf{H}}_{m}^{T} \Bar{\mathbf{H}}_{m}\mathbf{g} = \lambda \mathbf{H}_{m}^{T} \mathbf{g}}
\end{align}

After some algebraic manipulations, (\ref{eq:Generalized_eigenvalue_problem_b}) can be reformulated as a standard eigenvalue problem (see \cite{morgan1998} for this formula):
\begin{equation}
    \label{eq:Standard_eigenvalue_problem}
    (\mathbf{H}_{m} + h_{m+1,m}^{2}\mathbf{H}_{m}^{-T}\mathbf{e}_{m}\mathbf{e}_{m}^{T})\mathbf{g} = \lambda \mathbf{g},
\end{equation}
where $\lambda$ is a harmonic Ritz value and the corresponding harmonic Ritz vector is $\mathbf{y}=\mathbf{V}_{m}\mathbf{g}$. To follow definition \ref{HarmonicRitzpair}, we use the identity $\mathbf{V}_m^T \mathbf{V}_m = \mathbf{I}$ into equation \eqref{eq:Generalized_eigenvalue_problem_a} which gives
\begin{equation}
	((\mathbf{AZ}_{m}\mathbf{V}_m^T) \mathbf{V}_m)^T ((\mathbf{AZ}_{m}\mathbf{V}_m^T) \mathbf{V}_m \mathbf{g} - \lambda \mathbf{V}_{m} \mathbf{g}) = \mathbf{0}.
\end{equation}

Thus, $\mathbf{Y}_m = \{ \mathbf{y}_1, \cdots, \mathbf{y}_m \}$ corresponds to harmonic Ritz vectors of $\mathbf{AZ}_{m}\mathbf{V}_m^T$ with respect to $\operatorname{range}(\mathbf{V}_m)$. To underline the link with FGCRO-DR in section \ref{FGCRODR_theory}, we denote by \emph{strategy B} this deflation strategy.

Obviously, in exact arithmetic solutions to (\ref{eq:Generalized_eigenvalue_problem_b}) and (\ref{eq:Standard_eigenvalue_problem}) are identical. But it is not the case in finite precision since the operator $\Bar{\mathbf{H}}_{m}^{T} \Bar{\mathbf{H}}_{m}$ is usually ill-conditioned in the fully linearized turbulence case. Therefore, the accurate estimation of the eigenvectors could be strongly altered and lead to stagnation of the relative true residual of the GMRES process.

As mentioned by Giraud et al. \cite{giraud2010}, the flexible Arnoldi relation obtained at each restart within the FGMRES with deflated restarting (FGMRES-DR) framework given by
\begin{equation}
    \label{eq:Flexible_Arnoldi_Realtion_Deflation}
    \mathbf{AZ}_{k} = \mathbf{V}_{k+1}\bar{\mathbf{H}}_{k},
\end{equation}

\noindent holds with $\mathbf{Z}_{k} = \mathbf{Z}_{m}\mathbf{P}_{k}$, $\mathbf{V}_{k+1} = \mathbf{V}_{m+1}\mathbf{P}_{k+1}$ and $\bar{\mathbf{H}}_{k} = \mathbf{P}_{k+1}\bar{\mathbf{H}}_{m}\mathbf{P}_{k}$ where $\mathbf{P}_{k} \in \mathbb{R}^{m \times k}$ corresponds to the orthonormal matrix whose columns are spanned by the set of the $k$ retained eigenvectors of (\ref{eq:Standard_eigenvalue_problem}).

Also, the right-hand side of the least-squares problem is computed at each restart as $\mathbf{c} = \mathbf{V}_{m+1}^{T}\mathbf{r}_{m}$ which requires $(2N-1)(m+1)$ operations. Rollins and Fichtner \citep{rollin2008} have proposed an efficient way to compute $\mathbf{c}$ so that we can save some inner products. Indeed, they notice that the residual $\mathbf{r}_{m}$ is a linear combination of the columns of the deflation subspace $\mathbf{V}_{k+1}$. Consequently, we compute $ \mathbf{c} = [ (\mathbf{r}_{m}^T \mathbf{V}_{k+1})^T \quad \mathbf{0}_{1 \times (m-k)} ]^T $ with $(2N-1)(k+1)$ operations only. Also, they improved the construction of $\mathbf{P}_{k+1} \in \mathbb{R}^{(m+1)\times(k+1)}$ and in particular the last column of this matrix which is usually chosen as the vector $\mathbf{c} - \bar{\mathbf{H}}_{m}\mathbf{y}_{m}$. More specifically, they demonstrated that the vector $\mathbf{c} - \bar{\mathbf{H}}_{m}\mathbf{y}_{m}$ is colinear to the vector $[-\delta \mathbf{f}^{T} \; 1]^{T}$ with $\mathbf{f}=\mathbf{H}^{-T}_{m}\mathbf{e}_{m}$  and $\delta=h_{m+1,m}$. The explicit relation is given below:
\begin{equation}
    \label{eq:Rollin_vector}
    \mathbf{c} - \bar{\mathbf{H}}_{m}\mathbf{y}_{m} = \begin{pmatrix} -\delta \mathbf{f} \\ 1 \end{pmatrix} \left( \frac{\omega - \delta \mathbf{f}^{T}\mathbf{v}}{1+\delta^{2}\mathbf{f}^{T}\mathbf{f}} \right),
\end{equation}

\noindent where $\mathbf{v}$ and $\omega$ are the first $m$ rows of $\mathbf{c}$ and the last element of $\mathbf{c}$ respectively. This formulation reduces the propagation of rounding errors in the Arnoldi relation at restart.
\vspace{\baselineskip}

\begin{algorithm}[H]
	\caption{FGMRES-DR($m$,$m_{i}$,$k$)}
	\fontsize{10}{14}\selectfont
	Choose $\mathbf{x}_{0}$, $m$, $m_i$, $k$, $tol$. Let $\mathbf{r}_{0} = \mathbf{b}-\mathbf{Ax}_{0}$, $\beta=\|\mathbf{r}_0\|$, $\mathbf{v}_1=\mathbf{r}_0/\beta$ and $\mathbf{c}=\beta\mathbf{e}_1 \in \mathbb{R}^{m+1}$ \;
	
	Apply a full cycle of FGMRES($m$,$m_{i}$) (algorithm \ref{alg:FGMRES}) to build $\mathbf{V}_{m+1}$, $\mathbf{Z}_{m}$, $\bar{\mathbf{H}}_{m}$. \;
	
	Compute $\mathbf{x}_m = \mathbf{x}_0 + \mathbf{V}_m \mathbf{y}_m$, $ \rho_m = \| \mathbf{c}-\bar{\mathbf{H}}_m \mathbf{y}_m \|/\|\mathbf{b}\| $, where 
	$\mathbf{y}_m = {\operatorname{argmin}}_{\mathbf{y}} \| \mathbf{c}-\bar{\mathbf{H}}_m \mathbf{y} \|$ \;	
	Set $\mathbf{x}_0 = \mathbf{x}_m$ \;	
	\While{$\rho_m$ > $tol$} {
		Solve $\mathbf{f}=\mathbf{H}^{-T}_{m}\mathbf{e}_{m}$.  \;
		
		Compute the first $k$ eigenvectors $\mathbf{p}_{i}$ of $\mathbf{H}_{m} + h^{2}_{m+1,m}\mathbf{f} \, \mathbf{e}_{m}^{T}$ associated to smallest eigenvalues $\theta_i$.  \;
		
		Form $\mathbf{P}_{k} = \left[\mathbf{p}_1,\cdots,\mathbf{p}_k \right]$ and $ \mathbf{P}_{k+1} = \left[ \begin{bmatrix} \mathbf{P}_k \\\mathbf{0}_{1 \times k} \end{bmatrix}, \begin{bmatrix} -\beta \mathbf{f} \\ 1 \end{bmatrix} \right ]$  \;
		
		Obtain the reduced QR decomposition of $\mathbf{P}_{k+1} = \bar{\mathbf{P}}_{k+1} \Gamma_{k+1}$, $\bar{\mathbf{P}}_{k+1} \in \mathbb{R}^{(m+1) \times (k+1)}$ \;
		
		Form $\mathbf{V}_{k+1}^{new} = \mathbf{V}_{m+1}\Bar{\mathbf{P}}_{k+1}$, $\mathbf{Z}^{new}_{k} = \mathbf{Z}_{m}\Bar{\mathbf{P}}_{k}$ and $\Bar{\mathbf{H}}_{k}^{new} = \Bar{\mathbf{P}}_{k+1}^{T}\Bar{\mathbf{H}}_{m}\Bar{\mathbf{P}}_{k}$. \label{alg:step1} \;
		
		Let $\Bar{\mathbf{H}}_{k} = \Bar{\mathbf{H}}_{k}^{new}$, $\mathbf{V}_{k+1} = \mathbf{V}_{k+1}^{new}$ and $\mathbf{Z}_{k} =\mathbf{Z}^{new}_{k}$. \label{alg:step2} \;
		
		Apply ($m$-$k$) steps of the Arnoldi process from $\mathbf{V}_{k+1}$, $\mathbf{Z}_{k}$, and $\Bar{\mathbf{H}}_{k}$ to build $\mathbf{V}_{m+1}$, $\mathbf{Z}_{m}$ and $\Bar{\mathbf{H}}_{m}$. \;
		
		Set $\mathbf{c} = \mathbf{V}^{T}_{k+1}\mathbf{r}_{0}$, with $\mathbf{r}_{0} = \mathbf{b} - \mathbf{Ax}_0$. \;
		
		Solve the least-squares problem $\mathbf{y}_m = {\operatorname{argmin}}_{\mathbf{y}} \| \mathbf{c}-\bar{\mathbf{H}}_m \mathbf{y} \|$ \;
		
		Compute $\mathbf{x}_{m} = \mathbf{x}_{0} + \mathbf{Z}_{m}\mathbf{y}_{m}$ and $ \rho_m = \| \mathbf{c}-\bar{\mathbf{H}}_m \mathbf{y}_m \|/\|\mathbf{b}\| $ \;
		Set $\mathbf{x}_0 = \mathbf{x}_m$ \;
	}
	\label{alg:FGMRES_DR}
\end{algorithm}  \hfill \break

\subsection{Alternative deflation strategy for FGMRES-DR} \label{StrategyA_FGMRESDR}

From the expression of the solution $\mathbf{x}_{m} = \mathbf{x}_{0} + \mathbf{Z}_{m}\mathbf{y}_{m}$, where $\mathbf{Z}_{m} = [\mathbf{z}_{1},...,\mathbf{z}_{m}]$, it seems legitimate to write the Harmonic Ritz vectors as
$\mathbf{y}=\mathbf{Z}_{m}\mathbf{g}$, $\mathbf{g} \in \mathbb{R}^m$. From \eqref{eq:Eigenvalue_equation} and \eqref{eq:Harmonic_Ritz_problem} we have 
\begin{align}
	&(\mathbf{Ay} -\lambda \mathbf{y}) \perp \operatorname{span}(\mathbf{AZ}_{m}) \quad \forall \mathbf{y} \in \operatorname{span}(\mathbf{Z}_{m}) \nonumber \\
	\Leftrightarrow \quad & (\mathbf{AZ}_{m})^T (\mathbf{AZ}_{m} \mathbf{g} - \lambda \mathbf{Z}_{m} \mathbf{g}) = \mathbf{0} \nonumber  \\
	\Leftrightarrow \quad & \boxed{\Bar{\mathbf{H}}_{m}^{T} \Bar{\mathbf{H}}_{m}\mathbf{g} = \lambda  \bar{\mathbf{H}}^T_{m} \mathbf{V}^T_{m+1} \mathbf{Z}_m \mathbf{g}}
	\label{eq:Alternative_deflation_strategy_b}
\end{align}

Thus, $\mathbf{Y}_m = \{ \mathbf{y}_1, \cdots, \mathbf{y}_m \}$ corresponds to harmonic Ritz vectors of $\mathbf{A}$ or equivalently of $\mathbf{AZ}_{m} \mathbf{Z}_m^{\dagger}$ with respect to $\operatorname{range}(\mathbf{Z}_m)$.
Noticing that the last row of $\bar{\mathbf{H}}_m$ is null except the rightmost element $h_{m+1,m}$, and defining $\mathbf{f} = \mathbf{H}_m^{-T} \mathbf{e}_m$, we can reformulate \eqref{eq:Alternative_deflation_strategy_b} as
\begin{equation}
	[ \mathbf{H}_m + h_{m+1,m}^2\mathbf{f} \, \mathbf{e}_m^T ] \mathbf{g} = \lambda [ \mathbf{I}_m \quad h_{m+1,m}\mathbf{f} ] \mathbf{V}^T_{m+1} \mathbf{Z}_m \mathbf{g}.
\end{equation}

The product $\mathbf{V}^T_{m+1} \mathbf{Z}_m$ brings an additional cost to the deflation process. A block form of this term can be found as

\begin{equation}
	\label{eq:FGMRESDR_block_deflation}
	\mathbf{V}^T_{m+1} \mathbf{Z}_m =
	\begin{bmatrix}
		\mathbf{V}^T_{k+1} \mathbf{Z}_k & \mathbf{V}^T_{k+1} \mathbf{Z}_{m-k} \\
		\mathbf{V}^T_{m-k} \mathbf{Z}_k &  \mathbf{V}^T_{m-k} \mathbf{Z}_{m-k}
	\end{bmatrix}\!.
\end{equation}

Thanks to steps \ref{alg:step1} and \ref{alg:step2} in algorithm \ref{alg:FGMRES_DR} the leading block $\mathbf{V}^T_{k+1} \mathbf{Z}_k$ can be written as

\begin{equation}
	(\mathbf{V}^T_{k+1} \mathbf{Z}_k)^{(i)} = \Bar{\mathbf{P}}_{k+1}^{T} (\mathbf{V}^T_{m+1} \mathbf{Z}_m)^{(i-1)} \Bar{\mathbf{P}}_{k},
\end{equation}
where the superscript is related to the cycle index. Thus storing $\mathbf{V}^T_{m+1} \mathbf{Z}_m$ at the end of each cycle allows us to compute at a cheap cost the $(k+1)\times k$ block $\mathbf{V}^T_{k+1} \mathbf{Z}_k$ for the next cycle. To underline the link with FGCRO-DR in section \ref{FGCRODR_theory}, we denote by \emph{strategy A} this formulation.
%

\subsection{Application of GMRES-DR and FGMRES-DR to the solution of the fluid adjoint system} %
\label{GMRESDR_Fluid_Adjoint_System}

\begin{figure}[H]
	\captionsetup[subfigure]{belowskip=-6pt}
	\begin{subfigure}{.5\textwidth}
		\raggedright
		\includegraphics[trim={0.5cm 0.5cm 1.5cm 1.9cm},clip,width=0.99\textwidth]{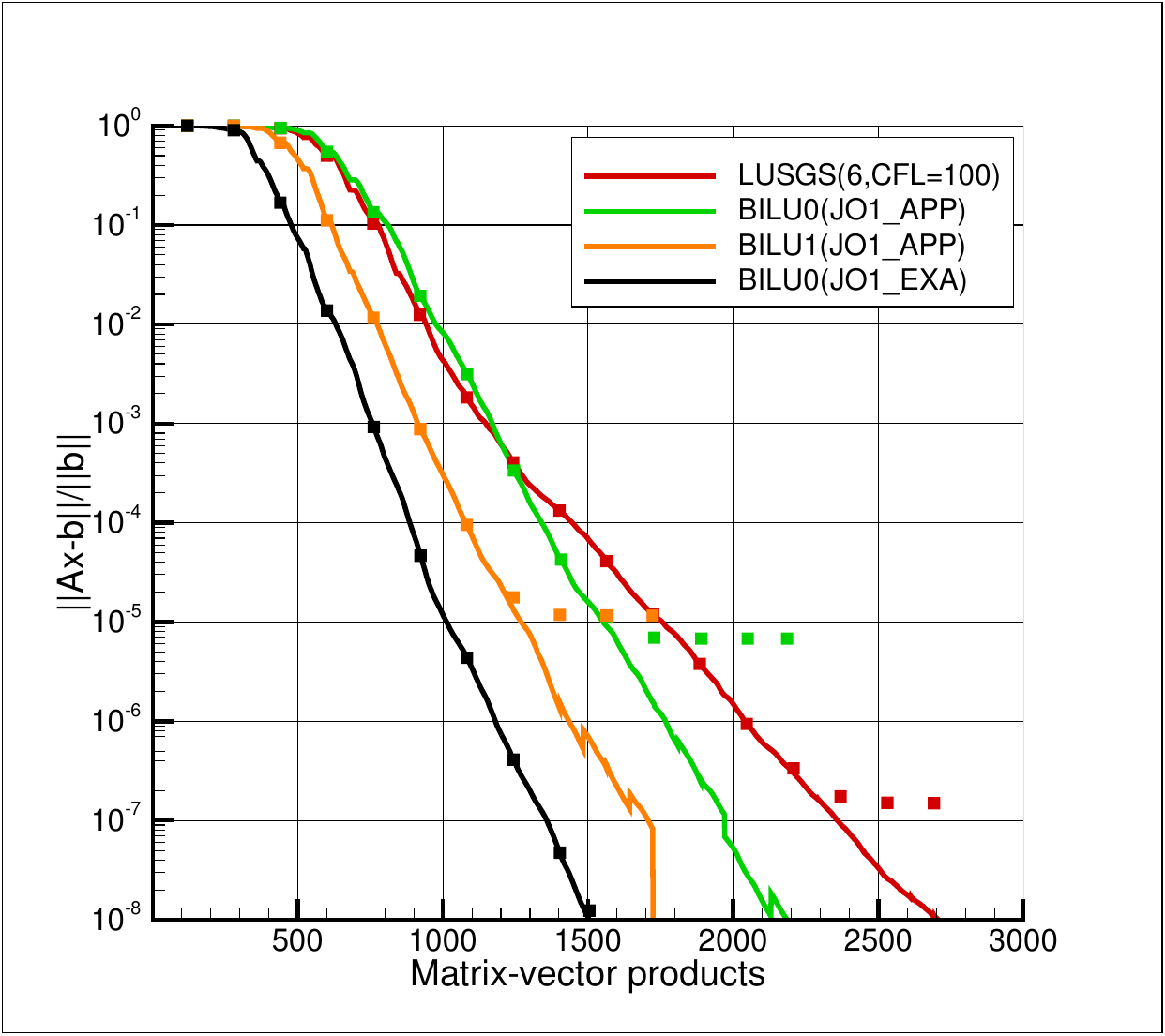}
		\caption{}
		\label{fig:cvg_GMRESDR_120_40_impact_preconditioning_strategy_matvecprod_cvg_erratic}
	\end{subfigure}
	\begin{subfigure}{.5\textwidth}
		\centering
		\includegraphics[trim={0.5cm 0.5cm 1.5cm 1.9cm},clip,width=0.99\textwidth]{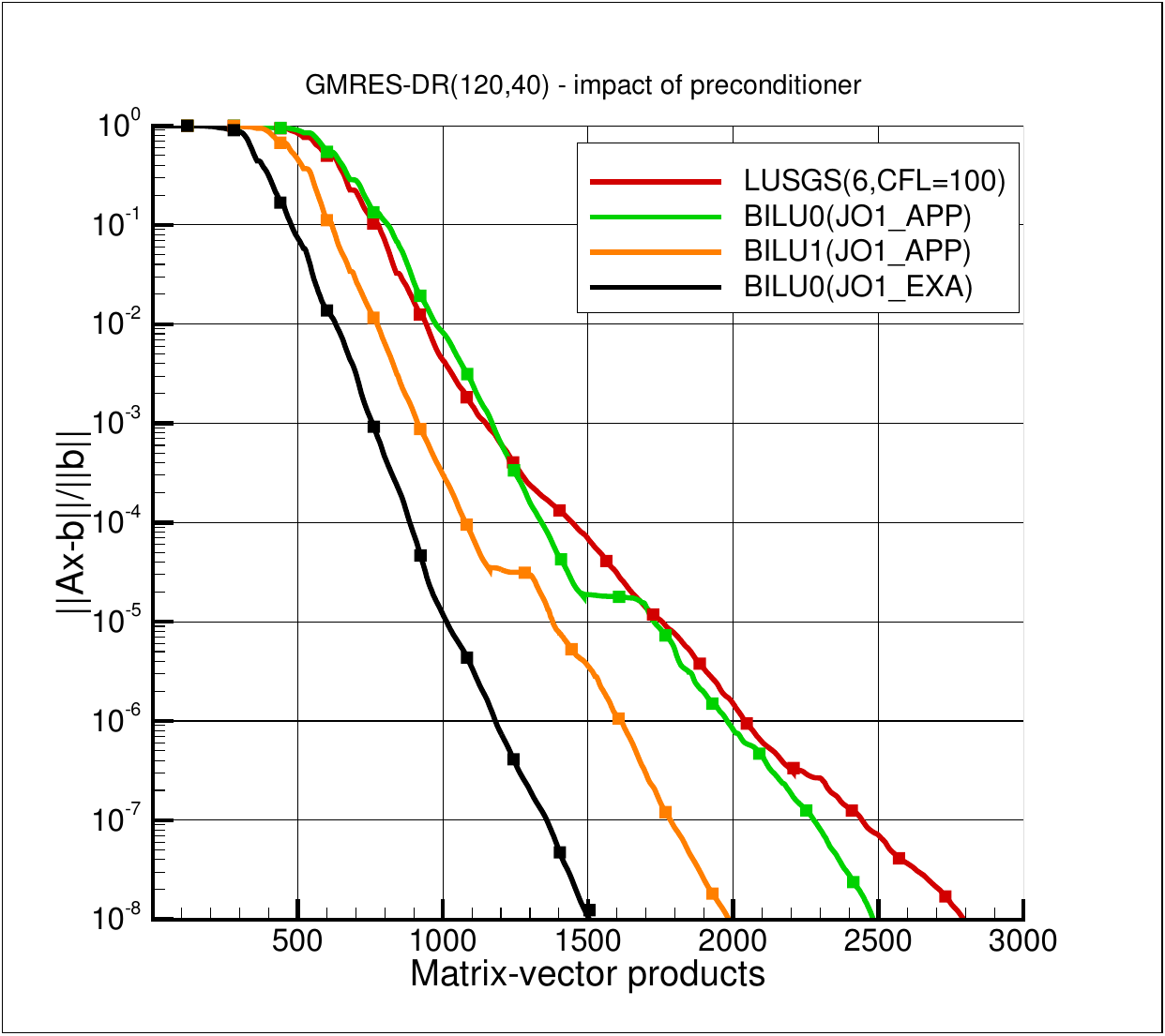}
		\caption{}		
		\label{fig:cvg_GMRESDR_120_40_impact_preconditioning_strategy_matvecprod_resratio}
	\end{subfigure}
	\vspace{0pt}
	\caption{Adjoint relative residual norm convergence history of GMRES-DR(120,40). Impact of various preconditioners. In \ref{fig:cvg_GMRESDR_120_40_impact_preconditioning_strategy_matvecprod_cvg_erratic} some stagnation of the true residual occurs due to propagation of rounding errors during the restarting process. In \ref{fig:cvg_GMRESDR_120_40_impact_preconditioning_strategy_matvecprod_resratio} a cold restart, corresponding to the convergence plateaus, allows to suppress this stagnation at the price of a higher computational cost.}
	\label{fig:cvg_GMRESDR_120_40_impact_preconditioning_strategy_matvecprod}
\end{figure}

\begin{figure}[H]
    \begin{subfigure}{.5\textwidth}
        \raggedright
        \includegraphics[trim={0.5cm 0.5cm 1.5cm 1.9cm},clip,width=0.99\textwidth]{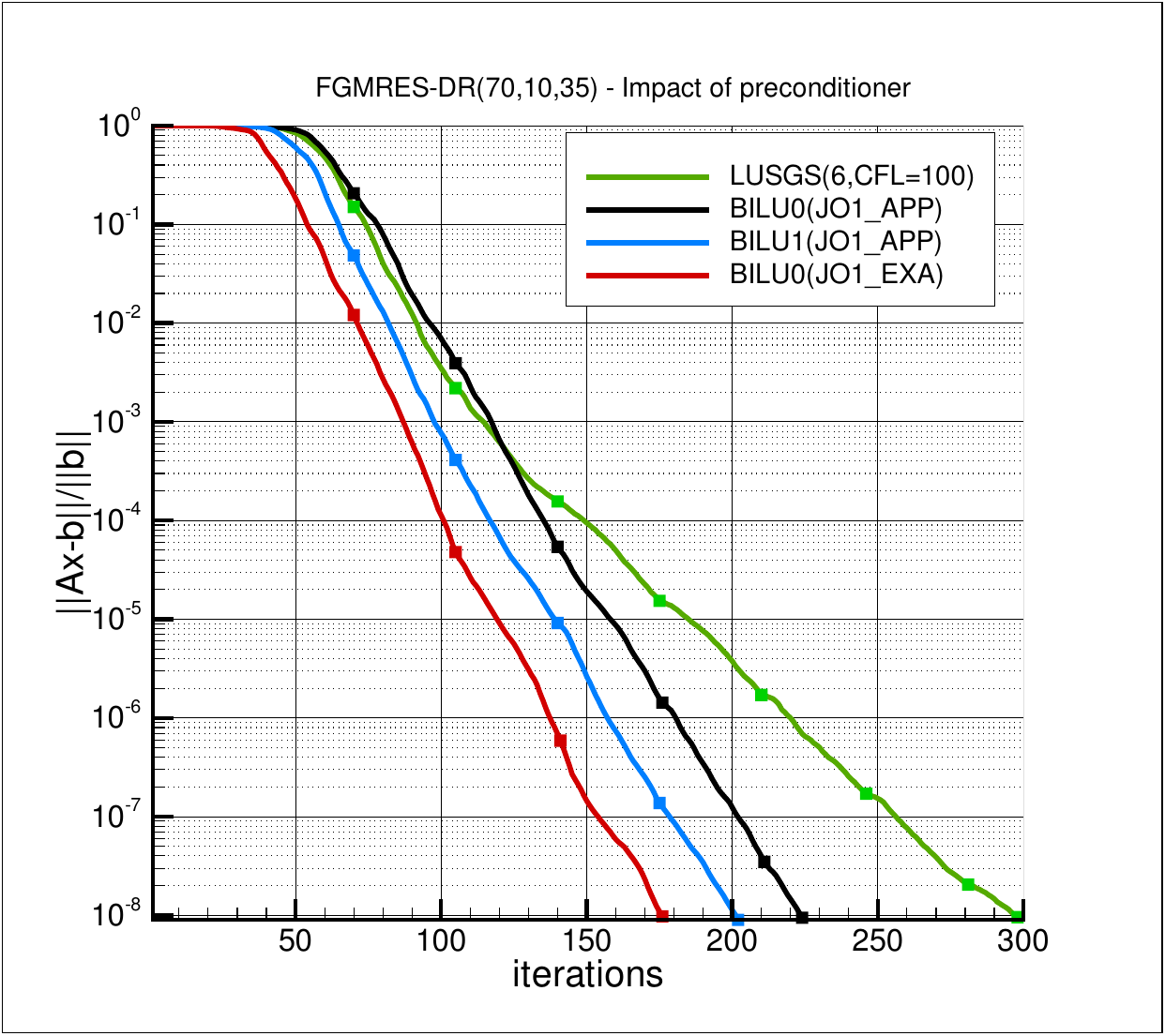}
    \end{subfigure}%
    \begin{subfigure}{.5\textwidth}
        \centering
        \includegraphics[trim={0.5cm 0.5cm 1.5cm 1.9cm},clip,width=0.99\textwidth]{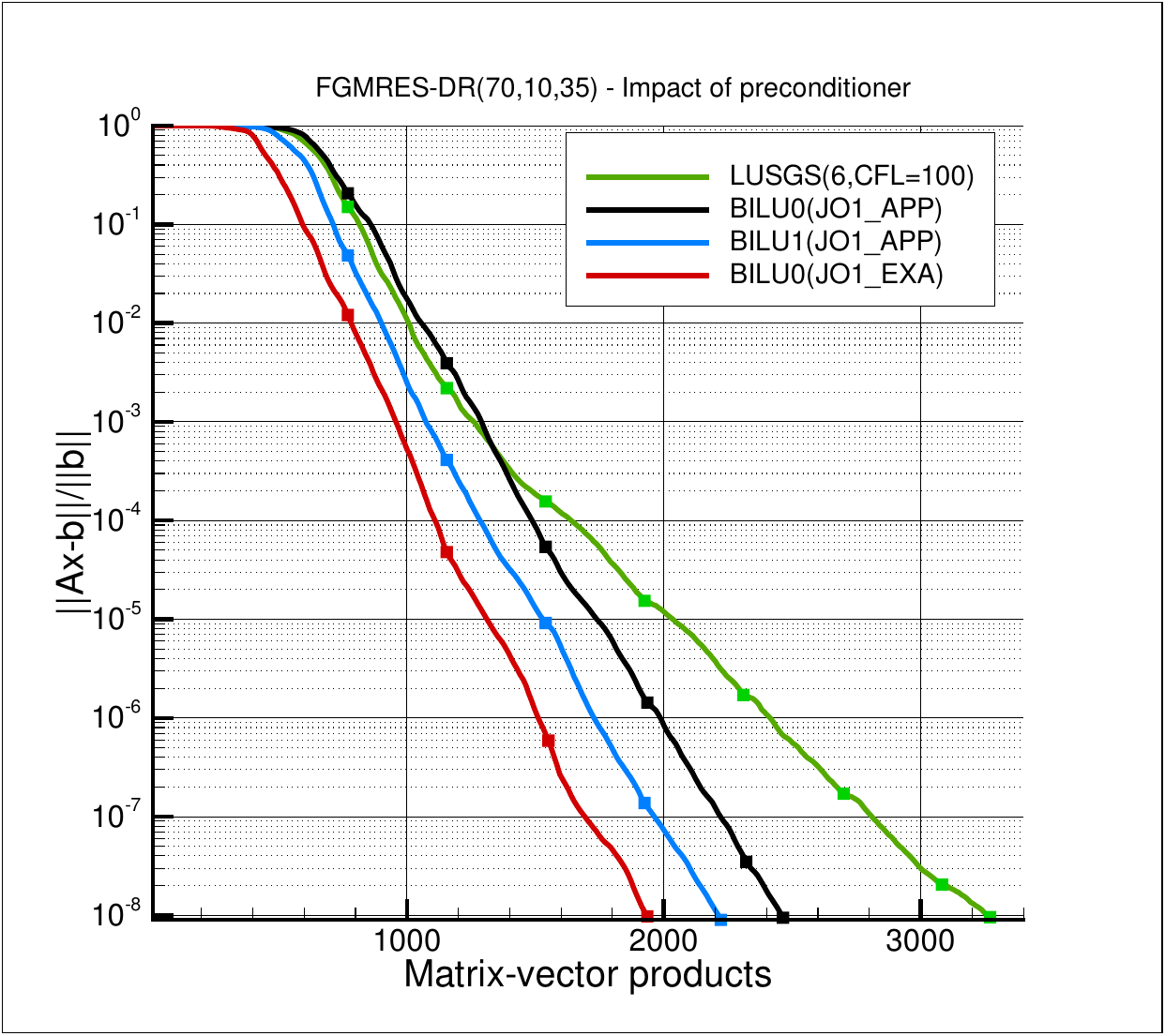}
    \end{subfigure}
    \caption{Adjoint relative residual norm convergence history of FGMRES-DR(70,10,35). Impact of various preconditioners. The true and least-squares residuals match with numerical accuracy}
    \label{fig:cvg_AOC_FGMRESDR_70_10_35_SAlin_impact_preconditioning}
\end{figure}

In this section we apply the GMRES-DR Krylov solver and its flexible variant (variable preconditioning), to the solution of the fluid adjoint system of the M6 wing. We recall that the discretization of the fluid RANS equations relies on an upwind Roe scheme associated to a Van Albada limiter and that the Spalart Allmaras turbulence model is fully linearized. Various stationary preconditioners are considered. First, we consider the legacy RAS+LU-SGS(nrelax, CFL), with nrelax = 6 and CFL = 100 selected as the best options, applied to the first-order approximate Jacobian operator. Then, the BILU(k) decomposition type preconditioners are used. More specifically, BILU(0) is applied to the first-order approximate flux Jacobian operators (approximate and exact) and BILU(1) is applied to the first-order approximate Jacobian matrix only. We start with a  GMRES($m$, $k$) solver with a krylov basis of size $m$ = 120 and a deflation subspace of size $k$ = 40. 
\bigbreak
Figure \ref{fig:cvg_GMRESDR_120_40_impact_preconditioning_strategy_matvecprod_cvg_erratic} plots the convergence of the relative least-squares residual for the above preconditioning strategies. The least-squares residuals all converge up to a prescribed threshold of $10^{-8}$. However, below $10^{-6}$ some erratic behavior happens for the BILU preconditioner applied to the first-order approximate Jacobian operator. The square symbols correspond to the true residuals computed at the end of each GMRES cycle. Except for the preconditioning strategy based on the BILU(0) factorization of the first-order Jacobian matrix, some residual stagnation is observed. This is explained by the propagation of rounding errors in the Arnoldi relation during the restart procedure. A simple remedy for this is to perform a cold restart when the true and least-squares residuals differ too much. In this work, we set this criterion at $\epsilon$ = 5\% of relative discrepancy: $\epsilon= | \mathbf{r}_{true} - \mathbf{r}_{leastsq} | / \mathbf{r}_{true}$ = 0.05. Convergence improvements are clearly illustrated in Figure \ref{fig:cvg_GMRESDR_120_40_impact_preconditioning_strategy_matvecprod_resratio} where the true and least-squares residuals now match accurately. Also, it is not surprising to see a plateau appear when restarting.
\bigbreak
The second set of numerical experiments consisted in the solution of the adjoint linear system but with a nested GMRES strategy where the inner GMRES is right preconditioned by the same variants of stationary preconditioners. Figure \ref{fig:cvg_AOC_FGMRESDR_70_10_35_SAlin_impact_preconditioning} shows the corresponding convergence histories of the relative and least-squares residuals in terms of number of iterations and matrix-vector products. The FGMRES-DR($m$, $m_i$, $k$) solver is applied with a size of Krylov basis $m$ = 70, a size $m_i$ = 10 of the inner Krylov basis and a size $k$ = 40 for the deflation subspace. This time, the true and least-squares residuals always match which confirms the superior numerical stability of the flexible Krylov solver. For both standard and flexible GMRES-DR, the best convergence rate is obtained using the BILU(0) factorization of the exact first-order Jacobian matrix. FGMRES-DR reaches the tolerance threshold in 175 iterations compared to 1500 for GMRES-DR. However, the number of matrix-vector products is lower for the latter (1500 compared to 1920), illustrating a higher computational cost associated to FGMRES-DR for this specific numerical test.

\subsection{Application of GMRES-DR and FGMRES-DR to the solution of the fluid-structure coupled-adjoint system}

We are now interested in the application of GMRES-DR and FGMRES-DR to the solution of the coupled-adjoint system \ref{Coupled_Adjoint_system}. The solution process is presented in algorithm \ref{alg:LBGS-algorithm}. The fluid-structure coupling is triggered when the ratio $\rho = \mathbf{r}_i / \mathbf{r}_{i-1}$ between the residual $\mathbf{r}_i$ at the end of the current cycle and the residual  $\mathbf{r}_{i-1}$ at the end of the previous cycle is less than a prescribed threshold. Relaxation of the structural adjoint vector (parameter $\theta_s$ at step \ref{alg:Relaxation_step} in algorithm \ref{alg:LBGS-algorithm}) can be activated to stabilize the convergence of the partitioned solver. A popular choice is to adapt $\theta_s$ dynamically by using an autoregressive predictor such as the well-known Aitken strategy \cite{irons1969version, kuttler2008fixed}.

\begin{figure}[H]
	\centering
	\includegraphics[trim={0.5cm 0.5cm 1.5cm 1.8cm},clip,width=0.6\textwidth]{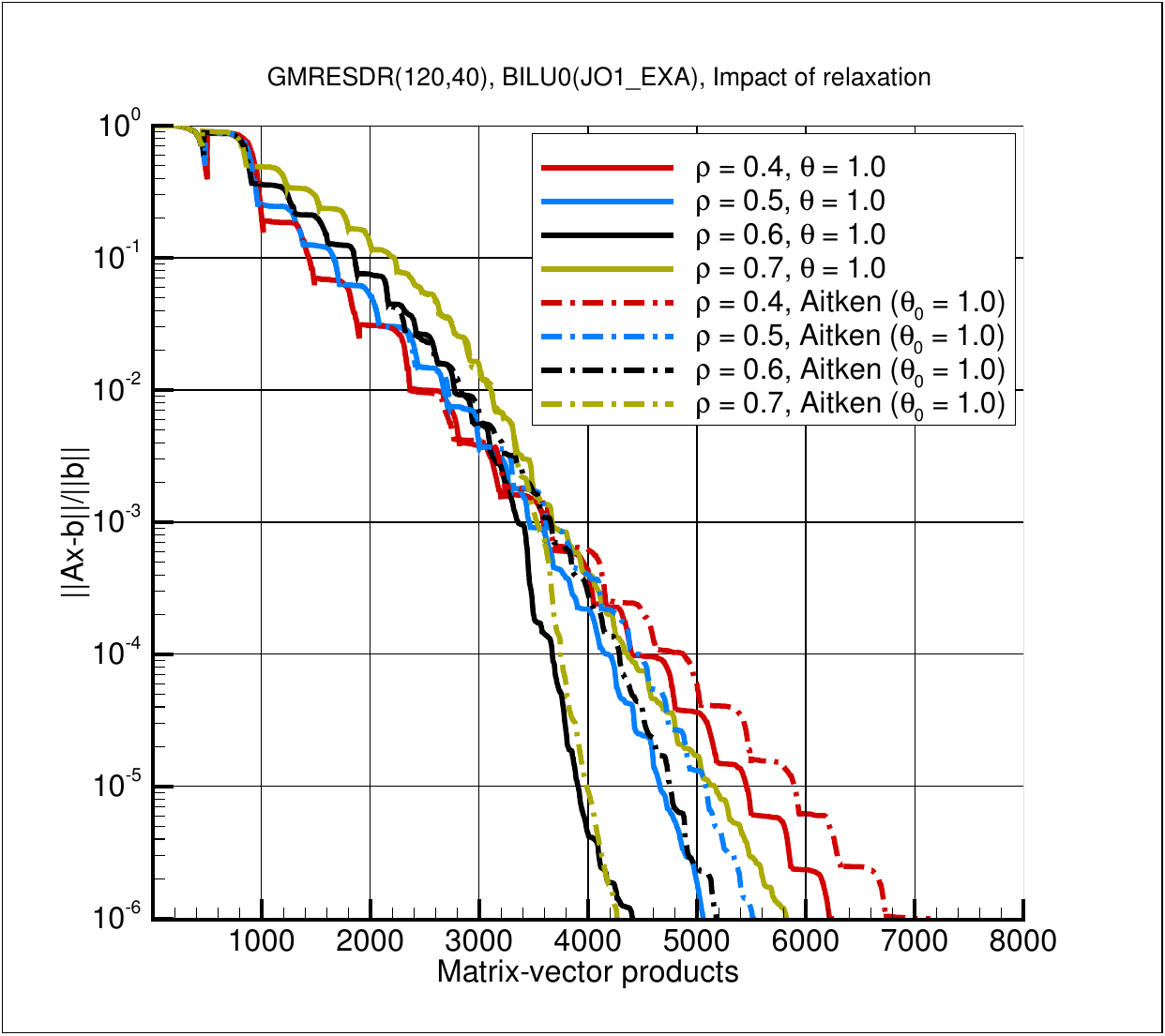}
	\caption{Coupled-adjoint relative residual norm convergence history of GMRES-DR(120,40) right-preconditioned by BILU0($\textbf{J}^{EXA}_{O1}$). Impact of constant and dynamic relaxation parameter.}
	\label{fig:cvg_AOC_GMRESDR_120_40_relaxation_parameter_study}
\end{figure}

 In order to select appropriate values for $\rho$ ans $\theta_s$, we carried out a parameter study with $\rho \in [0.4,0.5,0.6,0.7]$ in combination with a constant relaxation parameter $\theta_s=1.0$ on one hand, and a dynamic relaxation parameter $\theta_s^{*}$, with $\theta_s^0=1.0$, on the other hand. The associated convergence curves are presented in Figure \ref{fig:cvg_AOC_GMRESDR_120_40_relaxation_parameter_study}. For a constant relaxation parameter, the best convergence rate is achieved for a residual ratio $\rho=0.6$ (plain black line). The impact of dynamic relaxation (dash-dotted lines) is systematically negative except for $\rho=0.7$ where the efficiency of the combination of $\rho=0.6$ with $\theta_s=1.0$ is retrieved. In light of these results, it seems more appropriate to adopt a constant parameter relaxation strategy for our specific test case. However, higher-order Aitken techniques could be applied to better capture the nonlinear profile of the convergence, but this is out of the scope of this paper. Consequently, in the remaining of this document we will select  $\rho=0.6$ and $\theta_s=1.0$ for the partitioned solver.

\begin{figure}[H]
    \centering
    \includegraphics[trim={0.5cm 0.5cm 1.5cm 1.8cm},clip,width=0.6\textwidth]{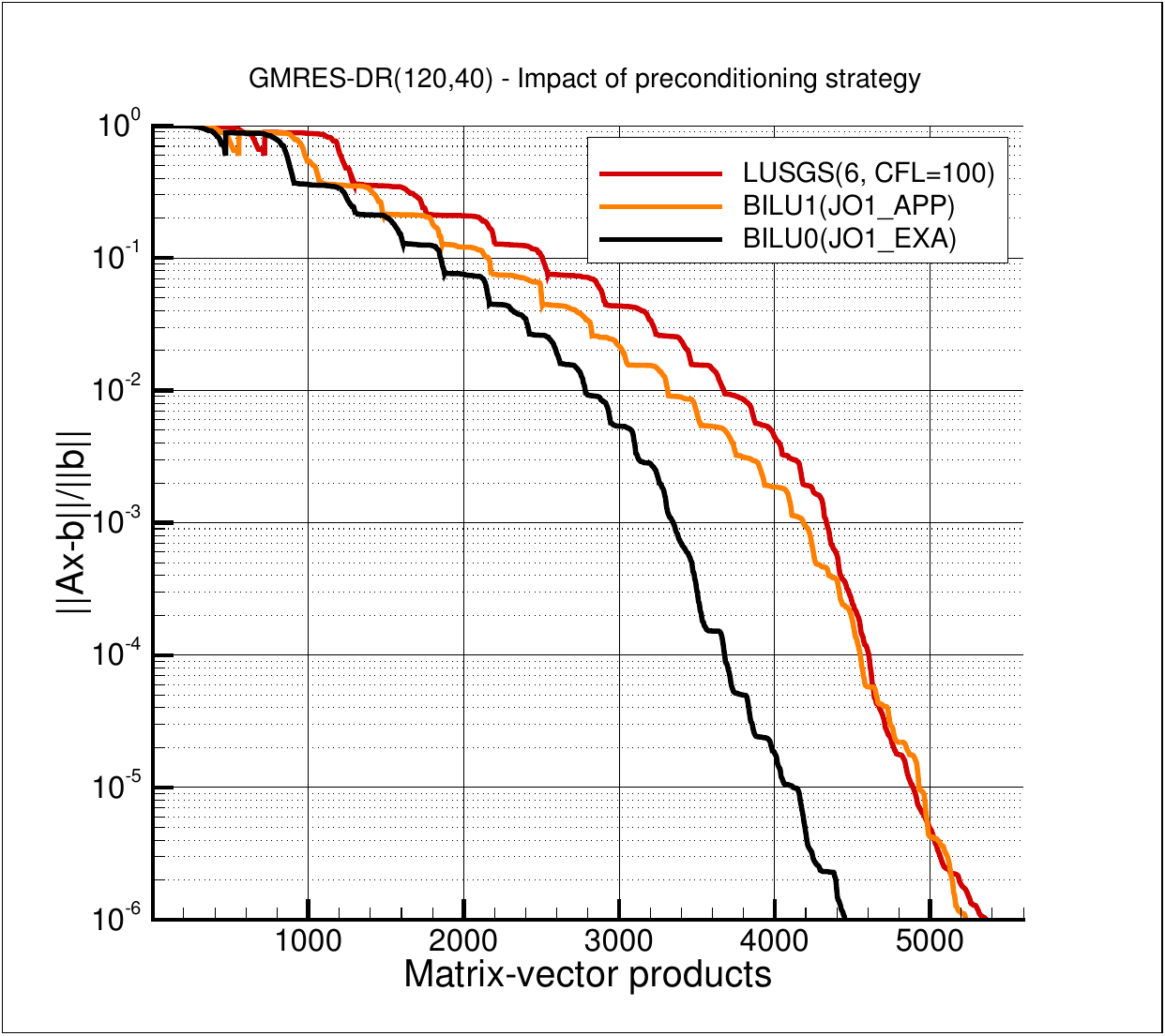}
    \caption{Coupled-adjoint relative residual norm convergence history of GMRES-DR(120,40). Impact of various preconditioners.}
    \label{fig:cvg_AOC_GMRESDR_120_40_SAlin_impact_preconditioning_matvecprod}
\end{figure}

In Figure \ref{fig:cvg_AOC_GMRESDR_120_40_SAlin_impact_preconditioning_matvecprod} we plot the convergence of the fluid relative least-squares residual for GMRES-DR(120,40) associated to several preconditioners.  
The preconditioners RAS+LUSGS(6,100) and BILU($1$) both applied to the approximate first-order Jacobian operator perform similarly. As before, BILU(0) applied to the exact first-order Jacobian shows the best efficiency. Contrary to the purely fluid case, the cycles of the coupled solver are often shorter which disadvantages the less efficient preconditioners. The same computation is repeated using a FGMRES-DR(70,10,35) Krylov solver. This time, BILU(1) applied to the approximate first-order Jacobian operator is not competitive and BILU(0) applied to the exact first-order Jacobian far outperforms other strategies. Again, the nested solver requires roughly 4 times less iterations than its non-flexible counterpart whereas the number of matrix-vector products is 10800 compared to 4400.

\vspace{-\baselineskip}

\begin{figure}[H]
    \begin{subfigure}{.5\textwidth}
        \raggedright
        \includegraphics[trim={1.cm 0.5cm 1.5cm 1.8cm},clip,width=0.99\textwidth]{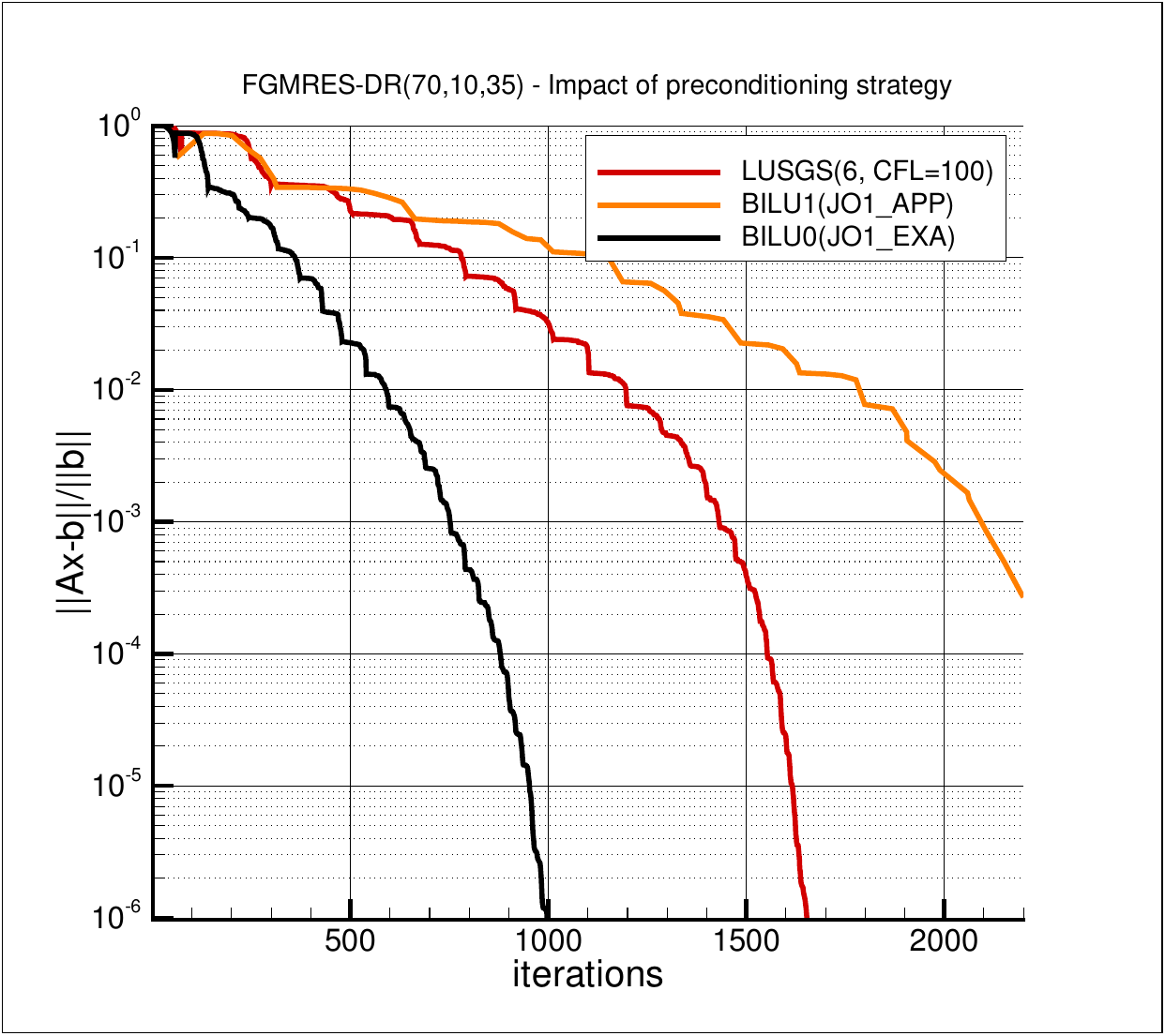}
    \end{subfigure}%
    \begin{subfigure}{.5\textwidth}
        \centering
        \includegraphics[trim={1.cm 0.5cm 1.5cm 1.8cm},clip,width=0.99\textwidth]{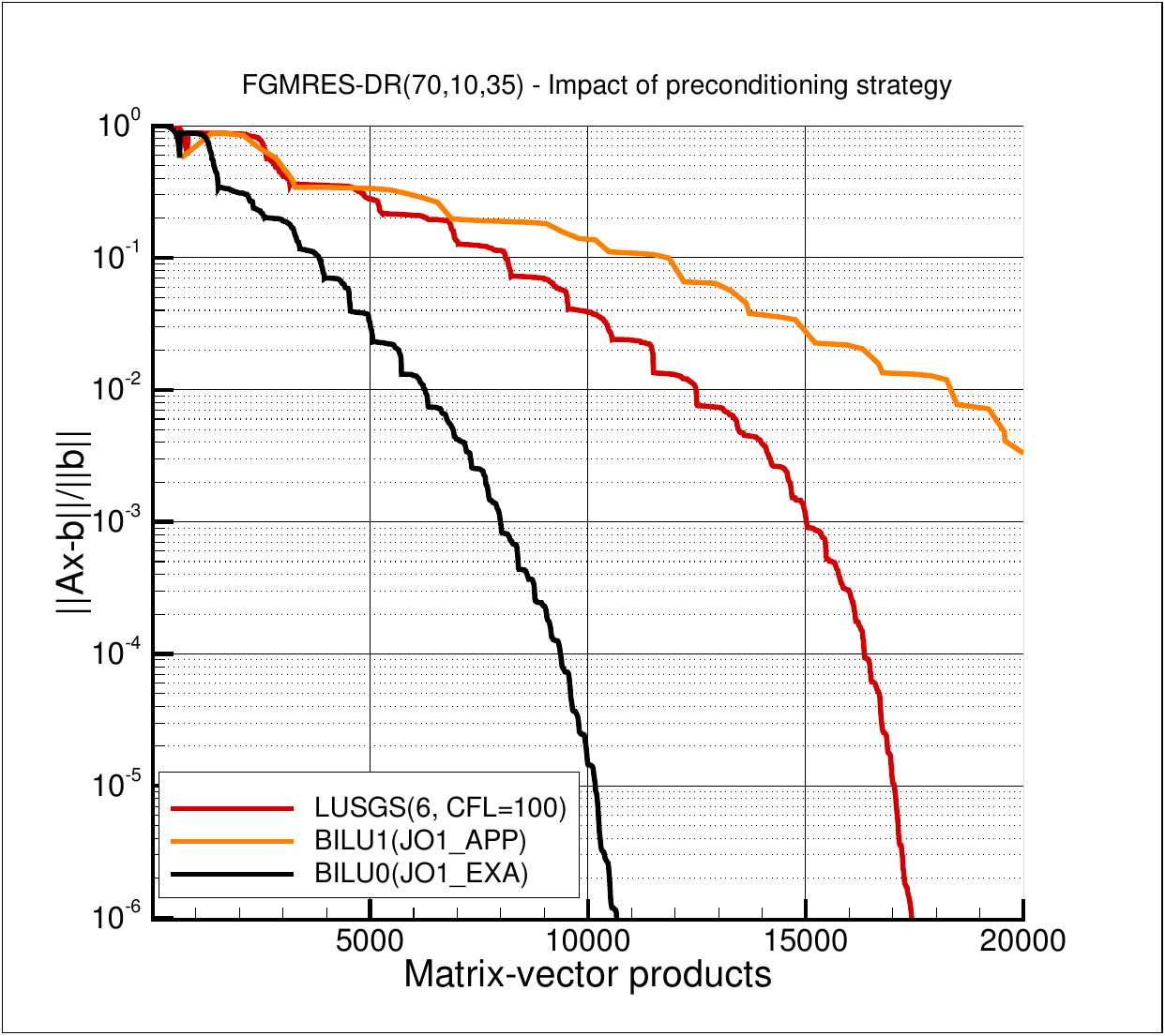}
    \end{subfigure}
    \caption{Coupled-adjoint relative residual norm convergence history of FGMRES-DR(70,10,35). Impact of various preconditioners.}
    \label{fig:cvg_AOC_FGMRESDR_70_10_35_SAlin_impact_preconditioning_iterations}
\end{figure}

\section{Generalized Conjugate Residual with inner Orthogonalization and Deflated Restarting: GCRO-DR} \label{GCRODR_theory}

As already explained, the current implementation of our partitioned solver does not take advantage of spectral information produced during the Krylov solver instances applied to the sequence of preceding right-hand sides. An advanced inner-outer GMRES-DR solver is used for the approximate solution of the fluid block between consecutive fluid-structure couplings but a simple restart is performed when the structural source term is updated. Unfortunately, deflated restarting based on subspace augmentation by appending approximate Ritz vectors to the Krylov subspace is not suitable for sequences of linear systems. Indeed, if the Ritz vectors are obtained from a previous linear system with another matrix or even another right-hand side, the concatenation of the recycled subspace to the current Krylov subspace does not form a Krylov subspace for the current problem. It is necessary to introduce a new Krylov solver that uses recycling of any given subspace without restriction. A famous one is the generalized conjugate residual with inner orthogonalization (GCRO). It belongs to the family of inner-outer methods \citep{axelsson1996iterative} where the outer method is based on GCR, a minimum residual norm method proposed by Eisenstat,  Elman and Schultz \citep{eisenstat1983variational}. The inner solver is GMRES applied to a projected system matrix.

\subsection{The GCRO Krylov solver}
As previously mentioned, we focus on the solution of a sequence of linear systems with a varying right-hand side only. This section briefly introduces the Generalized Conjugate Residual with inner Orthogonalization (GCRO) algorithm with deflated restarting. In this framework, deflation can re-use spectral information from a previous cycle or from a previous linear system.
We start by recalling the original formulation of the generalized conjugate residual method (GCR)\citep{eisenstat1983variational}. The idea is to introduce the concept of optimality for the solution residual. We want to solve
\begin{equation}
	\mathbf{A}\mathbf{x}^{(s)} = \mathbf{b}^{(s)}, \quad s=1,2,...
	\label{eq7}
\end{equation}

\noindent where $\mathbf{A} \in \mathbb{R}^{n \times n}$ and $\mathbf{b}^{(s)} \in \mathbb{R}^{n}$ changes from one system to the next.

The GCR method relies on a given full-rank matrix $\mathbf{U}_{k} \in \mathbb{R}^{n \times k}$ and an orthonormal matrix $\mathbf{C}_{k} \in \mathbb{R}^{n \times k}$ as the image of $\mathbf{U}_{k}$ by $\mathbf{A}$ satisfying the relations
\begin{align}
	\mathbf{A}\mathbf{U}_{k}         &= \mathbf{C}_{k}, \label{eq8} \\
	\mathbf{C}^{T}_{k}\mathbf{C}_{k} &= \mathbf{I}_{k}. \label{eq9}
\end{align}

For the sake of understanding, suppose that we want to solve the first system of the sequence of linear systems, i.e., $s=1$. Given an initial guess $\mathbf{x}_{0}^{(1)}$, the principle is to compute an approximation to the solution $\mathbf{x}^{(1)} \in \mathbf{x}_{0}^{(1)} + \operatorname{range}(\mathbf{U}_{k})$ that minimizes the corresponding residual norm over the approximation space $\operatorname{range}(\mathbf{U}_{k})$. More precisely, GCR solves the following minimization problem
\begin{equation}
	\mathbf{x}^{(1)} = \underset{\mathbf{x} \in \mathbf{x}^{(1)}_{0} + \operatorname{range}(\mathbf{U}_{k})}{\textnormal{argmin}} \| \mathbf{b}^{(1)} - \mathbf{Ax} \|_2,
	\label{eq10}
\end{equation}

The optimal solution of (\ref{eq10}) over the subspace $\mathbf{x}^{(1)}_{0} + \operatorname{range}(\mathbf{U}_{k})$ is defined by

\begin{equation}
	\mathbf{x}^{(1)} = \mathbf{x}^{(1)}_{0} + \mathbf{U}_{k}\mathbf{C}^{T}_{k}\mathbf{r}_{0}^{(1)}
	\label{eq11}
\end{equation}

\noindent where $ \mathbf{r}_{0}^{(1)} = \mathbf{b}^{(1)} - \mathbf{A}\mathbf{x}^{(1)}_{0}$. Consequently, the corresponding residual vector satisfies

\begin{equation}
	\mathbf{r}^{(1)}_k = \mathbf{b}^{(1)}-\mathbf{A}\mathbf{x}^{(1)} = (\mathbf{I} - \mathbf{A}\mathbf{U}_k\mathbf{C}_k^T)\mathbf{r}^{(1)}_{0} = (\mathbf{I} - \mathbf{C}_{k}\mathbf{C}^{T}_{k})\mathbf{r}^{(1)}_{0}, \quad \mathbf{r}^{(1)}_k \perp \operatorname{range}(\mathbf{C}_{k}).
	\label{eq12}
\end{equation}

The orthogonality of the residual $ \mathbf{r}^{(1)}_k $ to the subspace $\mathcal{AK}_{k}(\mathbf{A},\mathbf{r}^{(1)}_{0})$ spanned by the columns of $\mathbf{C}_{k}$ is known as the optimality property of the residual. In practice, GCR is not considered as a means to solve the linear system. Instead, for the first system or when no spectral information is available from a previous solve, it is replaced by GMRES which computes an implicit representation of the matrices $\mathbf{U}_{k}$ and $\mathbf{C}_{k}$ \citep{de1999truncation}.
\bigbreak
In order to simplify the notations, we temporarily omit the index of the current system in the equations. Given an orthonormal basis $\te{C}_{k}$ of an outer subspace and the corresponding residual $ \ve{r}_k=(\matI - \te{C}_{k}\te{C}^{T}_{k})\ve{r}_0 $ after $k$ steps of GCR, GCRO($m$) obtains the next iterates $\ve{r}_{k+1}$ and $\ve{x}_{k+1}$ by performing $m$ steps of GMRES applied to the projected operator $\te{A}_{\te{C}_k}=(\matI - \te{C}_{k}\te{C}^{T}_{k})\te{A}$, thereby maintaining optimality of the inner residual with respect to the outer space. Let $\te{V}_{m+1} \in \mathbb{R}^{n \times (m+1)}$ be an orthonormal basis for $\mathcal{K}_{m+1}(\te{A}_{\te{C}_k},\ve{r}_{k})$ with $\ve{v}_1=\ve{r}_{k}/\| \ve{r}_{k} \|$, if we apply an inner GMRES we have the following Arnoldi relation:
\begin{equation}
\te{A}_{\te{C}_k}\te{V}_m = (\matI - \te{C}_{k}\te{C}^{T}_{k})\te{A}\te{V}_m = \te{V}_{m+1}\bar{\te{H}}_{m} \quad \text{with} \quad \bar{\te{H}}_{m} = \te{V}_{m+1}^{T} \te{A}_{\te{C}_k}\te{V}_m, \quad \bar{\te{H}}_{m} \in \mathbb{R}^{(m+1) \times m}
\label{eq:ArnoldiGCRO}
\end{equation}

This equation can be expanded as
\begin{equation}
\mathbf{A}\mathbf{V}_m = \mathbf{C}_{k}\mathbf{B}_m + \mathbf{V}_{m+1}\bar{\mathbf{H}}_{m} \quad \text{with} \quad \mathbf{B}_m = \mathbf{C}_k^{T} \mathbf{A} \mathbf{V}_m, \quad \te{B}_{m} \in \mathbb{R}^{k \times m}
\label{eq:ExpandedArnoldiGCRO}
\end{equation}

At the end of the inner GMRES cycle, using \eqref{eq:ArnoldiGCRO}, the residual can be written as
\begin{equation}
\mathbf{r}_{k+1} = \mathbf{r}_k - \mathbf{A}_{\mathbf{C}_k}\mathbf{V}_{m}\mathbf{y}_m = \mathbf{r}_{k} - \mathbf{V}_{m+1}\bar{\mathbf{H}}_{m}\mathbf{y}_m,
\label{eq:GCRO_residual}
\end{equation}

\noindent and the corresponding solution reads
\begin{equation}
\mathbf{x}_{k+1} = \mathbf{x}_{k} + \mathbf{A}^{-1}\mathbf{A}_{\mathbf{C}_{k}}\mathbf{V}_{m}\mathbf{y}_{m},
\label{eq:GCRO_solution}
\end{equation}

\noindent where, using \eqref{eq8}, we have $ \mathbf{A}^{-1}\mathbf{A}_{\mathbf{C}_{k}}=\matI-\mathbf{U}_k \mathbf{C}^{T}_{k} \mathbf{A} $. Using \eqref{eq11} and \eqref{eq:GCRO_solution} we obtain the expression for the solution at the end of the first cycle:
\begin{align}
\mathbf{x}_{k+1} &= \mathbf{x}^{(1)}_{0} + \mathbf{U}_{k}\mathbf{C}^{T}_{k}\mathbf{r}_{0}^{(1)} + (\matI-\mathbf{U}_k \mathbf{C}^{T}_{k} \mathbf{A}) \mathbf{V}_m \mathbf{y}_m \nonumber \\
 &= \mathbf{x}^{(1)}_{0} + \mathbf{U}_{k}\mathbf{C}^{T}_{k} \left( \mathbf{r}_{0}^{(1)} - \mathbf{A} \mathbf{V}_m \mathbf{y}_m \right) + \mathbf{V}_m \mathbf{y}_m  \nonumber \\
 &= \mathbf{x}^{(1)}_{0} + \mathbf{U}_{k} \mathbf{z}_k + \mathbf{V}_m \mathbf{y}_m,
 \label{solution_GCR}
\end{align}
with
\begin{equation}
	\mathbf{z}_{k} = \mathbf{C}^{T}_{k}\mathbf{r}_{0}^{(1)} - \mathbf{B}_m \mathbf{y}_m, \quad \mathbf{B}_m = \mathbf{C}^{T}_{k}\mathbf{AV}_m.
\end{equation}

We easily show that the residual computed in the inner GMRES equals the true outer residual:
\begin{equation}
	\mathbf{r}_{k+1} = \mathbf{b} - \mathbf{Ax}_{k+1} = \mathbf{b} - \mathbf{Ax}_{k} - \mathbf{A}_{\mathbf{C}_k}\mathbf{V}_{m}\mathbf{y}_m = \mathbf{r}_{k} - \mathbf{V}_{m+1}\bar{\mathbf{H}}_{m}\mathbf{y}_m.
\end{equation}

Because $ \operatorname{range}(\mathbf{C}_{k}) \bot \operatorname{range}(\mathbf{V}_{m+1}) $ by construction, it follows that $\mathbf{x}_{k+1} \in \mathbf{x}^{(1)}_{0} + \operatorname{range}(\mathbf{U}_{k}) \oplus \operatorname{range}(\mathbf{V}_m)$ is the unique solution of the global residual minimization problem over the subspace spanned by the complementary subspaces  $\operatorname{range}(\mathbf{U}_{k})$ and $\operatorname{range}(\mathbf{V}_m)$ (see Theorem 2.2 and its proof in \cite{de1996nested}).

Now, introducing the orthogonal decomposition $ \bar{\mathbf{H}}_m = \bar{\mathbf{Q}}_m \mathbf{R}_m $, the vector of reduced coordinates $ \mathbf{y}_m $ is such that
\begin{equation}
\mathbf{y}_m = \underset{\mathbf{y} \in \mathbb{R}^m}{\operatorname{argmin}} \| \mathbf{r}_{k+1} \|_2 = \underset{\mathbf{y} \in \mathbb{R}^m}{\operatorname{argmin}} \| \mathbf{r}_k - \mathbf{A}_{\mathbf{C}_k}\mathbf{V}_{m}\mathbf{y} \|_2 = \mathbf{R}_m^{-1}\bar{\mathbf{Q}}_m^{T} \| \mathbf{r}_k \|_2 \mathbf{e}_1
\end{equation}

In the above Arnoldi process \eqref{eq:ExpandedArnoldiGCRO} of the inner GMRES($m$), the vectors $\mathbf{Av}_i$ are first orthogonalized against $\mathbf{C}_k$, thus constructing $\mathbf{V}_{m+1}$ such that $\mathbf{C}_k^{T}\mathbf{V}_{m+1}=\mathbf{0}$. 

\subsection{The GCRO-DR Krylov solver}

So far, we have not yet explained how the matrix $\mathbf{C}_{k}$ is built. A fundamental difference between GCRO and GCRO-DR is that GCRO is a nested Krylov solver where the outer spaces $\mathbf{U}_k$ and $\mathbf{C}_k$ keep growing and GMRES($m$) is applied at each iteration $k$, whereas GCRO-DR supplements given deflation subspaces derived from a previous cycle with ($m-k$) iterations of GMRES. 
To comply to the GCRO-DR formalism, we first slightly change the notation in \eqref{eq:ArnoldiGCRO} and \eqref{eq:ExpandedArnoldiGCRO} to point out that we perform $m-k$ steps (instead of $m$) of inner GMRES. This produces the Arnoldi relation
\begin{equation}
\mathbf{A}_{\mathbf{C}_k}\mathbf{V}_{m-k} = \mathbf{V}_{m-k+1}\bar{\mathbf{H}}_{m-k} \quad \Leftrightarrow \quad \bar{\te{H}}_{m-k} = \te{V}_{m-k+1}^{T} \te{A}_{\te{C}_k}\te{V}_{m-k}, \quad \bar{\te{H}}_{m-k} \in \mathbb{R}^{(m-k+1) \times (m-k)}
\label{eq:FlexibleGCRO1}
\end{equation}

Combining \eqref{eq8} and \eqref{eq:FlexibleGCRO1} we have
\begin{equation}
\mathbf{A} \left[ \mathbf{U}_k \quad \mathbf{V}_{m-k} \right] = [ \mathbf{C}_k \quad \mathbf{V}_{m-k+1} ]
\begin{bmatrix}
\mathbf{I}_{k} & \mathbf{B}_{m-k} \\
\mathbf{0}     & \bar{\mathbf{H}}_{m-k}
\end{bmatrix},
\label{eq:ExpandedArnoldiGCRO1}
\end{equation}
with $ \mathbf{B}_{m-k} = \mathbf{C}_k^{T} \mathbf{AV}_{m-k}, \mathbf{B}_{m-k} \in \mathbb{R}^{k \times (m-k)}$.

Numerical tests suggest that the rightmost matrix in \eqref{eq:ExpandedArnoldiGCRO1} is ill-conditioned. To reduce unnecessary ill-conditioning we can compute the diagonal matrix $\mathbf{D}_k = \operatorname{diag}(\|\textbf{u}_1\|^{-1}, \|\textbf{u}_2\|^{-1},\ldots, \|\textbf{u}_k\|^{-1})$ such that $\tilde{\mathbf{U}}_k = \mathbf{U}_k \mathbf{D}_k$, has unit columns. Defining
\begin{equation*}
    \hat{\mathbf{V}}_m = [\tilde{\mathbf{U}}_k \quad \mathbf{V}_{m-k}], \quad
    \hat{\mathbf{W}}_{m+1} = [\mathbf{C}_k \quad \mathbf{V}_{m-k+1}], \quad
    \bar{\mathbf{H}}_m = \begin{bmatrix}
        \mathbf{D}_k & \mathbf{B}_{m-k} \\
        \mathbf{0} & \bar{\mathbf{H}}_{m-k}
    \end{bmatrix},
\end{equation*}
the Arnoldi relation \eqref{eq:ExpandedArnoldiGCRO1} can be cast into
\begin{equation}
    \mathbf{A}\hat{\mathbf{V}}_m = \hat{\mathbf{W}}_{m+1} \bar{\mathbf{H}}_m.
   	\label{eq:GeneralizedArnoldi}
\end{equation}

From \eqref{solution_GCR} we have $\mathbf{x}_1 = \mathbf{x}_0 + \tilde{\mathbf{U}}_k \mathbf{z}_k + \mathbf{V}_{m-k} \mathbf{y}_{m-k} = \mathbf{x}_0 + \hat{\mathbf{V}}_m \mathbf{y}_m$ with $ \mathbf{y}_m = {\operatorname{argmin}}_{\mathbf{y}} \|\bar{\mathbf{H}}_{m}\mathbf{y}-\hat{\mathbf{W}}_{m+1}^T \mathbf{r}_0\|_2 $.
The expanded form of this least-squares problem reads
\begin{equation}
\mathbf{y}_m = [ \mathbf{z}_k, \mathbf{y}_{m-k} ] = {\operatorname{argmin}}_{\mathbf{z,y}} \left\| \begin{bmatrix}
	\mathbf{D}_k & \mathbf{B}_{m-k} \\	\mathbf{0} & \bar{\mathbf{H}}_{m-k} \end{bmatrix}
\begin{bmatrix} \mathbf{z} \\ \mathbf{y}  \end{bmatrix} - \begin{bmatrix} \mathbf{C}_k^T\mathbf{r}_0 \\ \mathbf{V}_{m-k+1}^T \mathbf{r}_0  \end{bmatrix} \right\|_2
\label{GCRO_LeastSquares}
\end{equation}

In \cite{de1996nested} and \cite{soodhalter2020} the authors suggest to solve \eqref{GCRO_LeastSquares} blockwise, first for $\mathbf{y}_{m-k}$ and then set $\mathbf{z}_k = \mathbf{D}_k^{-1} [\mathbf{C}_k^T \mathbf{r}_0 - \mathbf{B}_{m-k}\mathbf{y}_{m-k}]$. Recall that the first Krylov vector in $\mathbf{V}_{m-k+1}$ is $\mathbf{v}_{k+1}=(\mathbf{I}-\mathbf{C}_k\mathbf{C}_k^T)\mathbf{r}_0 $, which simplifies the expression for the right-hand side of the least-squares sub-problem to $\mathbf{V}_{m-k+1}^T \mathbf{r}_0 = \mathbf{V}_{m-k+1}^T (\mathbf{I}-\mathbf{C}_k\mathbf{C}_k^T) \mathbf{r}_0 = \beta \mathbf{e}_1^{m-k+1}$, where $\beta = \| (\mathbf{I}-\mathbf{C}_k\mathbf{C}_k^T)\mathbf{r}_0 \|_2$ and $\mathbf{e}_1^{m-k+1}$ is the first unit vector of $\mathbb{R}^{m-k+1}$.


We have now completed the process for the first linear system which formulates the solution from the complementary subspaces $\tilde{\mathbf{U}}_k$ and $\mathbf{V}_{m-k}$. At this stage, we introduce deflated restarting within the GCRO framework to build the outer space $\mathbf{C}_k$ for the next cycle. Like GMRES-DR, the deflation subspace is spanned by the $k$ harmonic Ritz vectors corresponding to harmonic Ritz values of smallest magnitude. Using definition \ref{HarmonicRitzpair} with $\mathcal{U} \equiv \operatorname{span}\lbrace \hat{\mathbf{V}}_m \rbrace $ and $\mathbf{B} \equiv \mathbf{A}$ we can write $\mathbf{y} = \hat{\mathbf{V}}_m \mathbf{g}$ which gives the orthogonality relation
\begin{equation}
	(\mathbf{A}\hat{\mathbf{V}}_m)^T(\mathbf{A}\hat{\mathbf{V}}_m \mathbf{g} -\theta \hat{\mathbf{V}}_m \mathbf{g}) = \mathbf{0}.
	\label{eq:Ritz_problem_GCRO}
\end{equation}

Now using \eqref{eq:GeneralizedArnoldi} in \eqref{eq:Ritz_problem_GCRO} and recalling that $\hat{\mathbf{W}}_{m+1}$ is orthonormal, we obtain the generalized eigenvalue problem 
\begin{equation}
	\label{eq:GCRO_eigenvalue_problem}
	\boxed{\bar{\mathbf{H}}^T_m \bar{\mathbf{H}}_m \mathbf{g} = \theta \bar{\mathbf{H}}^T_m \hat{\mathbf{W}}^T_{m+1} \hat{\mathbf{V}}_m \mathbf{g}}
\end{equation}

In practice, for better numerical accuracy, instead of computing the smallest eigenvalues of \eqref{eq:GCRO_eigenvalue_problem} we prefer to compute the largest eigenvalues $\theta^{-1}$ of $\bar{\mathbf{H}}^T_m \hat{\mathbf{W}}^T_{m+1} \hat{\mathbf{V}}_m \mathbf{g} = \theta^{-1}  \bar{\mathbf{H}}^T_m \bar{\mathbf{H}}_m \mathbf{g}$.
Since $\hat{\mathbf{W}}_{m+1} = [\mathbf{C}_k \quad \mathbf{V}_{m-k+1}]$, the matrix product $ \hat{\mathbf{W}}^T_{m+1} \hat{\mathbf{V}}_m $ has the following block structure:
\begin{equation}
	\label{eq:GCRO_block_structure1}
	\hat{\mathbf{W}}^T_{m+1} \hat{\mathbf{V}}_m =
	\begin{bmatrix}
		\mathbf{C}_k^T \tilde{\mathbf{U}}_k & \mathbf{0}_{k \times (m-k)} \\
		\mathbf{V}_{m-k+1}^T \tilde{\mathbf{U}}_k &  \begin{bmatrix} \mathbf{I}_{m-k} \\ \mathbf{0}_{1 \times (m-k)} \end{bmatrix}
	\end{bmatrix} \!\!.
\end{equation}

We can further simplify this expression. Let $\mathbf{v}_{k+1}=\mathbf{r}_k/\| \mathbf{r}_k \|$, we have $\hat{\mathbf{V}}_m = [\tilde{\mathbf{U}}_k, \mathbf{v}_{k+1}, \mathbf{V}_{m-k-1}]$ and $\hat{\mathbf{W}}_{m+1} = [\mathbf{C}_k, \mathbf{v}_{k+1},\mathbf{V}_{m-k}]$. It can also be shown that $\operatorname{range}(\tilde{\mathbf{U}}_k, \mathbf{v}_{k+1}) = \operatorname{range}(\mathbf{C}_k, \mathbf{v}_{k+1}) = \operatorname{range}(\mathbf{Y}_k, \mathbf{v}_{k+1})$, $\mathbf{Y}_k = \hat{\mathbf{V}}_m \mathbf{P}_k$ being the deflation subspace spanned by the Ritz vectors (see \cite{carvalho2011flexible}). It follows that the bottom left $(m-k) \times (k+1)$ block in \eqref{eq:GCRO_block_structure1} is null, which gives
\begin{equation}
\label{eq:GCRO_block_structure2}
\hat{\mathbf{W}}^T_{m+1} \hat{\mathbf{V}}_m =
\left[
	\begin{array}{cc|c}
	\mathbf{C}_k^T \tilde{\mathbf{U}}_k & \mathbf{0}_{k \times 1} & \mathbf{0}_{k \times (m-k-1)} \\
	\mathbf{v}_{1}^T \tilde{\mathbf{U}}_k & 1 & \mathbf{0}_{1 \times (m-k-1)} \\
	\cmidrule[.5pt]{1-3}
	\multicolumn{2}{c|}{\mathbf{V}_{m-k}^T [\tilde{\mathbf{U}}_k, \mathbf{r}_k/\| \mathbf{r}_k \|]}  & \left[ \! \begin{array}{c} \mathbf{I}_{m-k-1} \\ \mathbf{0}_{1 \times (m-k-1)} \end{array} \! \right]
	\end{array}
\right] = \hfill
\left[
	\begin{array}{cc|c}
	\mathbf{C}_k^T \tilde{\mathbf{U}}_k & \mathbf{0}_{k \times 1} & \mathbf{0}_{k \times (m-k-1)} \\
	\mathbf{v}_{1}^T \tilde{\mathbf{U}}_k & 1 & \mathbf{0}_{1 \times (m-k-1)} \\
	\cmidrule[.5pt]{1-3}
	\multicolumn{2}{c|}{\mathbf{0}_{(m-k) \times (k+1)}}  & \left[ \! \begin{array}{c} \mathbf{I}_{m-k-1} \\ \mathbf{0}_{1 \times (m-k-1)} \end{array} \! \right]
	\end{array}
\right] \!\!.
\end{equation}

In \cite{carvalho2011flexible} the authors suggest a number of improvements, mainly based on recursion formulas, to save computational cost in forming the block $\mathbf{C}_k^T \tilde{\mathbf{U}}_k$ in \eqref{eq:GCRO_block_structure2} at a cost independent of the problem size. The deflation subspace spanned by the Ritz vectors $\mathbf{Y}_k = \hat{\mathbf{V}}_m \mathbf{P}_k$ is then shifted by $\mathbf{A}$ to give
\begin{equation}
	\label{eq:shifted_ritz_subspace}
	\mathbf{A} \mathbf{Y}_k = \hat{\mathbf{W}}_{m+1} \bar{\mathbf{H}}_m \mathbf{P}_k.
\end{equation}

Introducing in \eqref{eq:shifted_ritz_subspace} the reduced QR factorization $\bar{\mathbf{H}}_m \mathbf{P}_k = \bar{\mathbf{Q}}_m \mathbf{R}$, $\bar{\mathbf{Q}}_m \in \mathbb{R}^{(m+1) \times k}$, we obtain $\mathbf{A} \mathbf{Y}_k = ( \hat{\mathbf{W}}_{m+1} \bar{\mathbf{Q}}_m) \mathbf{R} = \mathbf{C}_k \mathbf{R}$. The outer residual subspace for the next cycle is then given by 
\begin{equation}
	\label{eq:new_Ck}
	\mathbf{C}_k = \hat{\mathbf{W}}_{m+1} \bar{\mathbf{Q}}_m,
\end{equation}

\noindent and, from $ \mathbf{C}_k = \mathbf{A} \mathbf{Y}_k \mathbf{R}^{-1} = \mathbf{A} \mathbf{U}_k $, it follows directly the definition of the corresponding solution subspace 

\begin{equation}
	\label{eq:new_Uk}
	\mathbf{U}_k = \mathbf{Y}_k \mathbf{R}^{-1}.
\end{equation}

We note that \eqref{eq:new_Ck} implies that $\mathbf{C}_k$ is an orthonormal matrix since  $ \hat{\mathbf{W}}_{m+1} $ and $ \bar{\mathbf{Q}}_m $ are also orthonormal. We are now in a position to perform the next cycle $(i+1)$ and obtain the next iterates $\ve{r}_{i+1}^{(s)}$ and $\ve{x}_{i+1}^{(s)}$ for the current system $(s)$ by performing $m-k$ steps of GMRES applied to the projected operator $\te{A}_{\te{C}_k}=(\matI - \te{C}_{k}\te{C}^{T}_{k})\te{A}$.

We still have to detail the construction of the deflation subspaces $\mathbf{U}_k$ and $\mathbf{C}_k$ at the end of the initial GMRES cycle on one hand, and at the beginning of a subsequent linear system on the other hand.
\subsection*{\underline{First cycle of initial linear system}}
After $m$ iterations of GMRES we generate $\mathbf{V}_{m+1}$ and $\bar{\mathbf{H}}_{m}$ which are related by the Arnoldi relation $\mathbf{AV}_{m} = \mathbf{V}_{m+1}\bar{\mathbf{H}}_{m}$. We then solve the standard eigenvalue problem \eqref{eq:Standard_eigenvalue_problem} to obtain the matrix $\mathbf{P}_{k} \in \mathbb{R}^{m \times k}$  which columns correspond to the retained eigenvectors used for the construction of the harmonic Ritz vectors $\mathbf{Y}_k = \mathbf{V}_m \mathbf{P}_k$. Defining $[\mathbf{Q}, \mathbf{R}]$ as the reduced QR-factorization of $\bar{\mathbf{H}}_{m} \mathbf{P}_k$, we have

\begin{equation}
    \mathbf{AY}_k = \mathbf{AV}_m \mathbf{P}_k = \mathbf{V}_{m+1} \ols{\mathbf{H}}_m \mathbf{P}_k = (\mathbf{V}_{m+1}\mathbf{Q})\mathbf{R}.
    \label{eq:AVmPk_CkR}
\end{equation}

Since $\mathbf{V}_{m+1}$ and $\mathbf{Q}$ are orthonormal matrices and we are looking for a relation of the type  \eqref{eq8}, we directly deduce from \eqref{eq:AVmPk_CkR} the following identities to be used for the next cycle:
\begin{align}
    \mathbf{C}_k &= \mathbf{V}_{m+1} \mathbf{Q} \\
    \mathbf{U}_k &= \mathbf{Y}_k \mathbf{R}^{-1}
\end{align}

At the end of the GMRES($m$) cycle, the optimality property $\mathbf{r}_1^{(1)} \; \bot \; \mathcal{AK}_m(\mathbf{A},\mathbf{r}_0^{(1)})$ of the residual $\mathbf{r}_1^{(1)}=\mathbf{r}_0^{(1)}-\mathbf{AV}_m \mathbf{y}_m$ can be easily verified. Recall that the reduced coordinates $\mathbf{y}_m$ solve the least-squares problem $ \mathbf{y}_m = \underset{\mathbf{y} \in \mathbb{R}^m}{\operatorname{argmin}} \|\mathbf{V}_{m+1}^T \mathbf{r}_0^{(1)} - \bar{\mathbf{H}}_{m}\mathbf{y}\|_2 $, which gives $ \mathbf{y}_m = ( \bar{\mathbf{H}}_m^T \bar{\mathbf{H}}_m )^{-1} \bar{\mathbf{H}}_m^T \mathbf{V}_{m+1}^T \mathbf{r}_0^{(1)} = \bar{\mathbf{H}}_m^{\dagger} \mathbf{V}_{m+1}^T \mathbf{r}_0^{(1)} = \beta \bar{\mathbf{H}}_m^{\dagger}  \mathbf{e}_1$, with $\beta=\|\mathbf{r}_0\|$. Then, it follows that
\begin{align}
    (\mathbf{AV}_m)^{T} \mathbf{r}^{(1)}_1 &= (\mathbf{AV}_m)^{T} \mathbf{V}_{m+1} ( \beta \mathbf{e}_1 - \bar{\mathbf{H}}_m \mathbf{y}_m ) \nonumber \\
    &= \bar{\mathbf{H}}_m^T \mathbf{V}_{m+1}^T \mathbf{V}_{m+1} ( \beta \mathbf{e}_1 - \bar{\mathbf{H}}_m \mathbf{y}_m ) \nonumber \\
    &= \bar{\mathbf{H}}_m^T ( \beta \mathbf{e}_1 - \beta \bar{\mathbf{H}}_m  \bar{\mathbf{H}}_m^{\dagger}  \mathbf{e}_1 ) \nonumber \\
    \Rightarrow (\mathbf{AV}_m)^{T} \mathbf{r}^{(1)}_1 &= \mathbf{0}.
    \label{eq:optimality_GMRES1}
\end{align}
\indent Equation \eqref{eq:optimality_GMRES1} can be written for any system $(s)$, i.e., $(\mathbf{AV}_m)^{T} \mathbf{r}^{(s)}_1 = \mathbf{0}, s=1,2,\cdots$. Now, starting from  \eqref{eq:optimality_GMRES1} we obtain
\begin{align}
    \mathbf{P}_k^T (\mathbf{AV}_m)^{T} \mathbf{r}^{(1)}_1 &= \mathbf{0} \nonumber \\
    ( \mathbf{V}_{m+1} \bar{\mathbf{H}}_m \mathbf{P}_k )^T \mathbf{r}^{(1)}_1 &= \mathbf{0} \nonumber \\
    \Rightarrow \mathbf{R}^T \mathbf{C}_k^T \mathbf{r}^{(1)}_1 &= \mathbf{0},
    \label{eq:optimality_GMRES2}
\end{align}

\noindent which shows that $ \mathbf{C}_k^T \mathbf{r}^{(1)}_1 = \mathbf{0} $ since $\mathbf{R}$ is supposed to be nonsingular. Similarly, this is valid for any system in the sequence: $ \mathbf{C}_k^T \mathbf{r}^{(s)}_1 = \mathbf{0}, s=1,2,\cdots $. Incidentally, comparing \eqref{eq:optimality_GMRES1} and \eqref{eq:optimality_GMRES2} shows that $\mathcal{AK}_m(\mathbf{A},\mathbf{r}_0^{(1)})$ is spanned by the columns of $\mathbf{C}_k$. 

\subsection*{\underline{First cycle of subsequent linear systems}}
We now assume that two subspaces $\mathbf{U}_k^{(s-1)}$ and $\mathbf{C}_k^{(s-1)}$ are available from a previous linear system $(s-1)$ and we want to solve a new system $ \mathbf{A} \mathbf{x}^{(s)} = \mathbf{b}^{(s)} $. In this work, we deal with the special case where the system matrix remains unchanged from one linear system to the next. This implies $\mathbf{U}_k^{(s)} = \mathbf{U}_k^{(s-1)}$  and $\mathbf{C}_k^{(s)} = \mathbf{C}_k^{(s-1)}$.
In order to apply $m-k$ steps of GMRES to the projected operator $\te{A}_{\te{C}_k^{(s)}}=(\matI - \mathbf{C}_k^{(s)} \mathbf{C}_k^{(s)T})\te{A}$, we need to compute first the projected residual $\mathbf{r}_k^{(s)} = (\mathbf{I} - \mathbf{C}_k^{(s)} \mathbf{C}_k^{(s)T}) \mathbf{r}_0^{(s)} $ and the associated solution vector $\mathbf{x}_k^{(s)} = \mathbf{x}_0^{(s)} + \mathbf{U}_k^{(s)} \mathbf{C}_k^{(s)T}\mathbf{r}_0^{(s)} $, where $\mathbf{x}_0^{(s)}$ and $\mathbf{r}_0^{(s)}$ denote the solution and residual vector at convergence of the previous linear system. Once the Arnoldi process \eqref{eq:GeneralizedArnoldi} is completed, the solution at the end of the first cycle of system $(s)$ reads (see \eqref{eq:GCRO_residual} and \eqref{eq:GCRO_solution}):
\begin{align}
	\mathbf{r}_{1}^{(s)} &= \mathbf{r}_{k}^{(s)} - \mathbf{V}_{m-k+1}\bar{\mathbf{H}}_{m-k}\mathbf{y}_{m-k} \\
	\mathbf{x}_{1}^{(s)} &= \mathbf{x}_{k}^{(s)} + \mathbf{A}^{-1}\mathbf{A}_{\mathbf{C}_{k}^{(s)}}\mathbf{V}_{m-k}\mathbf{y}_{m-k} \nonumber \\
	& = \mathbf{x}_{k}^{(s)} + (\matI-\mathbf{U}_k^{(s)}\mathbf{C}^{(s)T}_{k} \mathbf{A}) \mathbf{V}_{m-k} \mathbf{y}_{m-k} \nonumber \\
	& = \mathbf{x}_{k}^{(s)} + (\mathbf{V}_{m-k}-\mathbf{U}_k^{(s)}\mathbf{B}_{m-k})  \mathbf{y}_{m-k}.
\end{align}

The above expressions make explicitly appear the correction to the outer iterates $\mathbf{x}_k^{(s)}$ and $\mathbf{r}_k^{(s)}$. Equivalent formulas can be obtained by considering the generalized Arnoldi relation \eqref{eq:GeneralizedArnoldi}:
\begin{align}
	\mathbf{r}_{1}^{(s)} &= \mathbf{r}_{0}^{(s)} - \hat{\mathbf{W}}_{m+1} \bar{\mathbf{H}}_m\mathbf{y}_{m} \\
	\mathbf{x}_{1}^{(s)} &= \mathbf{x}_{0}^{(s)} + \hat{\mathbf{V}}_{m} \mathbf{y}_{m}
\end{align}
The vector of reduced coordinates $\mathbf{y}_m$ is the solution of \eqref{GCRO_LeastSquares} and its partition $\mathbf{y}_{m-k}$ is the solution of ${\operatorname{argmin}}_{\mathbf{y}} \left\| \bar{\mathbf{H}}_{m-k} \mathbf{y} -\beta \mathbf{e}_1^{m-k+1} \right\|_2$, $\beta = \| (\mathbf{I}-\mathbf{C}_k^{(s)}\mathbf{C}_k^{(s)T})\mathbf{r}_0^{(s)} \|_2$.

\subsection*{\underline{Reformulation of the generalized eigenvalue problem}}
At the end of each cycle of GCRO-DR we compute the Ritz eigenpairs from the following generalized eigenvalue problem:
\begin{equation}
	\bar{\mathbf{H}}^T_m \bar{\mathbf{H}}_m \mathbf{g} = \theta \bar{\mathbf{H}}^T_m \hat{\mathbf{W}}^T_{m+1} \hat{\mathbf{V}}_m \mathbf{g}
	\label{GCRODR_Ritz_eigenproblem}
\end{equation}

However, direct construction of both sides of this equation as product of matrices may lead to highly ill-conditioned operators. In the following we propose to reformulate \eqref{GCRODR_Ritz_eigenproblem} by exploiting the specific structure of $\bar{\mathbf{H}}_m$ and $\hat{\mathbf{W}}^T_{m+1} \hat{\mathbf{V}}_m$. First, noticing that the last row of $\bar{\mathbf{H}}_m$ is null except the last term $h_{m+1,m}$, we can write $\bar{\mathbf{H}}_m^T = [ \mathbf{H}_m^T \quad h_{m+1,m}\mathbf{e}_m ] $. Thus, we have $\bar{\mathbf{H}}^T_m \bar{\mathbf{H}}_m = \mathbf{H}^T_m \mathbf{H}_m + h_{m+1,m}^2\mathbf{e}_m \mathbf{e}_m^T$. Inserting these relations in \eqref{GCRODR_Ritz_eigenproblem} gives
\begin{equation}
	[\mathbf{H}^T_m \mathbf{H}_m + h_{m+1,m}^2\mathbf{e}_m \mathbf{e}_m^T ] \mathbf{g} = \theta [ \mathbf{H}_m^T \quad h_{m+1,m}\mathbf{e}_m ] \hat{\mathbf{W}}^T_{m+1} \hat{\mathbf{V}}_m \mathbf{g}.
\end{equation}

Left multiplying by $\mathbf{H}_m^{-T}$ and defining $\mathbf{f} = \mathbf{H}_m^{-T} \mathbf{e}_m$, we obtain
\begin{equation}
	[ \mathbf{H}_m + h_{m+1,m}^2\mathbf{f} \, \mathbf{e}_m^T ] \mathbf{g} = \theta [ \mathbf{I}_m \quad h_{m+1,m}\mathbf{f} ] \hat{\mathbf{W}}^T_{m+1} \hat{\mathbf{V}}_m \mathbf{g}.
\end{equation}
\indent From \eqref{eq:GCRO_block_structure2} we know that the last row of $\hat{\mathbf{W}}_{m+1}^T \hat{\mathbf{V}}_{m}$ is null which finally leads to the following eigenvalue problem

\begin{equation}
	\label{eq:Reformulated_GCRO_eigenproblem}
	\boxed{[ \mathbf{H}_m + h_{m+1,m}^2 \mathbf{f} \, \mathbf{e}_m^T ] \mathbf{g} = \theta 	
		\left[
		\begin{array}{cc}
			\hat{\mathbf{W}}_{k+1}^T \hat{\mathbf{V}}_{k+1} & \mathbf{0}_{(k+1) \times (m-k-1)} \\
			\mathbf{0}_{(m-k-1) \times (k+1)}  & \mathbf{I}_{m-k-1}
		\end{array} \right] \mathbf{g}},
\end{equation}

\noindent where the structure of the block $\hat{\mathbf{W}}_{k+1}^T \hat{\mathbf{V}}_{k+1}$ has been highlighted in \eqref{eq:GCRO_block_structure2}.

\begin{algorithm}[H]
    \label{GCRODR_algorithm}
    \caption{GCRO-DR($m$,$k$) for a sequence of linear systems $\mathbf{A}\mathbf{x}^{(s)}=\mathbf{b}^{(s)}$.}
    \fontsize{10}{14}\selectfont
    Choose $m$, the size of the Krylov subspace and $k$, the size of the deflated/recycled subspace. Let \textit{tol} be the relative convergence threshold. Choose an initial guess $\mathbf{x}_{0}$. Compute $\mathbf{r}_{0}=\mathbf{b}-\mathbf{Ax}_{0}$, and set $i=1$. \;
    \eIf{ $\mathbf{C}_{k}$ and $\mathbf{U}_{k}$ are available (from a previous linear system)}
    {
   	  $\mathbf{x}_{1} = \mathbf{x}_{0} + \mathbf{U}_{k}\mathbf{C}^{T}_{k}\mathbf{r}_{0}$  \Comment*[f]{Do not update $\mathbf{U}_k$ and $\mathbf{C}_k$ since only $\mathbf{b}$ varies}  \;
   	  $\mathbf{r}_{1} = \mathbf{r}_{0} - \mathbf{C}_{k}\mathbf{C}^{T}_{k}\mathbf{r}_{0}$ \;
    }{
      $\mathbf{v}_{1} = \mathbf{r}_{0}/\|\mathbf{r}_{0}\|$  \;
      $\mathbf{c} = \|\mathbf{r}_{0}\| \mathbf{e}_{1}$  \;
      
      Apply $m$ steps of GMRES to generate $\mathbf{V}_{m+1}$ and $\bar{\mathbf{H}}_{m}$.  \;
      
      solve $\operatorname{min}_{\mathbf{y}} \|\mathbf{c} -\ols{\mathbf{H}}_{m}\mathbf{y}\|_{2}$ for $\mathbf{y}_m$ \;
      $\mathbf{x}_{1} = \mathbf{x}_{0} + \mathbf{V}_{m}\mathbf{y}_{m}$  \;
      $\mathbf{r}_{1} = \mathbf{V}_{m+1}(\mathbf{c}-\bar{\mathbf{H}}_{m}\mathbf{y}_{m})$  \;
      Select the $k$ deflation eigenvectors of $\mathbf{H}_{m}\mathbf{g} = \lambda \mathbf{g}$ associated to the smallest eigenvalues and concatenate into $\mathbf{P}_k$\;
      Let [Q,R] be the reduced QR factorization of $\bar{\mathbf{H}}_{m}\mathbf{P}_{k}$.  \;
      $\mathbf{C}_{k} = \mathbf{V}_{m+1}\mathbf{Q}$  \;
      $\mathbf{U}_{k} = \hat{\mathbf{V}}_{m}\mathbf{P}_{k}\mathbf{R}^{-1}$  \;
    }
    \While{$\|\mathbf{r}_{i}\|_{2} > \textit{tol}$}  {
      $i = i+1$  \;
      Apply ($m-k$) steps of Arnoldi with the projected operator $(\mathbf{I}- \mathbf{C}_{k}\mathbf{C}^{T}_{k})\mathbf{A}$ starting with $\mathbf{v}_{1}=\mathbf{r}_{i-1}/\|\mathbf{r}_{i-1}\|_{2}$ to generate $\mathbf{V}_{m-k+1}$, $\bar{\mathbf{H}}_{m-k}$, and $\mathbf{B}_{m-k}$.  \;
      $\hat{\mathbf{W}}_{m+1} = [\mathbf{C}_{k} \quad \mathbf{V}_{m-k+1}]$  \;
      $\hat{\mathbf{V}}_m = [\tilde{\mathbf{U}}_k \quad \mathbf{V}_{m-k}]$  \;
      $\bar{\mathbf{H}}_{m} = 
     \begin{bmatrix}
        \mathbf{D}_{k} & \mathbf{B}_{m-k} \\
        \mathbf{0}     & \bar{\mathbf{H}}_{m-k} 
     \end{bmatrix}$  \;

     $ \mathbf{y}_m = {\operatorname{argmin}}_{\mathbf{y}} \|\bar{\mathbf{H}}_{m}\mathbf{y}-\hat{\mathbf{W}}_{m+1}^T \mathbf{r}_{i-1}\|_2 $ \;
     $\mathbf{x}_{i} = \mathbf{x}_{i-1} + \hat{\mathbf{V}}_m \mathbf{y}_m$  \;
     $\mathbf{r}_{i} = \mathbf{r}_{i-1} - \hat{\mathbf{W}}_{m+1}\bar{\mathbf{H}}_m\mathbf{y}_{m}$  \;
     
	Solve $\mathbf{f}=\mathbf{H}^{-T}_{m}\mathbf{e}_{m}$  \;     
	Select the $k$ deflated eigenvectors of $[ \mathbf{H}_m + h_{m+1,m}^2 \mathbf{f} \, \mathbf{e}_m^T ] \mathbf{g} = \theta 	
	\begin{bmatrix}
		\hat{\mathbf{W}}_{k+1}^T \hat{\mathbf{V}}_{k+1} & \mathbf{0}_{(k+1) \times (m-k-1)} \\
		\mathbf{0}_{(m-k-1) \times (k+1)}  & \mathbf{I}_{m-k-1} 
	\end{bmatrix} \mathbf{g}$ associated to the smallest eigenvalues and concatenate into $\mathbf{P}_k$\;
     
     $\mathbf{Y}_k = \hat{\mathbf{V}}_{m} \mathbf{P}_{k} $ \;
     Let [Q,R] be the reduced QR factorization of $\bar{\mathbf{H}}_{m}\mathbf{P}_{k}$.  \;
     $\mathbf{C}_k = \hat{\mathbf{V}}_{m+1} \mathbf{Q}$  \;
     $\mathbf{U}_{k} = \mathbf{Y}_k\mathbf{R}^{-1}$  \;
    }
    Recycle $\mathbf{U}_{k}$ and $\mathbf{C}_{k}$ for the next linear system  \;
\end{algorithm}   
%
\subsection{Application of GCRO-DR to the solution of the fluid adjoint system}

\begin{figure}[H]
	\captionsetup[subfigure]{belowskip=-6pt}
	\begin{subfigure}{.5\textwidth}
		\raggedright
		\includegraphics[trim={0.5cm 0.5cm 1.5cm 1.9cm},clip,width=0.99\textwidth]{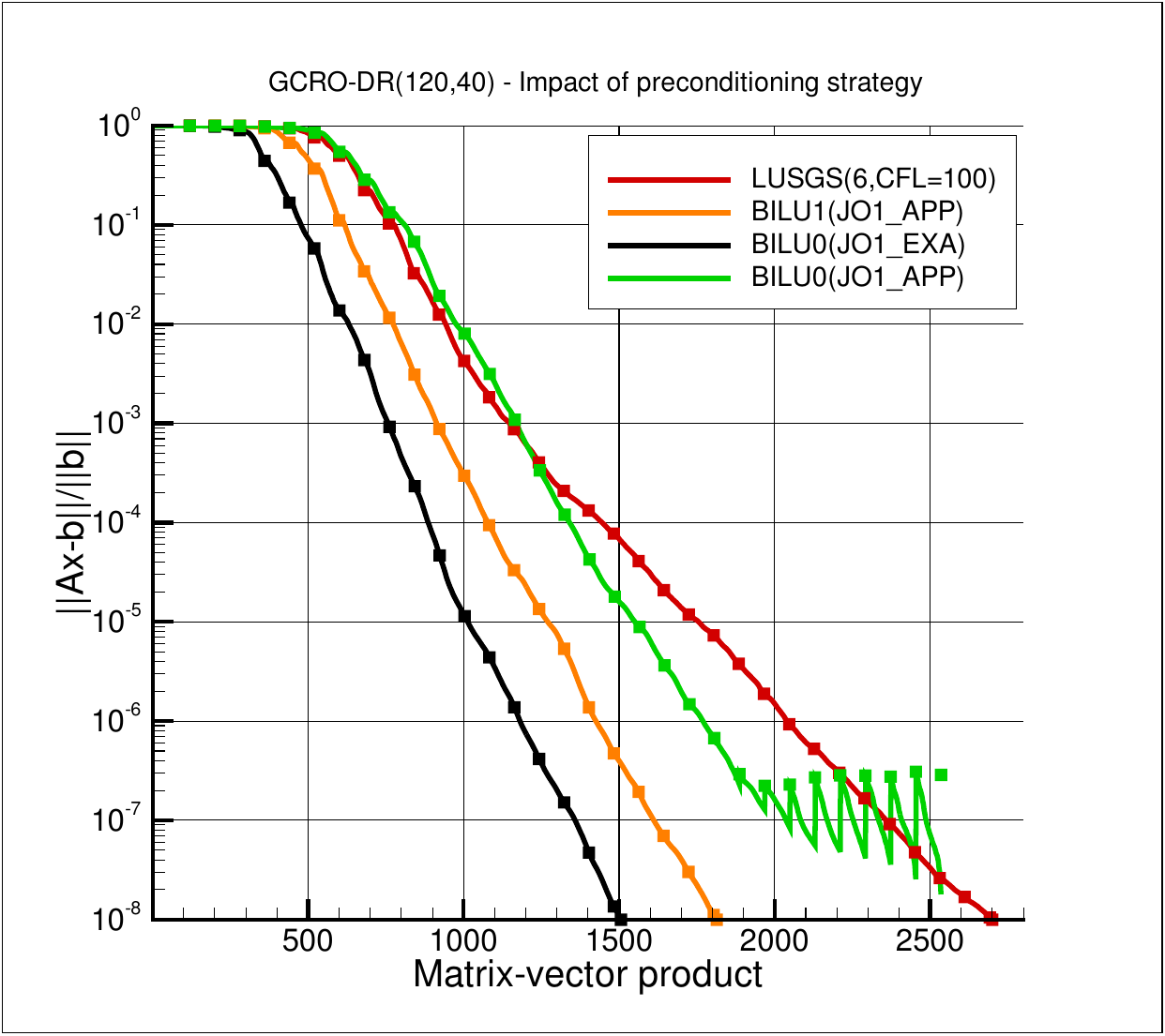}
		\caption{}
		\label{fig:cvg_GCRODR_120_40_impact_preconditioning_strategy_matvecprod_cvg_erratic}
	\end{subfigure}
	\begin{subfigure}{.5\textwidth}
		\centering
		\includegraphics[trim={0.5cm 0.5cm 1.5cm 1.9cm},clip,width=0.99\textwidth]{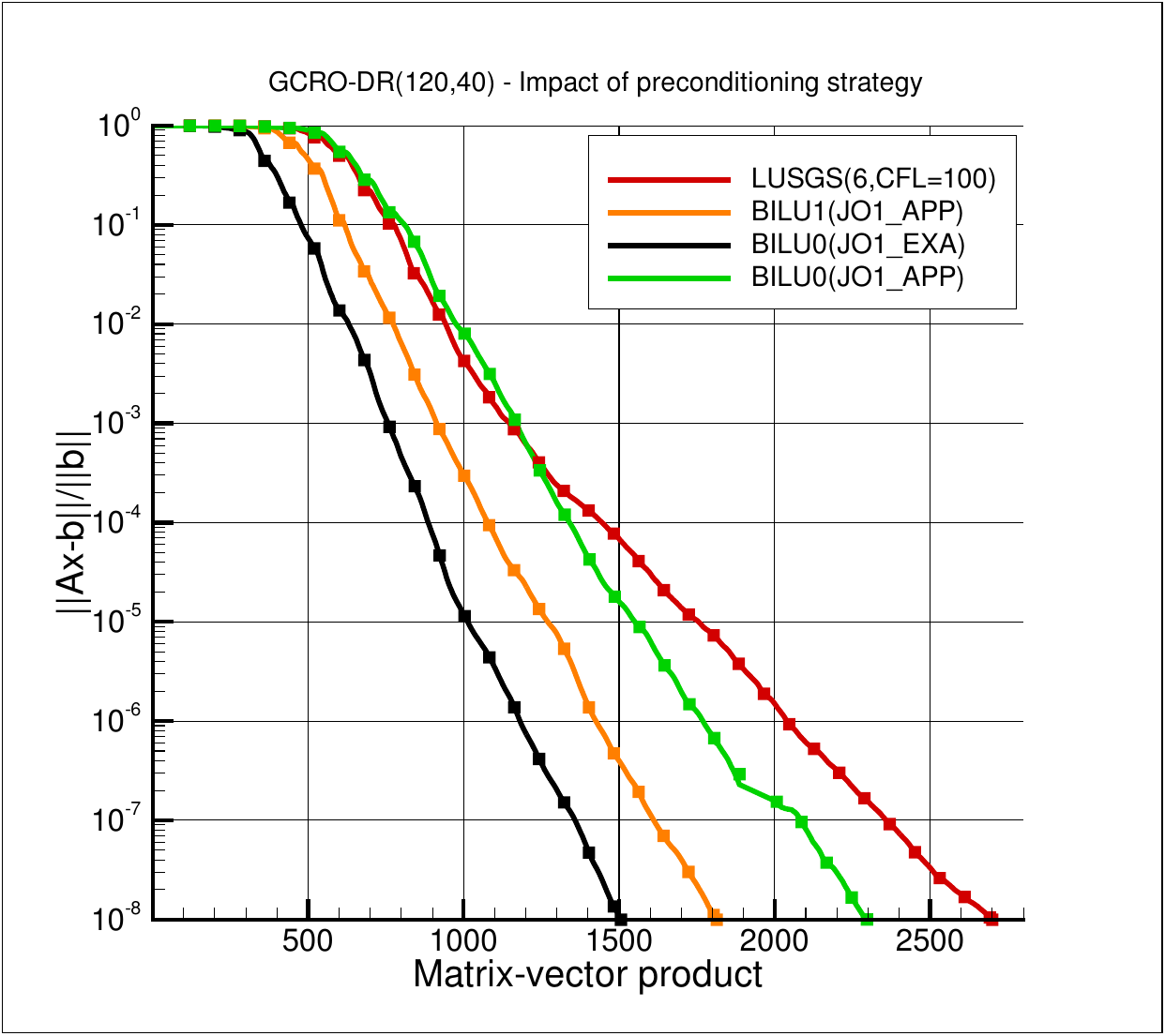}
		\caption{}		
		\label{fig:cvg_GCRODR_120_40_impact_preconditioning_strategy_matvecprod_cvg_resratio}
	\end{subfigure}
	\vspace{0pt}
	\caption{Adjoint relative residual norm convergence history of GGCRO-DR(120,40). Impact of various preconditioners. In \ref{fig:cvg_GCRODR_120_40_impact_preconditioning_strategy_matvecprod_cvg_erratic} an erratic convergence of the least-squares residual associated to a stagnation of the true residual occurs for the BILU(0) factorization applied to the first-order approximate Jacobian matrix.  In \ref{fig:cvg_GCRODR_120_40_impact_preconditioning_strategy_matvecprod_cvg_resratio} a cold restart allows to restore the convergence.}
	\label{fig:cvg_GCRODR_120_40_impact_preconditioning_strategy_matvecprod}
\end{figure}

\begin{figure}[H]
	\centering
	\includegraphics[trim=0.25cm 0.5cm 1.5cm 1.5cm,clip,width=0.6\textwidth]{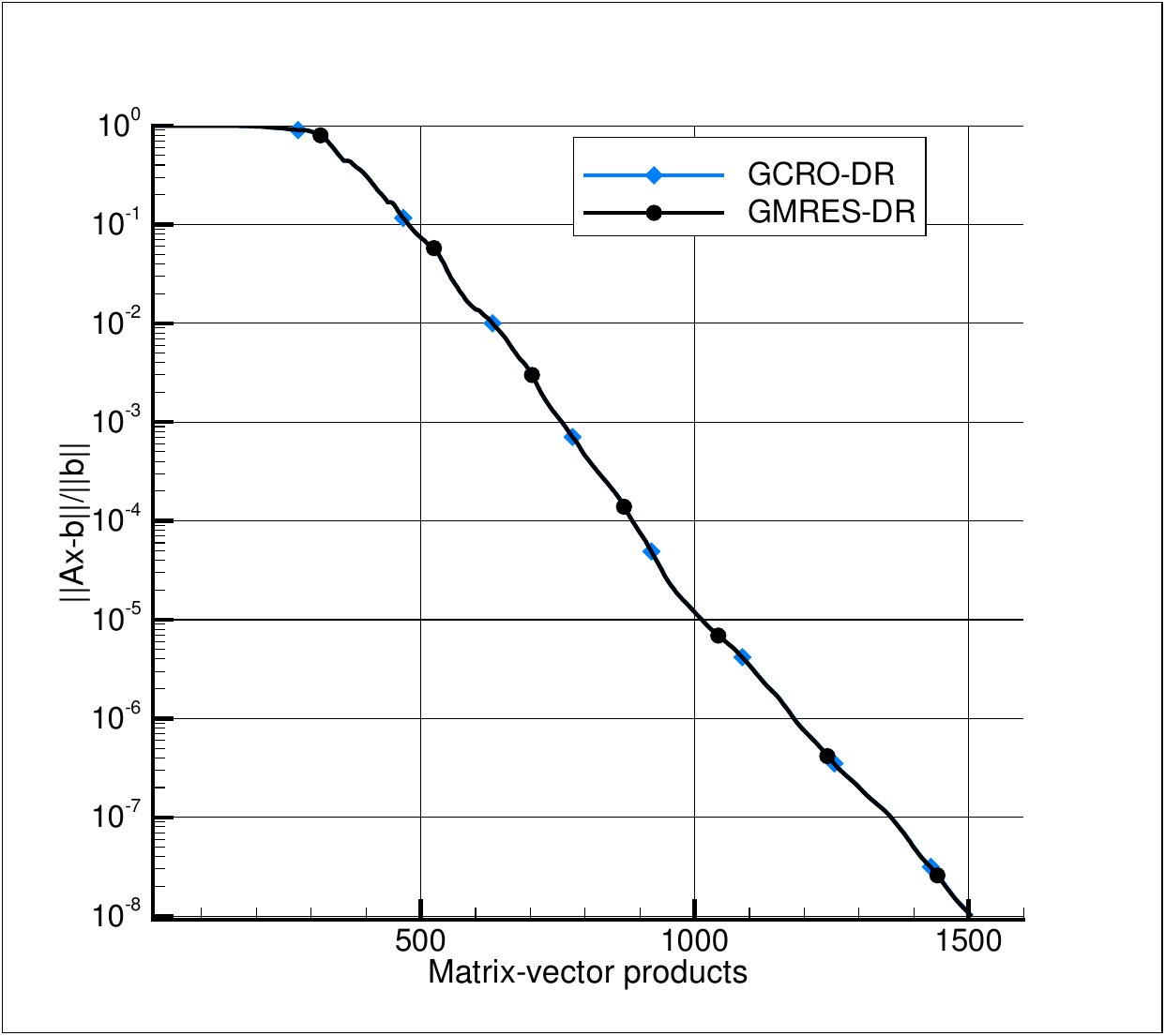}
	\caption{Illustration of the algebraic and numerical equivalence of GMRES-DR and GCRO-DR for the preconditioner BILU(0) applied to the exact first-order Jacobian matrix. We recall the solver parameters: $m=120$ and $k=40$.}
	\label{fig:cvg_OPT_comparison_GMRESDR_GCRODR}
\end{figure}

As for GMRES-DR in section \ref{GMRESDR_Fluid_Adjoint_System}, we first apply GCRO-DR(120,40) to the purely fluid adjoint system. The convergence curves are plotted in Figure \ref{fig:cvg_GCRODR_120_40_impact_preconditioning_strategy_matvecprod}. This Figure should be compared to the corresponding Figure \ref{fig:cvg_GMRESDR_120_40_impact_preconditioning_strategy_matvecprod_cvg_erratic} for GMRES-DR. First, we note that GCRO-DR performs much better than GMRES-DR in terms of numerical robustness since the true and least-squares relative residuals match perfectly except for the preconditioner BILU(0) applied to the approximate first-order Jacobian matrix. However, the discrepancy only appears below a threshold of $5\times10^{-7}$. Again, this loss of accuracy can be suppressed by resorting to a cold restart (see Figure \ref{fig:cvg_GCRODR_120_40_impact_preconditioning_strategy_matvecprod_cvg_resratio}). In theory, GMRES-DR and GCRO-DR are algebraically equivalent \cite{carvalho2011flexible} but in practice our numerical experiments show that this equivalence is only observed for relative residuals above $1\times10^{-5}$, except for the most accurate preconditioner BILU(0) applied to the exact first-order Jacobian matrix. This perfect numerical equivalence is illustrated in Figure \ref{fig:cvg_OPT_comparison_GMRESDR_GCRODR}.
\vspace{6pt}

\subsection{Application of GCRO-DR to the solution of the fluid-structure coupled-adjoint system}

\begin{figure}[H]
	\centering
	\includegraphics[trim={1.cm 0.5cm 1.5cm 1.95cm},clip,width=0.63\textwidth]{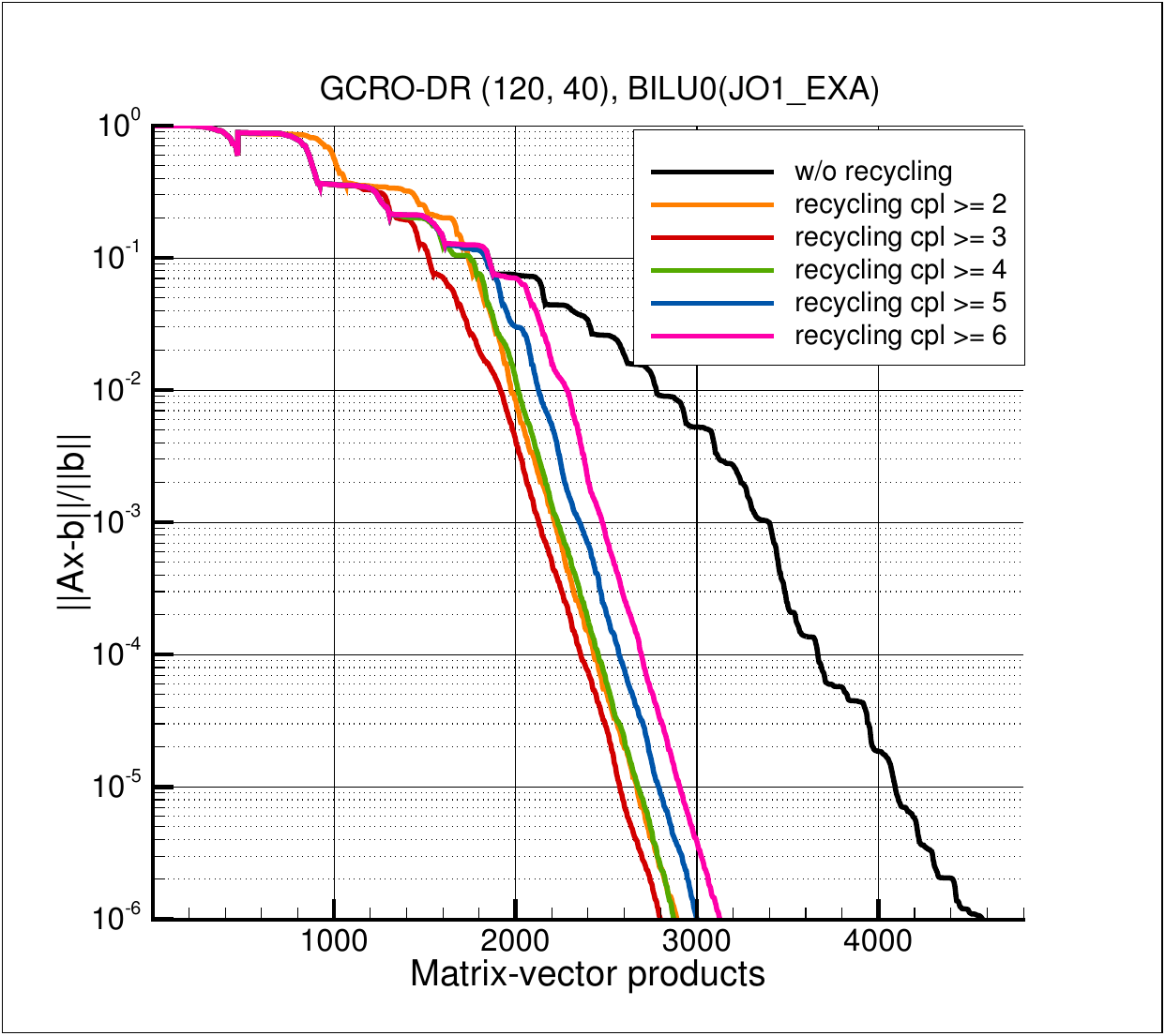}
	\caption{Coupled-adjoint relative residual norm convergence history of GCRO-DR(120,40). Impact of recycling strategy. Recycling spectral information starting from cycle 2 and above always improves convergence. A maximum saving of 1780 matrix-vector products out of 4580 ($\sim$39\%) is achieved.}
	\label{fig:parameter_study_GCRODR_subspace_recycling_matvecprod}
\end{figure}

As exposed in section \ref{GCRODR_theory}, unlike GMRES-DR which is based on an augmentation deflation strategy, GCRO-DR is able to accelerate the solution of a sequence of linear systems with a slowly varying system matrix and/or right-hand side. The strategy relies on recycling an approximate invariant subspace from one system to the next. In the context of coupled-adjoint linear systems, only the right-hand side of the fluid block of equations changes during the partitioned solution strategy. In Figure \ref{fig:parameter_study_GCRODR_subspace_recycling_matvecprod} we illustrate the acceleration of convergence for different recycling strategies. The black plain line corresponds to the convergence of the standard GMRES-DR or GCRO-DR solver preconditioned by BILU(0) applied to the first-order exact Jacobian matrix. For each new fluid-structure cycle a cold restart is performed, discarding any previous spectral information. The corresponding plateaus are clearly observable at the beginning of each cycle. The remaining curves present the convergence for a subspace recycling activated from an increasing fluid-structure cycle index. As can be seen, even recycling starting from the second cycle improves convergence compared to the standard Krylov solver. The best option seems here to recycle from the second fluid-structure coupling. This is very promising because one would have expected that the quality of the spectral information to be recycled, i.e., the distance between the true invariant subspace and its approximation, would not be small enough to get such an acceleration. A maximum saving of 1780 matrix-vector products out of 4580 ($\sim$39\%) is achieved.

\begin{figure}[H]
	\captionsetup[subfigure]{belowskip=-6pt}
	\begin{subfigure}{.5\textwidth}
		\raggedleft
		\includegraphics[trim={0.25cm 0.5cm 0.8cm 1.8cm},clip,width=0.99\textwidth]{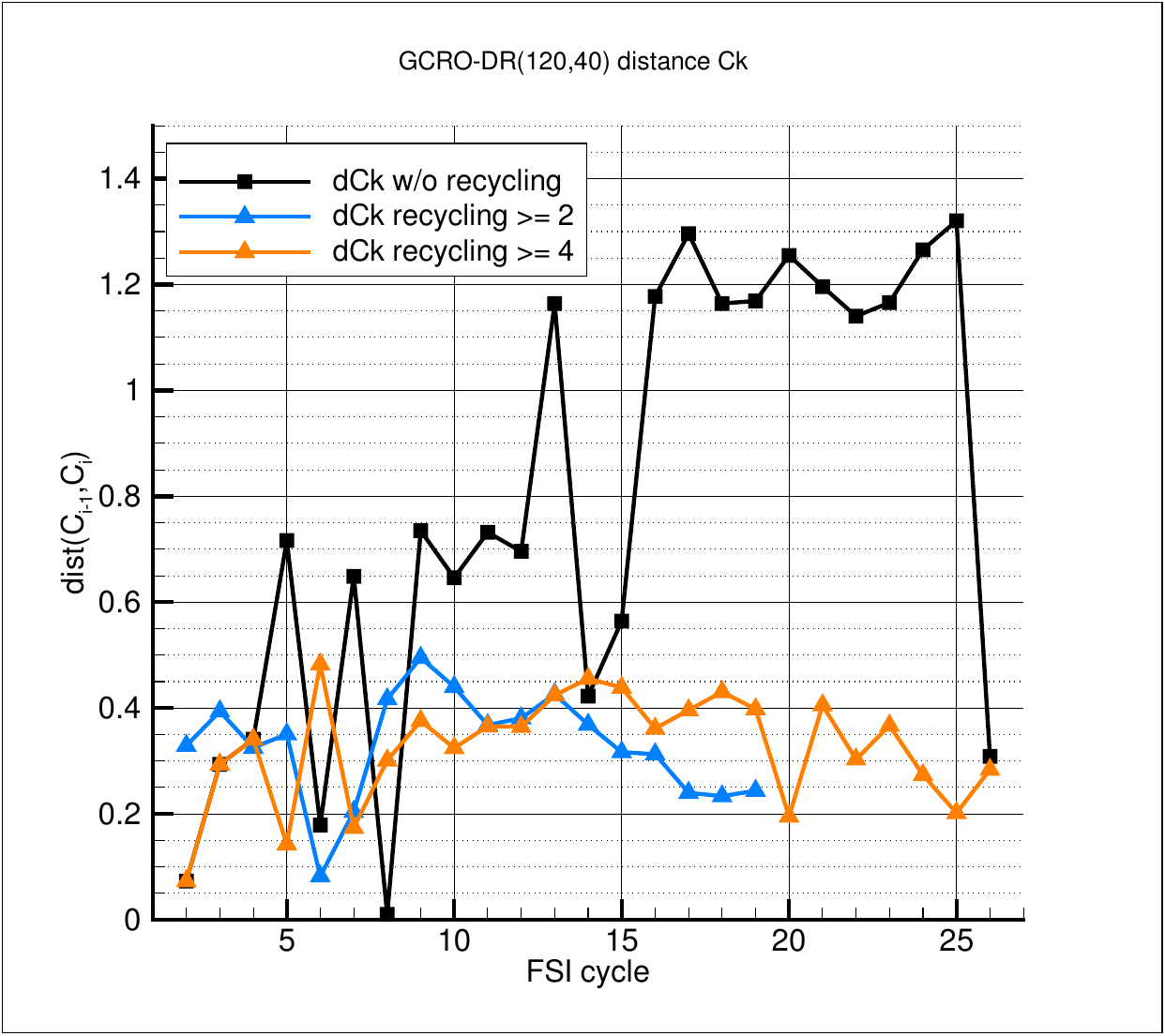}
		\caption{}
		\label{fig:GCRO-DR_120_40_distances}
	\end{subfigure}
	\begin{subfigure}{.5\textwidth}
		\raggedleft
		\includegraphics[trim={0.25cm 0.5cm 0.8cm 1.8cm},clip,width=0.99\textwidth]{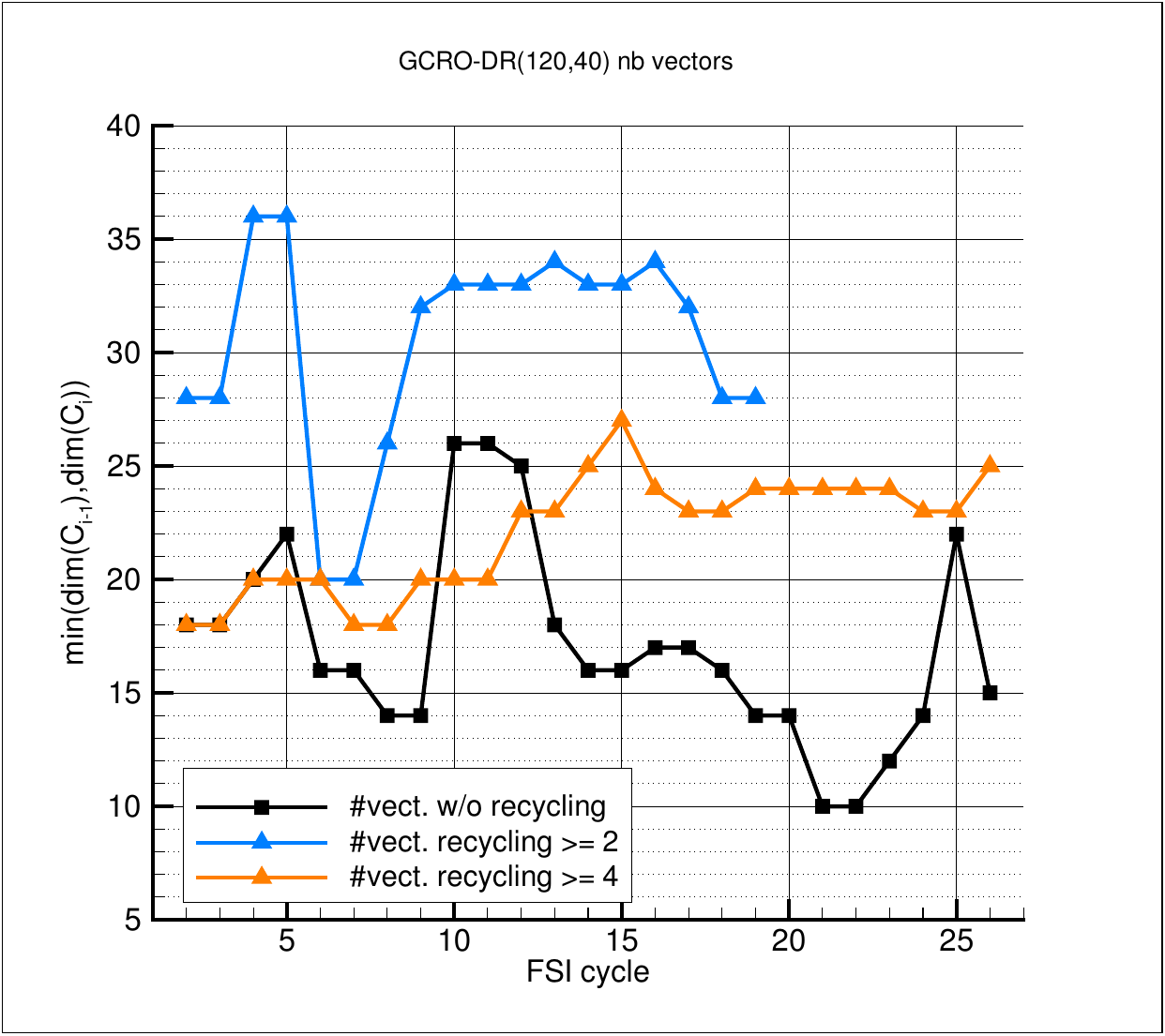}
		\caption{}		
		\label{fig:GCRO-DR_120_40_Ck_vectors}
	\end{subfigure}
	\vspace{0pt}
	\caption{Impact of subspace recycling for GCRO-DR(120,40). Monitoring of distance between approximate invariant subspaces $\mathbf{C}^{(i-1)}$ and $\mathbf{C}^{(i)}$ is reported in \ref{fig:GCRO-DR_120_40_distances}.  The number of vectors used to compute the subspace distance is also reported in \ref{fig:GCRO-DR_120_40_Ck_vectors}.}
	\label{fig:GCRO-DR_120_40_subspace_recycling_monitoring}
\end{figure}

In an attempt to understand the convergence of the approximate invariant subspace, we propose to compute the distance between $\mathbf{C}^{(i-1)}$ and $\mathbf{C}^{(i)}$ which is the best estimation of the true distance to the invariant subspace of the system matrix that we can compute. Since the size of the recycling subspace varies (in this case 30\% of the size of the Krylov space for the last fluid cycle) we use the Grassmann distance formula between vector spaces of different dimensions \cite{ye2016schubert}:

\begin{equation}
	d_{p}(\mathbf{C}^{(i-1)}, \mathbf{C}^{(i)}) = \left( \sum^{p}_{i=1} \theta^{2}_i \right)^{1/2}, \quad p = \textrm{min}(k_{i-1}, k_{i}),
	\label{eq:Grassmann_distance}
\end{equation} 

\noindent where $\theta_i$, the principal angles between columns of $\mathbf{C}^{(i-1)}$ and $\mathbf{C}^{(i)}$, are computed via the singular value decomposition of ${\mathbf{C}^{(i-1)}}^T \mathbf{C}^{(i)}$. In order to interpret this distance criterion during convergence, we plot a normalized version of \eqref{eq:Grassmann_distance}: $\tilde{d}_{p}=d_{p}/\sqrt{p}$.
In Figures \ref{fig:GCRO-DR_120_40_distances} and \ref{fig:GCRO-DR_120_40_Ck_vectors} we report the histories of the subspace distance $d_p$ and of the smaller subspace dimension $p$. In case of no recycling, the distances show higher values meaning that the Krylov solver continues to search after an optimal projection space up to convergence. When recycling is activated, the last distances are about three times lower which confirms that spectral information was actually propagated between cycles.

\begin{figure}[H]
	\centering
	\includegraphics[trim={1.cm 0.5cm 0.2cm 1.95cm},clip,width=0.65\textwidth]{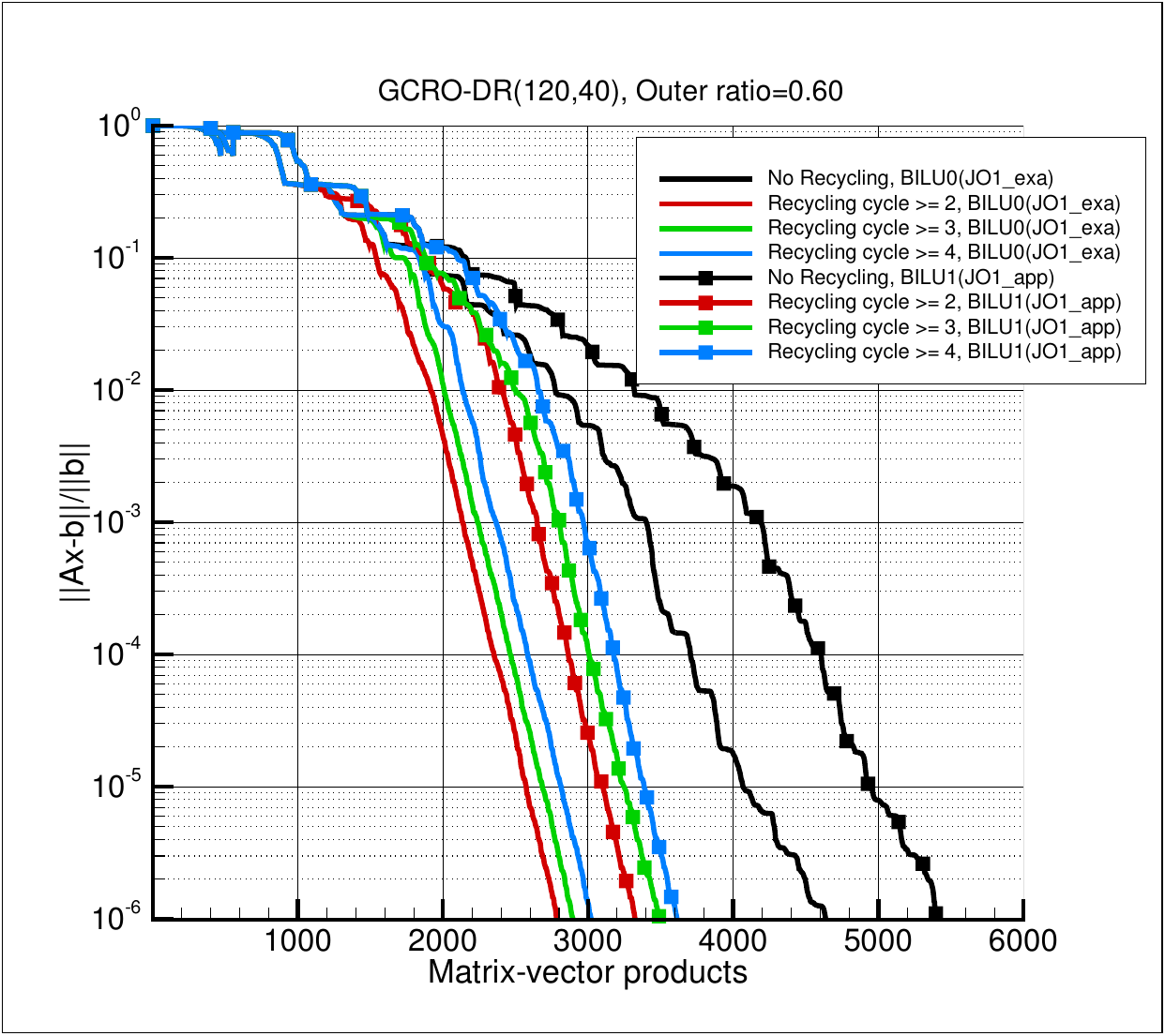}
	\caption{Coupled-adjoint relative residual norm convergence history of GCRO-DR(120,40). Impact of recycling strategy for two BILU preconditioners.}
	\label{fig:parameter_study_GCRODR_subspace_recycling_matvecprod_BILU0_vs_BILU1}
\end{figure}

To conclude our numerical experiments, Figure \ref{fig:parameter_study_GCRODR_subspace_recycling_matvecprod_BILU0_vs_BILU1} compares two instances of GCRO-DR with subspace recycling combined to the BILU(1) preconditioner applied to the first-order approximate Jacobian operator on one hand, and to BILU(0) applied to the first-order exact flux Jacobian operator on the other hand. Even for the simpler preconditioner a remarkable acceleration is observed. 

\section{Flexible Generalized Conjugate Residual with inner Orthogonalization and Deflated Restarting: FGCRO-DR} \label{FGCRODR_theory}

\subsection{The FGCRO-DR Krylov solver}

Similar to GCRO-DR, FGCRO-DR relies on a given full-rank matrix $\mathbf{Z}_{k} \in \mathbb{R}^{n \times k}$ and an orthonormal matrix $\mathbf{C}_{k} \in \mathbb{R}^{n \times k}$ as the image of $\mathbf{Z}_{k}$ by $\mathbf{A}$ satisfying the relations
\begin{align}
	\mathbf{A}\mathbf{Z}_{k}         &= \mathbf{C}_{k},  \\
	\mathbf{C}^{T}_{k}\mathbf{C}_{k} &= \mathbf{I}_{k}.
\end{align}

The first cycle of FGCRO-DR for solving the initial linear system in the sequence consists in applying $m$ steps of the flexible Arnoldi process to build $\mathbf{V}_{m+1}$, $\mathbf{Z}_{m}$ and $\Bar{\mathbf{H}}_{m}$ (see algorithm \ref{alg:FGMRES}). At this point, we assume that the following flexible Arnoldi relation exists: $\mathbf{AZ}_{m} = \mathbf{V}_{m+1}\bar{\mathbf{H}}_{m}$. Then, we solve the least-squares problem $\mathbf{y}_m = {\operatorname{argmin}}_{\mathbf{y}} \| \mathbf{c}-\bar{\mathbf{H}}_m \mathbf{y} \|$ with $\mathbf{c} = \| \mathbf{r}_0 \| \mathbf{e}_1$ and compute the solution $\mathbf{x}_{m} = \mathbf{x}_{0} + \mathbf{Z}_{m}\mathbf{y}_{m}$.

To initiate the next cycle, we solve the standard eigenvalue problem \eqref{eq:Standard_eigenvalue_problem} and retain the eigenpairs associated to the $k$ eigenvalues with smallest  magnitudes. The corresponding eigenvectors form the columns of the matrix $\mathbf{P}_k$. Following definition \ref{HarmonicRitzpair}, we then define $\mathbf{Y}_k$ as the set of $k$ harmonic Ritz vectors of $\mathbf{AZ}_{m}\mathbf{V}_m^T$ with respect to $\operatorname{range}(\mathbf{V}_m)$. Defining $[\mathbf{Q}, \mathbf{R}]$ as the reduced QR-factorization of $\bar{\mathbf{H}}_{m} \mathbf{P}_k$, we have

\begin{equation}
	\mathbf{AZ}_m \mathbf{P}_k = \mathbf{V}_{m+1} \ols{\mathbf{H}}_m \mathbf{P}_k = (\mathbf{V}_{m+1}\mathbf{Q})\mathbf{R}.
	\label{eq:AZmPk_CkR}
\end{equation}

 We define the subspaces for the next cycle as $\mathbf{C}_{k} = \mathbf{V}_{m+1}\mathbf{Q}$ and $\mathbf{Z}_{k} = \mathbf{Z}_m \mathbf{P}_k \mathbf{R}^{-1}$, assuming that $\mathbf{R}$ is nonsingular. Since $\mathbf{V}_{m+1}$ and $\mathbf{Q}$ are orthonormal matrices, so is $\mathbf{C}_k$.
At the end of the FGMRES($m$) cycle, the optimality property $\mathbf{r}_1^{(1)} \; \bot \; \operatorname{span}\lbrace \mathbf{AZ}_m \rbrace$ of the residual $\mathbf{r}_1^{(1)}=\mathbf{r}_0^{(1)}-\mathbf{AZ}_m \mathbf{y}_m$ can be easily verified. The following theorem holds:
\vspace{6pt}
\begin{thm}
	The vector $\mathbf{z}_i \in \operatorname{span}\lbrace \mathbf{Z}_m \rbrace$ satisfies $\mathbf{z}_i = \underset{\mathbf{z}}{\operatorname{argmin}} \| \mathbf{r}_0 - \mathbf{Az} \|_2 \Leftrightarrow \mathbf{r}_i \; \bot \; \operatorname{span}\lbrace \mathbf{AZ}_m \rbrace$
\end{thm}

Combining this result with relation \eqref{eq:AZmPk_CkR}, we have
\begin{align}
	\mathbf{P}_k^T (\mathbf{AZ}_m)^{T} \mathbf{r}^{(1)}_1 &= \mathbf{0} \nonumber \\
	( \mathbf{V}_{m+1} \bar{\mathbf{H}}_m \mathbf{P}_k )^T \mathbf{r}^{(1)}_1 &= \mathbf{0} \nonumber \\
	\Rightarrow \mathbf{R}^T \mathbf{C}_k^T \mathbf{r}^{(1)}_1 &= \mathbf{0},
\end{align}

\noindent which shows that $ \mathbf{C}_k^T \mathbf{r}^{(1)}_1 = \mathbf{0} $, with $\mathbf{C}_{k} = \mathbf{V}_{m+1}\mathbf{Q}$.

To complete the next cycle $(i+1)$ and obtain the next iterates $\ve{r}_{i+1}^{(s)}$ and $\ve{x}_{i+1}^{(s)}$ for the current system $(s)$ we perform $m-k$ steps of FGMRES applied to the projected operator $\te{A}_{\te{C}_k}=(\matI - \te{C}_{k}\te{C}^{T}_{k})\te{A}$, i.e., we solve the projected system $\textbf{PA}\textbf{x} = \textbf{Pb}$ where $\textbf{P}=\textbf{I} -\textbf{C}_k\textbf{C}^{T}_k$ is the orthogonal projector on $\textrm{range}(\textbf{C}_k)^{\perp}$. This results in the generalized Arnoldi relation

\begin{equation}
	\textbf{AZ}_m = \textbf{V}_{m+1} \ols{\textbf{H}}_m,
	\label{eq:flexible_generalized_Arnoldi_relation}
\end{equation}
\noindent where
\begin{equation}
	\mathbf{Z}_m = \left[ \mathbf{Z}_k \quad  \mathbf{Z}_{m-k} \right],  \quad \mathbf{V}_{m+1} = \left[ \mathbf{C}_k \quad \mathbf{V}_{m-k+1} \right] \quad\textrm{and} \quad \ols{\mathbf{H}}_{m} = \begin{bmatrix}
		\mathbf{I}_k & \mathbf{B}_{m-k}  \\
		\mathbf{0} & \ols{\mathbf{H}}_{m-k} 
	\end{bmatrix}\!\!,
	\label{eq:flexible_generalized_Arnoldi_matrices}
\end{equation}
\noindent with $ \mathbf{B}_{m-k} = \mathbf{C}_k^{T} \mathbf{AZ}_{m-k}, \mathbf{B}_{m-k} \in \mathbb{R}^{k \times (m-k)}$.
To build $\mathbf{V}_{m+1}$ the Arnoldi recurrence restarts from $\mathbf{v}_{k+1} = \mathbf{r}_1 / \| \mathbf{r}_1 \|$ with $\mathbf{r}_1 = (\mathbf{I} - \mathbf{C}_{k}\mathbf{C}^{T}_{k})\mathbf{r}_{0}$.

To initiate the subsequent cycle or to solve the next linear system, we compute $k$ eigenvectors of a generalized eigenvalue problem that we will specify in section \ref{Flexible_deflation_FGCRODR} hereafter. These eigenvectors are concatenated in matrix $\mathbf{P}_k$. Again, introducing the reduced QR-factorization of $\bar{\mathbf{H}}_{m} \mathbf{P}_k$, we define the new residual and solution subspaces
\begin{align}
	\mathbf{C}_k &= \mathbf{V}_{m+1} \mathbf{Q} \\
	\mathbf{Z}_k &= \mathbf{Z}_m \mathbf{P}_k \mathbf{R}^{-1}.  \label{eq:Zk}
\end{align}

To initiate the solving of the next linear system we compute the projected initial solution and residual with
\begin{align}
    \mathbf{x}_{1} &= \mathbf{x}_{0} + \mathbf{Z}_{k}\mathbf{C}^{T}_{k}\mathbf{r}_{0} \\
	\mathbf{r}_{1} &= (\mathbf{I} - \mathbf{C}_{k}\mathbf{C}^{T}_{k})\mathbf{r}_{0}
\end{align}

\subsection{Flexible deflation strategies} \label{Flexible_deflation_FGCRODR}

From the definition of harmonic Ritz vectors, there exits different ways of deriving the deflated subspace according to the choice of operator $\textbf{B}$ and vector space $\mathcal{U}$. In the technical report \cite{carvalho2010report}, Carvalho et al. proposed three formulations (labeled A, B and C) for the projected generalized eigenvalue problem. We review these variants below and give additional insight for strategy B. \\

\noindent \underline{\textbf{Strategy A}}: 

\vspace{6pt}
This strategy has already been considered for FGMRES-DR in section \ref{StrategyA_FGMRESDR}, see relation \eqref{eq:Alternative_deflation_strategy_b}. Defining $\mathbf{B} \equiv \mathbf{AZ}_{m} \mathbf{Z}_m^{\dagger}$ and $\mathcal{U} \equiv \operatorname{range}(\mathbf{Z}_m)$, the harmonic Ritz pair $(\lambda,\mathbf{y}=\mathbf{Z}_m \mathbf{g})$ satisfies
\begin{equation}
	\label{eq:FGCRODR_StrategyA}
	((\mathbf{AZ}_{m} \mathbf{Z}_m^{\dagger}) \mathbf{Z}_m)^T ((\mathbf{AZ}_{m} \mathbf{Z}_m^{\dagger}) \mathbf{Z}_m \mathbf{g} - \lambda \mathbf{Z}_{m} \mathbf{g}) = \mathbf{0}.
\end{equation}

Thus, $\mathbf{Y}_m = \{ \mathbf{y}_1, \cdots, \mathbf{y}_m \}$ corresponds to the harmonic Ritz vectors of $\mathbf{AZ}_{m} \mathbf{Z}_m^{\dagger}$ with respect to $\operatorname{range}(\mathbf{Z}_m)$. Also, inserting \eqref{eq:flexible_generalized_Arnoldi_relation} in \eqref{eq:FGCRODR_StrategyA}, the eigenpair $(\lambda, \mathbf{g})$ satisfies the generalized eigenvalue problem
\begin{equation}
	\boxed{\Bar{\mathbf{H}}_{m}^{T} \Bar{\mathbf{H}}_{m}\mathbf{g} = \lambda  \bar{\mathbf{H}}^T_{m} \mathbf{V}^T_{m+1} \mathbf{Z}_m\mathbf{g}}.
\end{equation}

\noindent \underline{\textbf{Strategy B}}:

\vspace{6pt}
This strategy has also already been considered for FGMRES-DR in section \ref{DeflatedRestarting_FGMRESDR}, see relations \eqref{eq:Generalized_eigenvalue_problem_a} and \eqref{eq:Generalized_eigenvalue_problem_b}. Letting $\mathbf{V}_{m} = \left[ \textbf{C}_k \quad \textbf{V}_{m-k} \right]$ and choosing $\mathbf{B} \equiv \mathbf{AZ}_{m}\mathbf{V}_m^{T}$ and $\mathcal{U} \equiv \operatorname{range}(\mathbf{V}_m)$, the harmonic Ritz pair $(\lambda,\mathbf{y}=\mathbf{V}_m \mathbf{g})$ satisfies
\begin{equation}
	((\mathbf{AZ}_{m}\mathbf{V}_m^{T}) \mathbf{V}_m)^T ((\mathbf{AZ}_{m}\mathbf{V}_m^{T}) \mathbf{V}_m \mathbf{g} - \lambda \mathbf{V}_m \mathbf{g}) = \mathbf{0}.
\end{equation}

Thus, $\mathbf{Y}_m = \{ \mathbf{y}_1, \cdots, \mathbf{y}_m \}$ corresponds to the harmonic Ritz vectors of $\mathbf{AZ}_{m}\mathbf{V}_m^{T}$ with respect to $\operatorname{range}(\mathbf{V}_m)$. Also, the eigenpair $(\lambda, \mathbf{g})$ satisfies the generalized eigenvalue problem
\begin{equation}
\boxed{\Bar{\mathbf{H}}_{m}^{T} \Bar{\mathbf{H}}_{m}\mathbf{g} = \lambda \mathbf{H}_{m}^{T} \mathbf{g}}.
\label{eq:FGCRODR_generalized_eigenvalue_problem}
\end{equation}

This generalized eigenvalue problem can be reformulated as a standard one, see equation \eqref{eq:Standard_eigenvalue_problem}. We have $\hat{\mathbf{H}}_m\mathbf{g}=\lambda\mathbf{g}$, with $\hat{\mathbf{H}}_m = [ \mathbf{H}_m + h^2 \mathbf{f} \, \mathbf{e}_m^T ]$ and $\mathbf{H}_m = [\mathbf{I}_m \quad \mathbf{0}_{m \times 1}]\ols{\mathbf{H}}_m$. The block upper triangular structure of $\ols{\mathbf{H}}_m$ is highlighted in \eqref{eq:flexible_generalized_Arnoldi_matrices} with a leading $k \times k$ identity block. More specifically we have
\begin{equation}
\hat{\mathbf{H}}_m =
\begin{bmatrix}
	\mathbf{I}_k & \tilde{\mathbf{B}}_{m-k}  \\
	\mathbf{0} & \tilde{\mathbf{H}}_{m-k} 
\end{bmatrix}
\label{eq:FGCRODR_generalized_matrix}
\end{equation}
with $\tilde{\mathbf{B}}_{m-k} = \mathbf{B}_{m-k} + h^2 \mathbf{f}_{[1:k]} \mathbf{e}_{m-k}^T$ and $\tilde{\mathbf{H}}_{m-k} = \mathbf{H}_{m-k} + h^2 \mathbf{f}_{[k+1:m]} \mathbf{e}_{m-k}^T$, $\mathbf{e}_{m-k} = [0 \cdots 0, 1]^T \in \mathbb{R}^{m-k}$.

We assume here that the Hessenberg sub-matrix $\ols{\mathbf{H}}_{m-k}$ is unreduced, meaning that there are no zero elements on the subdiagonal (the problem of multiple eigenvalues is discussed in \cite{morgan_Implicit_2000}). The eigenvalues of \eqref{eq:FGCRODR_generalized_matrix} are the combined eigenvalues of the diagonal blocks of $\hat{\mathbf{H}}_m$ and therefore satisfy $\operatorname{det}(\hat{\mathbf{H}}_m - \lambda \mathbf{I}_m) = \operatorname{det}((1-\lambda)\mathbf{I}_k) \operatorname{det}(\tilde{\mathbf{H}}_{m-k} -\lambda\mathbf{I}_{m-k}) = 0$. Thus, we have a unit eigenvalue of algebraic multiplicity $k$.
Now, if $\lambda=1$ is an eigenvalue of the upper diagonal block $\mathbf{I}_k$, with associated eigenvectors $\mathbf{g}_k=\mathbf{I}_k$, then it is also an eigenvalue of the full matrix $\hat{\mathbf{H}}_m$, with the same eigenvectors augmented with zeros, which gives $\mathbf{g}_m^{(\lambda=1)} = [ \mathbf{g}_k \quad \mathbf{0}_{(m-k) \times k}]^T = [ \mathbf{I}_k \quad \mathbf{0}_{(m-k) \times k}]^T$. This shows that the $k$ eigenvectors associated to the unit eigenvalues are linearly independent.
We can also obtain an explicit expression for the complementary eigenvectors associated to the eigenvalues of the lower diagonal block of $\hat{\mathbf{H}}_m$. Let $(\lambda_{m-k},\mathbf{g}_{m-k})$ be an eigenpair of $\tilde{\mathbf{H}}_{m-k}$. The full eigenvalue problem reads
\begin{equation}
\begin{bmatrix}
	\mathbf{I}_k & \tilde{\mathbf{B}}_{m-k}  \\
	\mathbf{0} & \tilde{\mathbf{H}}_{m-k} 
\end{bmatrix}
\begin{pmatrix}
	\mathbf{x}_k \\
	\mathbf{g}_{m-k}
\end{pmatrix}
= \begin{pmatrix}
	\mathbf{x}_k +\tilde{\mathbf{B}}_{m-k} \mathbf{g}_{m-k} \\
	\lambda_{m-k} \mathbf{g}_{m-k}
\end{pmatrix}\!\!.
\end{equation}

We can make $\mathbf{x}_k +\tilde{\mathbf{B}}_{m-k} \mathbf{g}_{m-k} = \lambda_{m-k} \mathbf{x}_k$ by choosing $\mathbf{x}_k=-(1-\lambda_{m-k})^{-1}  \tilde{\mathbf{B}}_{m-k} \mathbf{g}_{m-k}$. By assumption $\lambda_{m-k} \neq 1$ since it is not an eigenvalue of the leading block of $\hat{\mathbf{H}}_m$. Thanks to these relations, the deflated subspace associated to strategy B is obtained efficiently. In light of this attractive property, Jolivet and Tournier proposed a block version of GCRO-DR combined with this deflation strategy \cite{jolivet2016block}.\\

\noindent \underline{\textbf{Strategy C}}:

\vspace{6pt}
This third strategy was adopted in \cite{carvalho2011flexible} where the FGCRO-DR algorithm maintains an additional set of vectors $\mathbf{W}_m$ such that the deflated harmonic Ritz vectors read $\mathbf{Y}_k = \mathbf{W}_m \mathbf{P}_k$. Similar to $\mathbf{Z}_k$ in \eqref{eq:Zk}, a set $\mathbf{W}_k = \mathbf{W}_m \mathbf{P}_k \mathbf{R}^{-1}$ is propagated between cycles. The complete subspace is built by appending the Krylov vectors generated during the Arnoldi process: $\mathbf{W}_{m} = \left[ \textbf{W}_k \quad \textbf{V}_{m-k} \right]$. This particular subspace was used to prove the algebraic equivalence of FGMRES-DR and FGCRO-DR under a specific colinearity constraint. As a consequence, strict algebraic equivalence of GMRES-DR and GCRO-DR was also demonstrated by the authors. Note that this property was stated in \cite{parks2006recycling} without demonstration.

Using definition \ref{HarmonicRitzpair} with $\mathcal{U} \equiv \operatorname{range}(\mathbf{W}_m) $ and $\mathbf{B} \equiv \mathbf{AZ}_{m}\mathbf{W}_m^{\dagger}$ we can write $\mathbf{y} = \mathbf{W}_m \mathbf{g}$ which gives the orthogonality constraint
\begin{equation}
	((\mathbf{AZ}_{m}\mathbf{W}_m^{\dagger}) \mathbf{W}_m)^T ((\mathbf{AZ}_{m}\mathbf{W}_m^{\dagger}) \mathbf{W}_m \mathbf{g} - \lambda \mathbf{W}_m \mathbf{g}) = \mathbf{0}.
\end{equation}

Now, using \eqref{eq:flexible_generalized_Arnoldi_relation} and recalling that $\mathbf{V}_{m+1}$ is orthonormal, we obtain the generalized eigenvalue problem 
\begin{equation}
	\boxed{\bar{\mathbf{H}}^T_m \bar{\mathbf{H}}_m \mathbf{g} = \lambda \bar{\mathbf{H}}^T_m \mathbf{V}^T_{m+1} \mathbf{W}_m \mathbf{g}}
\end{equation}

To reduce ill-conditioning, an alternative formulation for this generalized eigenvalue problem has been proposed in equation \eqref{eq:Reformulated_GCRO_eigenproblem}. \\


\subsection{Application of FGCRO-DR to the fluid adjoint system} 

\begin{figure}[H]
	\captionsetup[subfigure]{belowskip=-6pt}
	\begin{subfigure}{.5\textwidth}
		\raggedright
		\includegraphics[trim={1.0cm 0.5cm 1.cm 1.9cm},clip,width=0.99\textwidth]{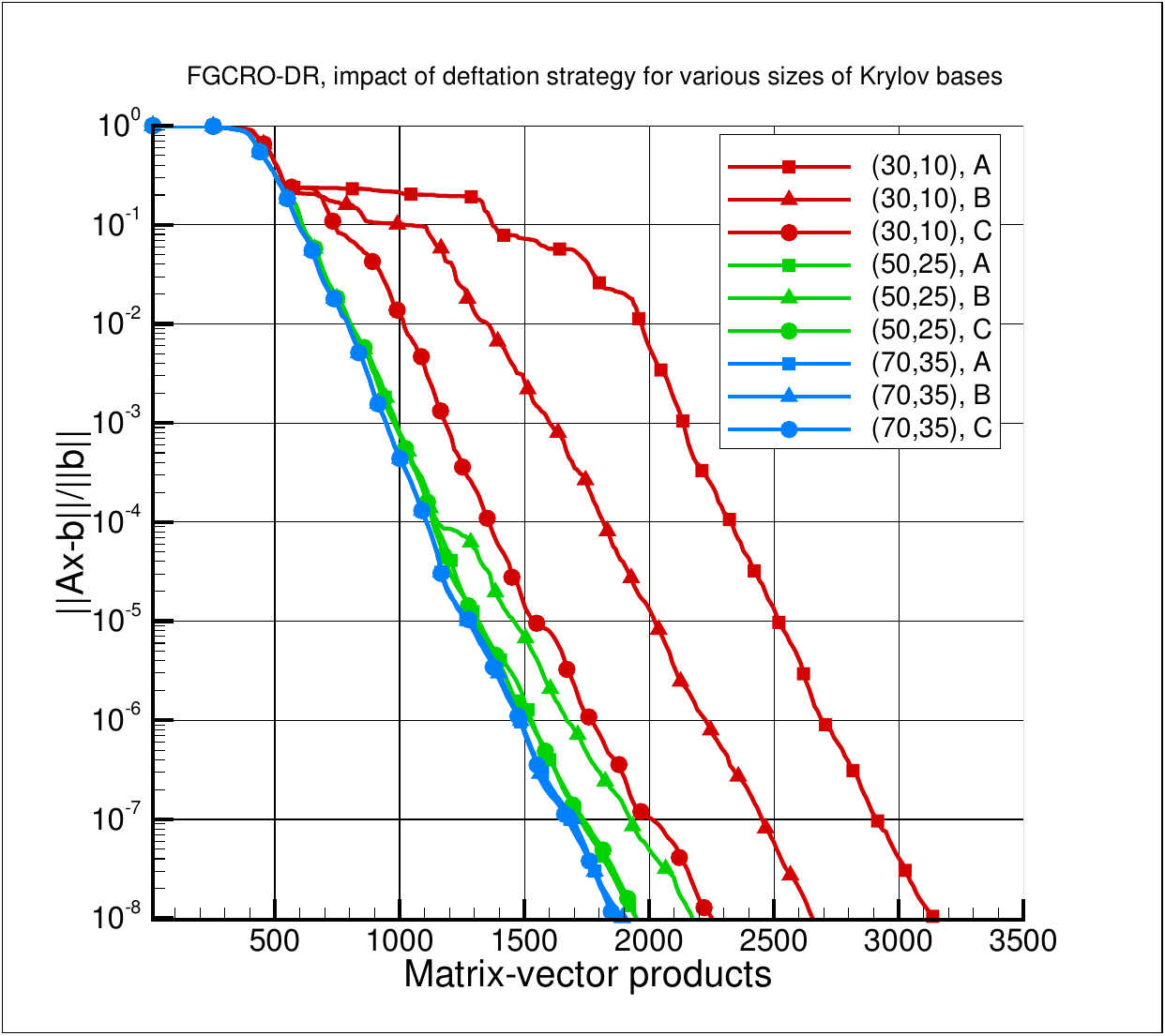}
		\caption{}
		\label{fig:cvg_OPT_FGCRODR_impact_deflation_strategy_various_sizes_Krylov_space_matvecprod}
	\end{subfigure}
	\begin{subfigure}{.5\textwidth}
		\centering
		\includegraphics[trim={1.0cm 0.5cm 1.cm 1.9cm},clip,width=0.99\textwidth]{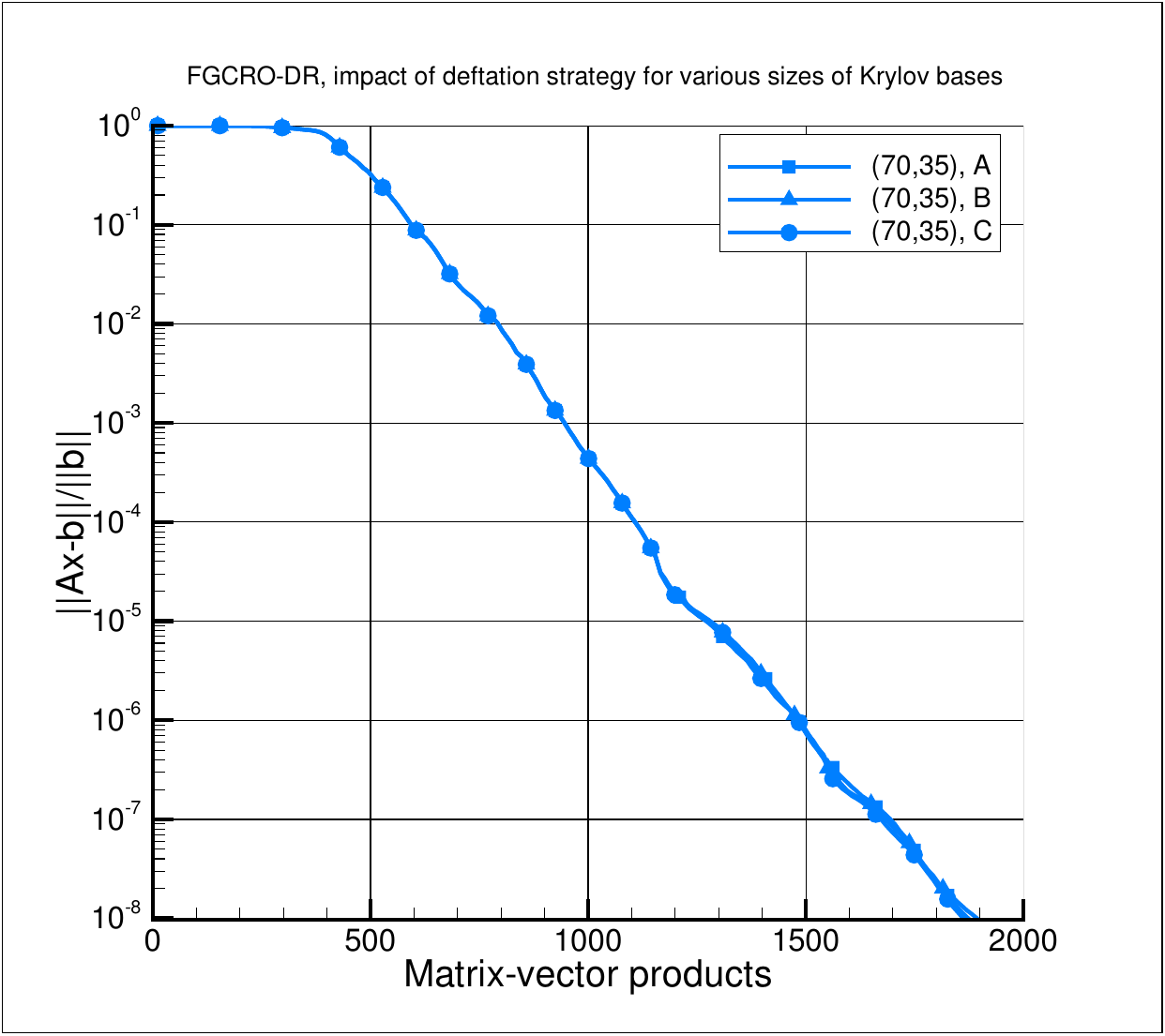}
		\caption{}		
		\label{fig:cvg_OPT_FGCRODR_70_10_35_impact_deflation_strategy_matvecprod}
	\end{subfigure}
	\vspace{0pt}
	\caption{Impact of the deflation strategy on the convergence of FGCRO-DR for a varying size of the external Krylov subspace. The size of the inner Krylov subspace is 10. Deflation strategy C outperforms strategies A and B. Strategy A seems the most sensitive to the size of the Krylov basis. For the larger Krylov subspace FGCRO-DR(70,35) performs equally well regardless of the chosen strategy.}
	\label{fig:cvg_OPT_FGCRODR_impact_deflation_strategy_various_sizes_Krylov_space}
\end{figure}

Our first numerical experiments compare the three deflation strategies detailed in section \ref{Flexible_deflation_FGCRODR} for the solving of the fluid adjoint linear system. We consider increasing sizes of Krylov bases of 30, 50 and 70 vectors. The deflated subspace size is half of the Krylov basis size except for the smaller basis where we keep only one third of the approximation subspace. The variable preconditioner consists in a non-restarted inner GMRES Krylov solver with a basis of size 10. The innermost stationary preconditioner is BILU(0) applied to the first-order exact Jacobian matrix. The corresponding convergence histories are presented in the left-hand side of Figure \ref{fig:cvg_OPT_FGCRODR_impact_deflation_strategy_various_sizes_Krylov_space}. Clearly, strategy C seems the most effective regardless the size of the Krylov space. However, when the number of vector increases, the three deflation strategies perform similarly. Indeed, the right-hand side of Figure \ref{fig:cvg_OPT_FGCRODR_impact_deflation_strategy_various_sizes_Krylov_space} shows very close convergence profiles of FGCRO-DR(70,35) for the three deflation strategies. There is no theoretical evidence that a deflation strategy should perform better than another one. As a consequence, these conclusions are likely to be case dependent.

\subsection{Application of FGCRO-DR to the solution of the fluid-structure coupled-adjoint system}

In this section we study the impact of recycling spectral information between fluid-structure cycles in the context of GCRO-DR with variable preconditioning. From our previous numerical experiments with FGMRES-DR and FGCRO-DR applied to the fluid adjoint system, we choose to recycle a subspace of size half the size of the Krylov subspace, i.e., the same size as the one of the deflated subspace. In \cite{carvalho2010report} the authors suggest that larger recycled subspaces give better performance in the context of a flexible Krylov solver used to solve linear systems arising from a specific class of problems in quantum chromodynamics. For instance, for an outer space of size $m=20$, they carried out a parameter study by varying the number of deflated vectors from 1 to $m-1$. Regardless of the selected recycling strategy, the lowest computational cost was achieved for a size of the recycled subspace higher than 50\% of the size of the Krylov subspace. In another technical report \cite{carvalho2011report} the same authors solve a sequence of twelve elliptic partial differential equation problems for an increasing grid size. Their numerical experiments show that a reduction in the range of 40\% to 45\% can be achieved by recycling an approximate invariant subspace of half the size of the Krylov vector space compared to a standard FGMRES-DR applied to each system separately. These results seem to confirm the similar gains observed for the non-flexible GCRO-DR Krylov solver.
\bigbreak
We recall that in our context, the sequence of linear adjoint systems corresponds to a varying right-hand side arising from the update of the structural source term during the partitioned fluid-structure solution strategy. After a number of fluid-structure cycles, the right-hand side is expected to converge and we end up with a constant linear system to solve.
We then apply FGCRO-DR($m=70$,$m_i=10$,$k=35$), $k$ being the number of deflated and recycled vectors, combined with the three deflation strategies A, B and C. Table \ref{tab:strategies_deflation_performances} collects the corresponding numbers of matrix-vector products. Using recycling helps to improve the convergence rate of flexible GCRO-DR in this application since a reduction of approximately 16\% to 19\% in terms of matrix-vector products is obtained.

\begin{table}[H]
	\centering
	\renewcommand{\arraystretch}{1.1}			
	\begin{tabular}{ %
		|>{\centering}m{3cm}|>{\centering}m{2.5cm}|>{\centering}m{2.5cm}|>{\centering\arraybackslash}m{2.5cm}|}		
		\hline 
		Starting recycling cycle & \makecell{Strategy A \\ \#Mvps} & \makecell{Strategy B \\ \#Mvps} & \makecell{Strategy C \\ \#Mvps}\\ 
		\hline 
		N/A & 10815 & 10815 & 10815 \\
		\hline 
		5 & \textbf{8940} & 10363 & 9768 \\
		\hline 
		6 & 9829 & 9107 & \textbf{8754} \\ 
		\hline 
		7 & 9348 & \textbf{9108} & 9341 \\ 
		\hline 
		8 & \textbf{8790} & 9862 & \textit{12018} \\ 
		\hline 
		9 & \textbf{9250} & 9279 & 8673 \\ 
		\hline
	\end{tabular} 
	\caption{Total number of matrix-vector products \#Mvps for strategies A, B and C, indexed by the starting recycling cycle. The lowest \#Mvps in each row is highlighted in bold and should be compared to the reference \#Mvps=10815. }
	\label{tab:strategies_deflation_performances}
\end{table}

\begin{figure}[H]
	\captionsetup[subfigure]{belowskip=-6pt}
	\begin{subfigure}{.5\textwidth}
		\raggedleft
		\includegraphics[trim={1.1cm 0.5cm 1.cm 1.9cm},clip,width=0.99\textwidth]{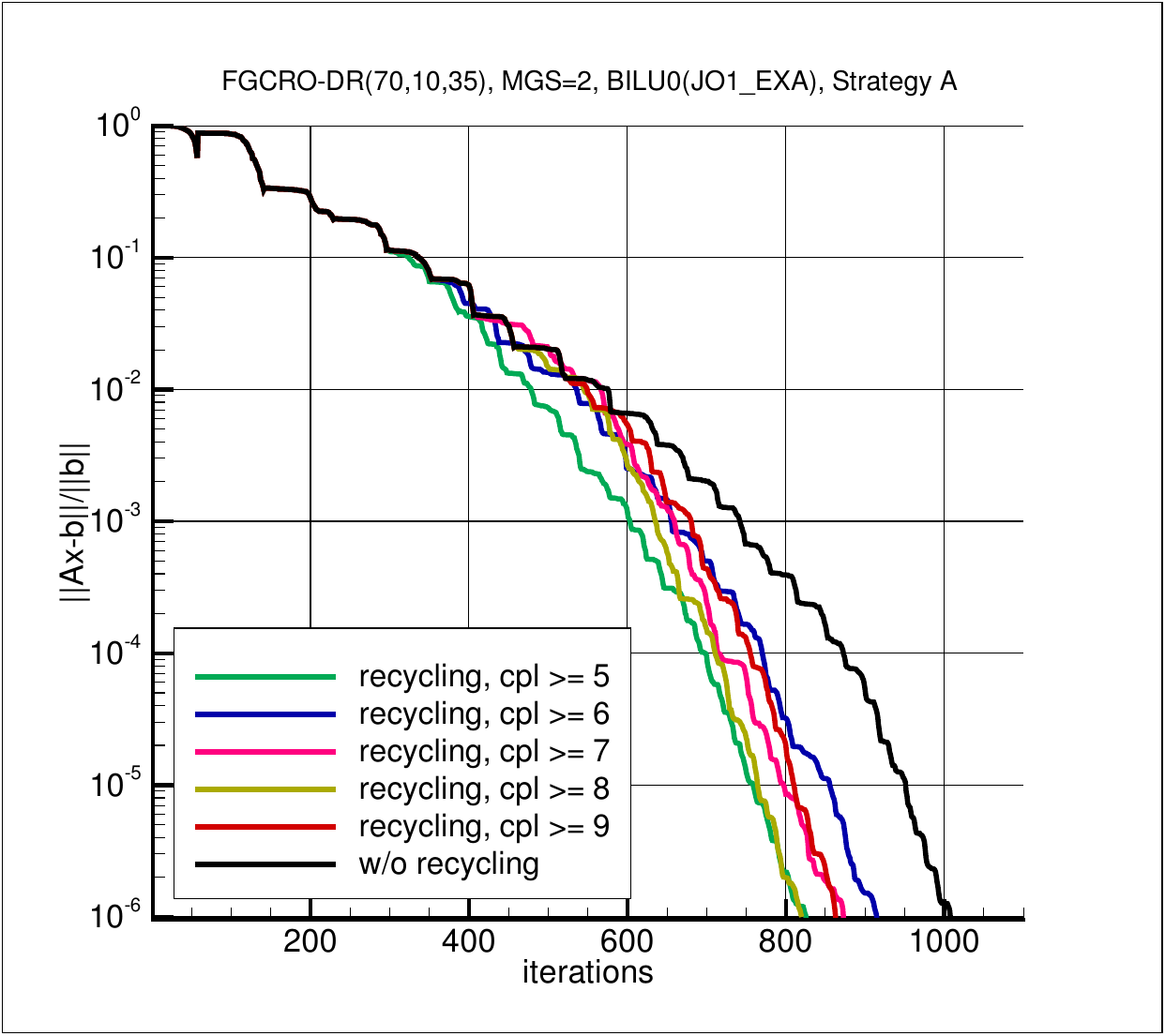}
		\caption{}
		\label{fig:cvg_AOC_FGCRODR_70_35_RATIO_060_MGS2_impact_deflation_strategy_A_iterations}
	\end{subfigure}
	\begin{subfigure}{.5\textwidth}
		\raggedleft
		\includegraphics[trim={1.1cm 0.5cm 1.cm 1.9cm},clip,width=0.99\textwidth]{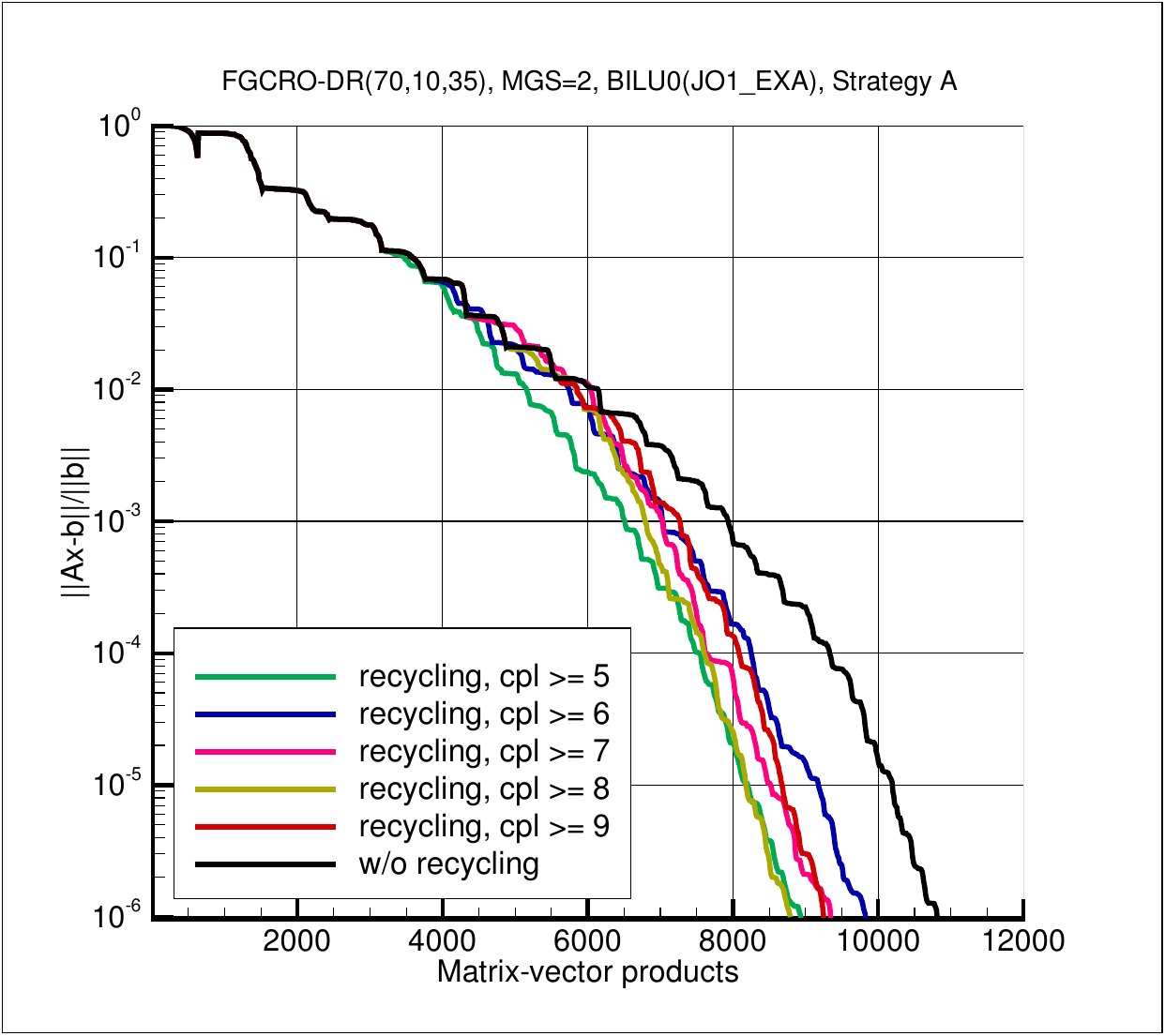}
		\caption{}		
		\label{fig:cvg_AOC_FGCRODR_70_35_RATIO_060_MGS2_impact_deflation_strategy_A_matvecprod}
	\end{subfigure}
	\vspace{0pt}
	\caption{Impact of approximate invariant subspace recycling strategy A on the relative residual convergence of the coupled-adjoint fluid block for FGCRO-DR(70,10,35), with innermost preconditioner BILU0($\textbf{J}^{EXA}_{O1}$). Recycling spectral information starting from cycle 5 and above always improves convergence. A maximum saving of 2065 matrix-vector products out of 10815 ($\sim$19\%) is achieved.}
	\label{fig:res_hist_M6W_AOC_fgcrodr_strategy_A_MGS2}
\end{figure}

\begin{figure}[H]
	\captionsetup[subfigure]{belowskip=-6pt}
	\begin{subfigure}{.5\textwidth}
		\raggedleft
		\includegraphics[trim={0.25cm 0.5cm 0.8cm 1.8cm},clip,width=0.99\textwidth]{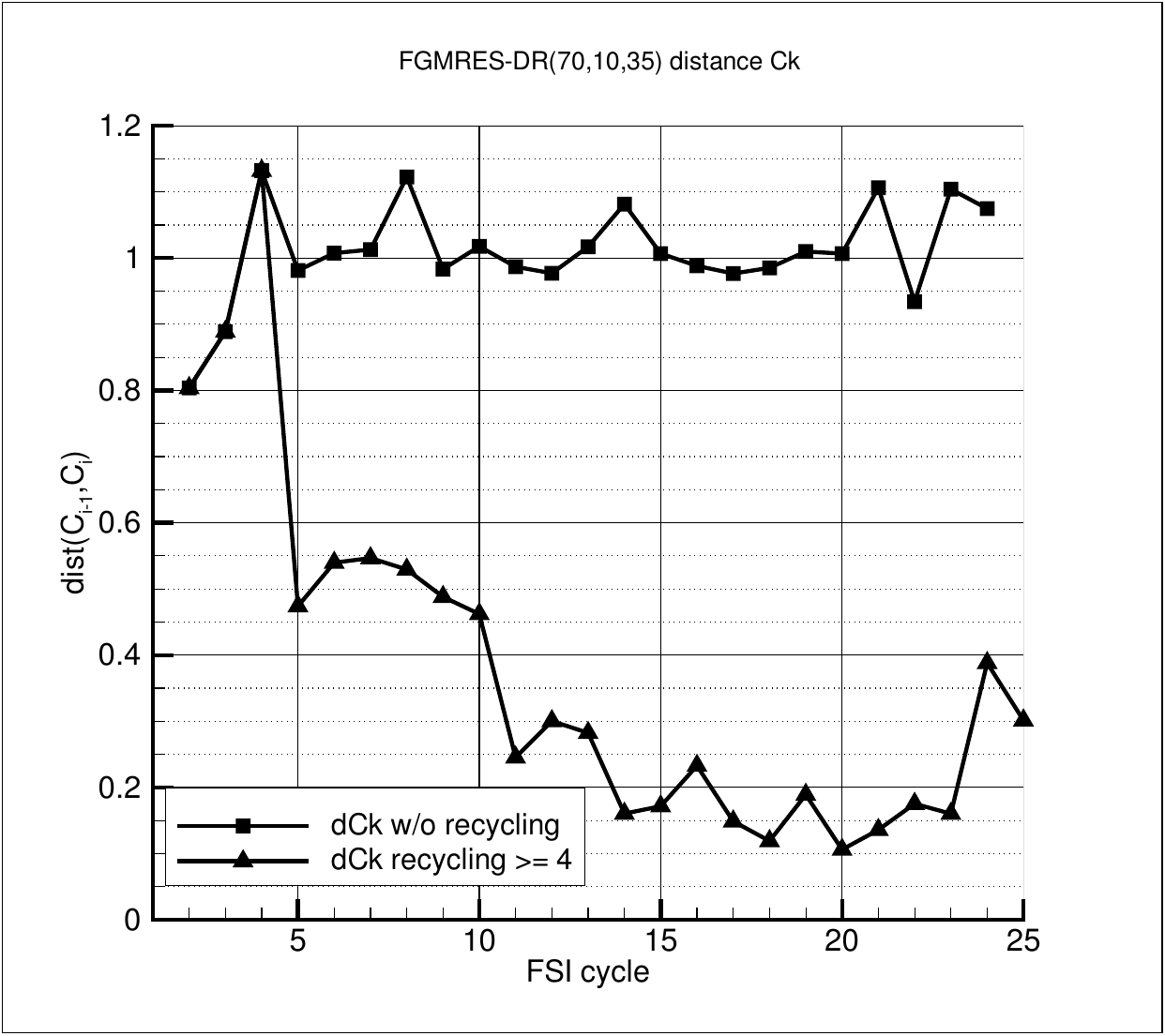}
		\caption{}
		\label{fig:FGMRES-DR_70_10_35_Ck_distances}
	\end{subfigure}
	\begin{subfigure}{.5\textwidth}
		\raggedleft
		\includegraphics[trim={0.25cm 0.5cm 0.8cm 1.8cm},clip,width=0.99\textwidth]{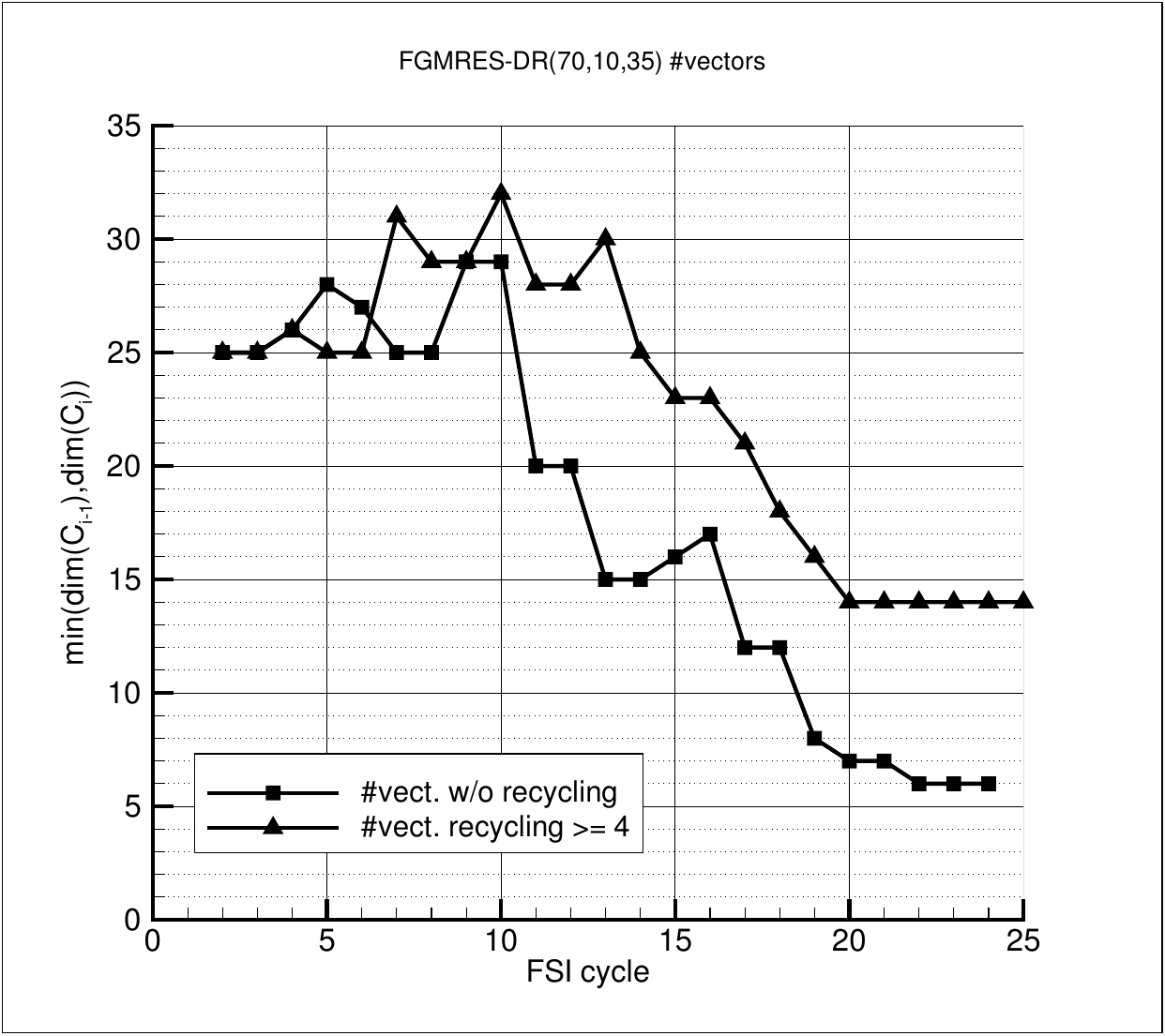}
		\caption{}		
		\label{fig:FGMRES-DR_70_10_35_Ck_vectors}
	\end{subfigure}
	\vspace{0pt}
	\caption{Impact of subspace recycling for FGMRES-DR(70,10,35). Monitoring of distance between approximate invariant subspaces $\mathbf{C}^{(i-1)}$ and $\mathbf{C}^{(i)}$ is reported in \ref{fig:FGMRES-DR_70_10_35_Ck_distances}. The number of vectors used to compute the subspace distance is also reported in \ref{fig:FGMRES-DR_70_10_35_Ck_vectors}.}
	\label{fig:FGMRES-DR_70_10_35_subspace_recycling_monitoring}
\end{figure}

\begin{figure}[H]
	\captionsetup[subfigure]{belowskip=-6pt}
	\begin{subfigure}{.5\textwidth}
		\raggedleft
		\includegraphics[trim={1.1cm 0.5cm 1.cm 1.9cm},clip,width=0.99\textwidth]{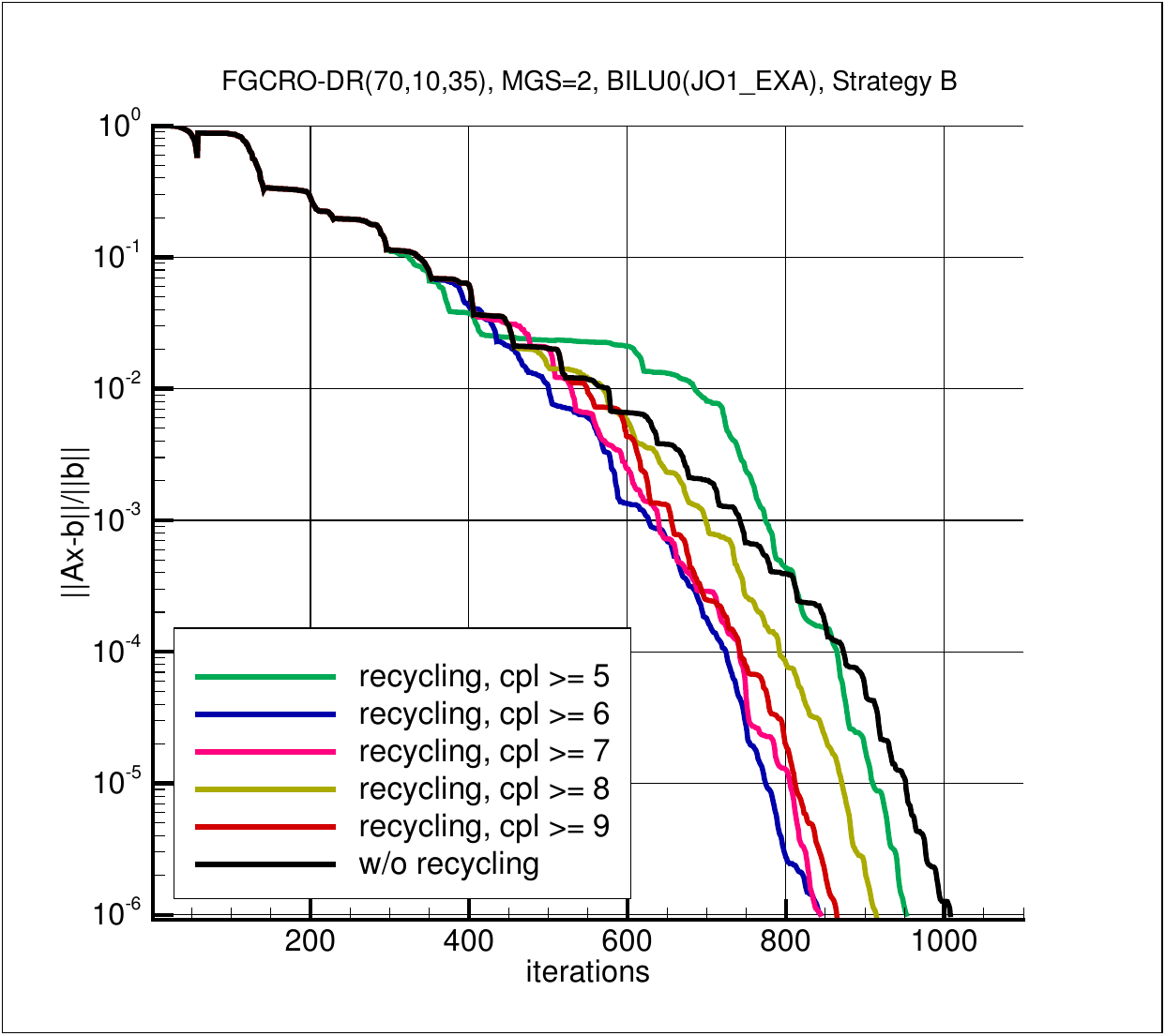}
		\caption{}
		\label{fig:cvg_AOC_FGCRODR_70_35_RATIO_060_MGS2_impact_deflation_strategy_B_iterations}
	\end{subfigure}
	\begin{subfigure}{.5\textwidth}
		\raggedleft
		\includegraphics[trim={1.1cm 0.5cm 1.cm 1.9cm},clip,width=0.99\textwidth]{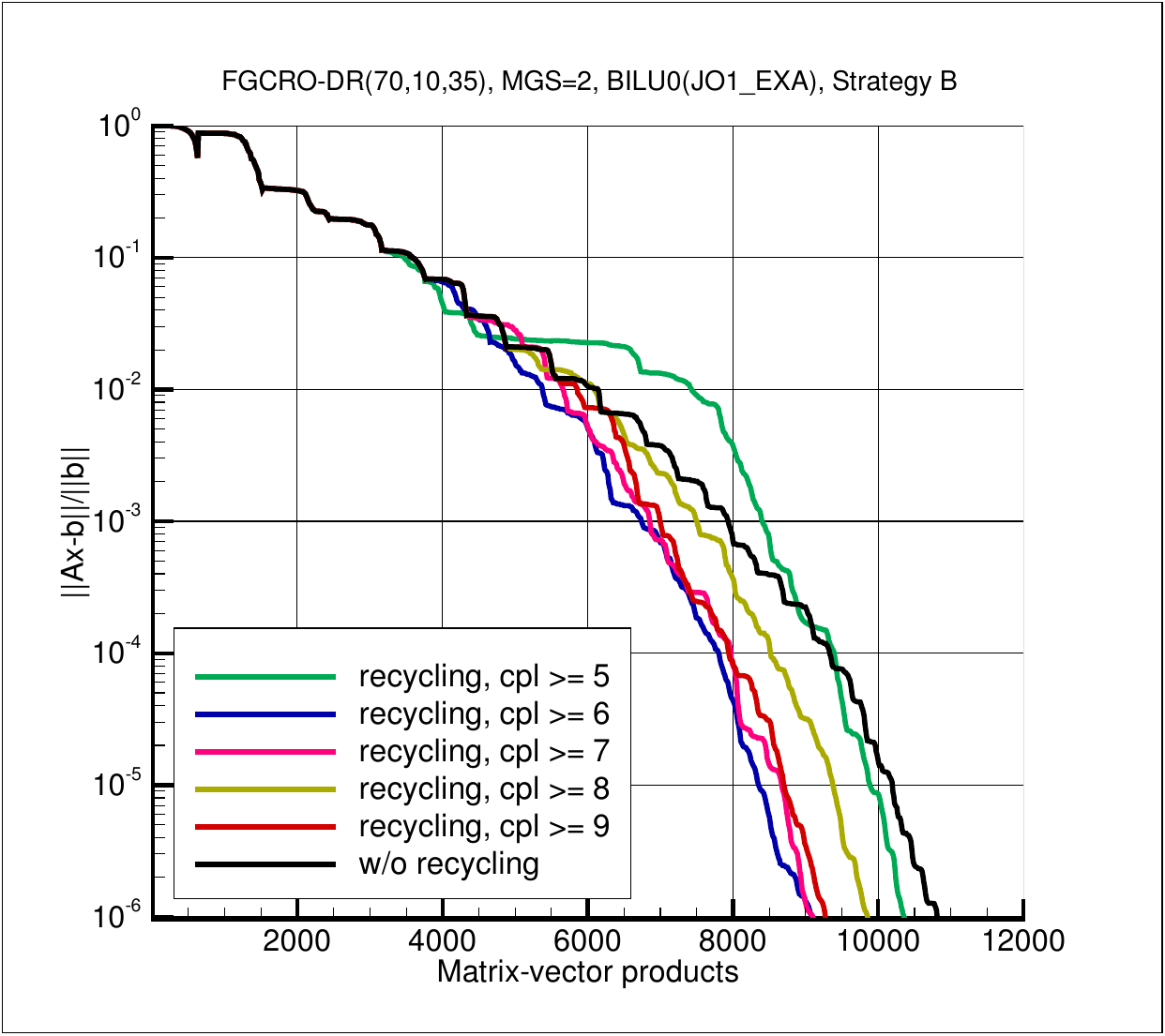}
		\caption{}		
		\label{fig:cvg_AOC_FGCRODR_70_35_RATIO_060_MGS2_impact_deflation_strategy_B_matvecprod}
	\end{subfigure}
	\vspace{0pt}
	\caption{Impact of approximate invariant subspace recycling strategy B on the relative residual convergence of the coupled-adjoint fluid block for FGCRO-DR(70,10,35), with innermost preconditioner BILU0($\textbf{J}^{EXA}_{O1}$). Recycling spectral information starting from cycle 5 and above always improves convergence. A maximum saving of 1765 matrix-vector products out of 10815 ($\sim$16\%) is achieved.}
	\label{fig:res_hist_M6W_AOC_fgcrodr_strategy_B_MGS2}
\end{figure}

Figures \ref{fig:res_hist_M6W_AOC_fgcrodr_strategy_A_MGS2}, \ref{fig:FGMRES-DR_70_10_35_subspace_recycling_monitoring} and \ref{fig:res_hist_M6W_AOC_fgcrodr_strategy_B_MGS2} show the convergence histories of FGCRO-DR(70,10,35) for strategy A, B and C respectively. We have reported convergences for a varying starting recycling fluid-structure cycle. Regardless the deflation strategy, the plateau following a fluid-structure coupling restart is always noticeably reduced.

A comparison of distances between FGCRO-DR without recycling and FGCRO-DR with a subspace recycling starting from cycle 5 is reported in Figure \ref{fig:FGMRES-DR_70_10_35_Ck_distances}. Clearly the subspace distance is much smaller when recycling is activated, whereas it keeps constant for the solver with cold restart after a fluid-structure coupling. The dimension $p$ in \eqref{eq:Grassmann_distance} is plotted if Figure \ref{fig:FGMRES-DR_70_10_35_Ck_vectors}. With recycling, the dimension of the approximation Krylov space is larger which favors a better convergence.

\begin{figure}[H]
	\captionsetup[subfigure]{belowskip=-6pt}
	\begin{subfigure}{.5\textwidth}
		\raggedleft
		\includegraphics[trim={1.1cm 0.5cm 1.cm 1.9cm},clip,width=0.99\textwidth]{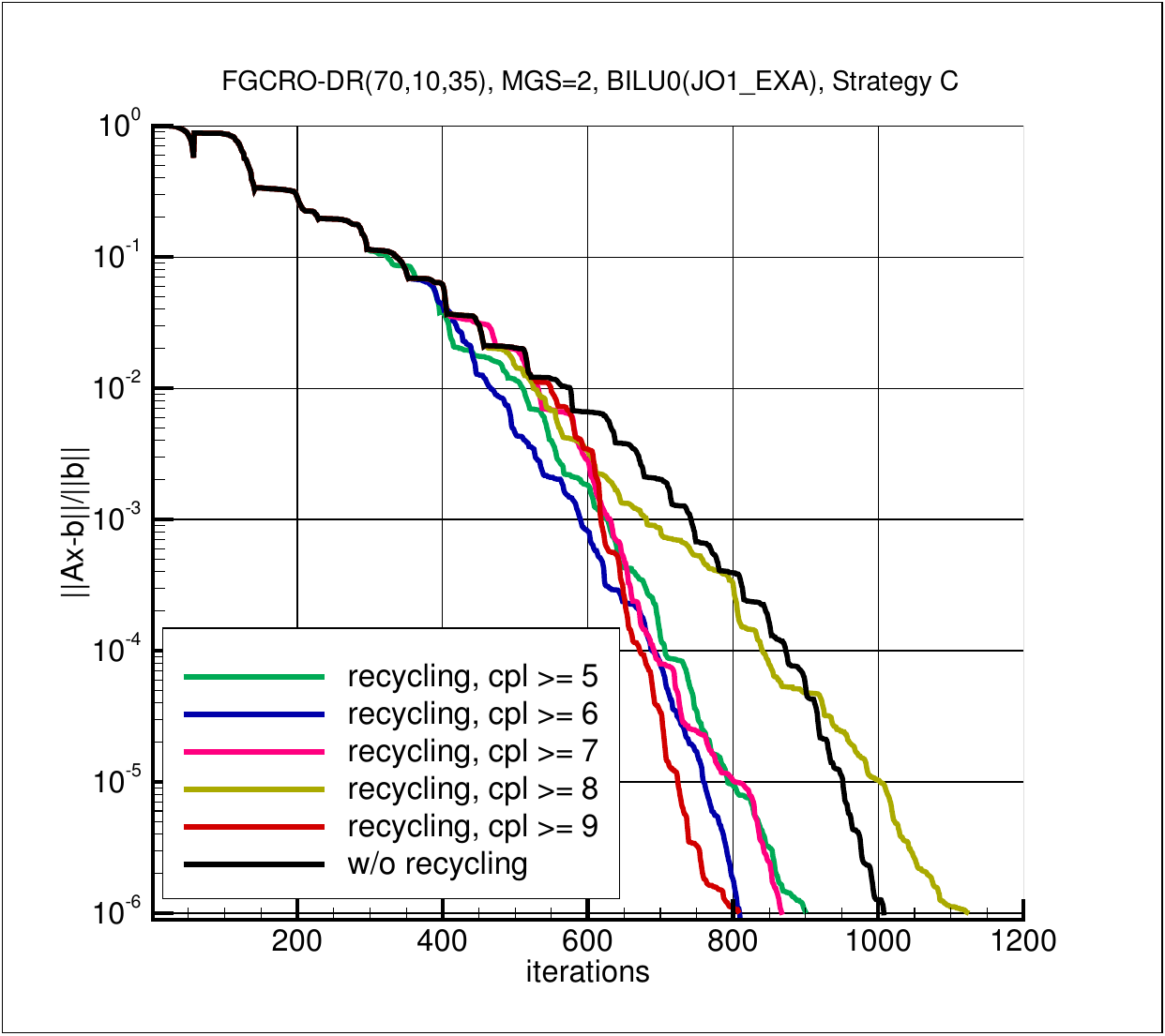}
		\caption{}
		\label{fig:cvg_AOC_FGCRODR_70_35_RATIO_060_MGS2_impact_deflation_strategy_C_iterations}
	\end{subfigure}
	\begin{subfigure}{.5\textwidth}
		\raggedleft
		\includegraphics[trim={1.1cm 0.5cm 1.cm 1.9cm},clip,width=0.99\textwidth]{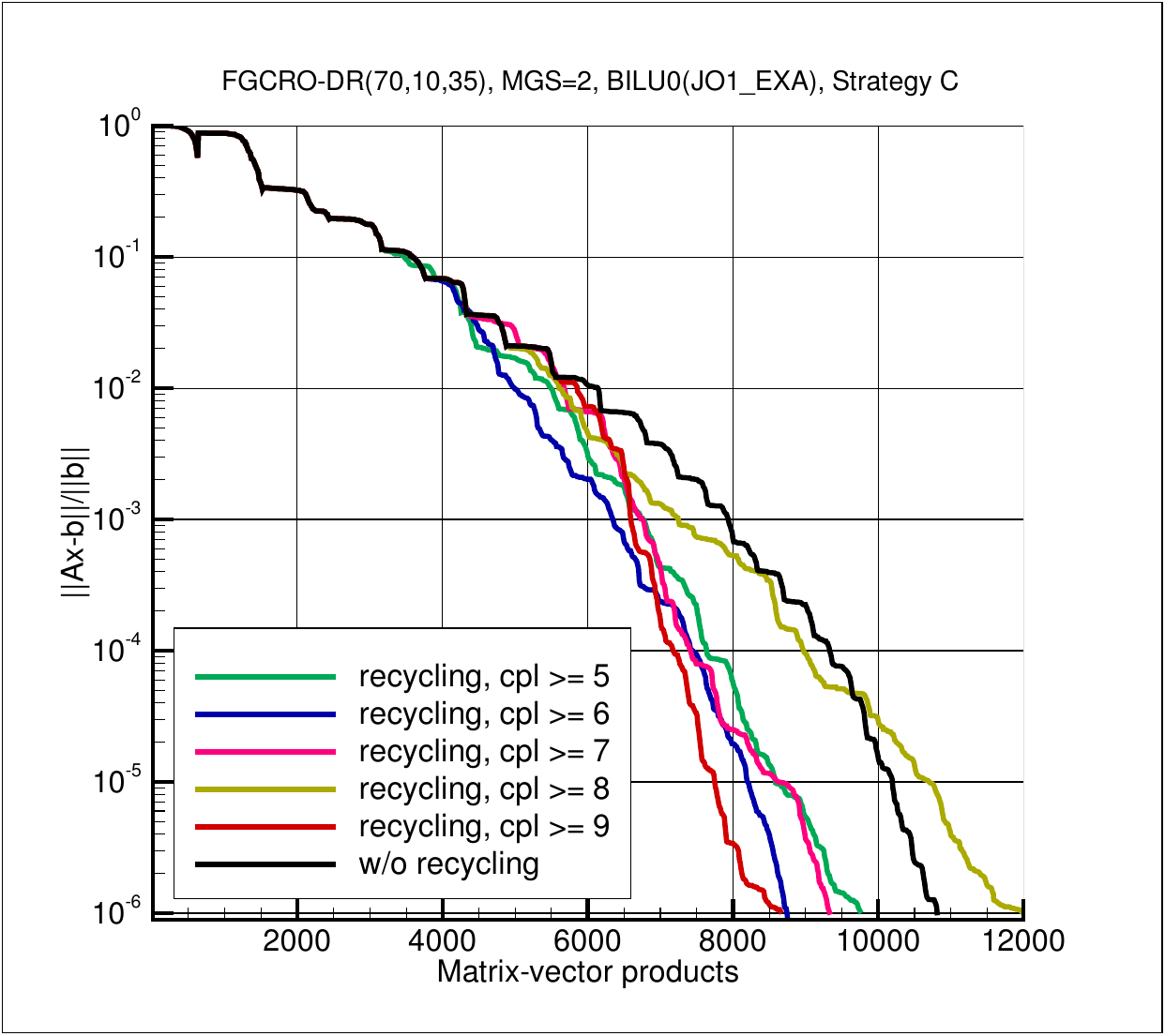}
		\caption{}		
		\label{fig:cvg_AOC_FGCRODR_70_35_RATIO_060_MGS2_impact_deflation_strategy_C_matvecprod}
	\end{subfigure}
	\vspace{0pt}
	\caption{Impact of approximate invariant subspace recycling strategy C on the relative residual convergence of the coupled-adjoint fluid block for FGCRO-DR(70,10,35), with innermost preconditioner BILU0($\textbf{J}^{EXA}_{O1}$). Recycling spectral information starting from cycle 5 and above always improves convergence, except when we start recycling from fluid-structure cycle 8. A maximum saving of 2050 matrix-vector products out of 10815 ($\sim$19\%) is achieved.}
	\label{fig:res_hist_M6W_AOC_fgcrodr_strategy_C_MGS2}
\end{figure}

\newpage

\section{Conclusion}

In this paper we have accelerated the fluid-structure coupled-adjoint partitioned solver by considering techniques borrowed from approximate invariant subspace recycling strategies adapted to sequences of linear systems with varying right-hand sides. Indeed, in a partitioned framework, the structural source term attached to the fluid block of equations affects the right-hand side with the nice property of quickly converging. This recycling and deflation strategies considered in this work were inspired by the theoretical developments related to GCRO-DR and its flexible variant detailed in \cite{carvalho2010report,carvalho2011flexible}. Our objective was to make this paper as self-contained as possible and in that respect, we choose to recall theoretical details of both GCRO-DR and FGCRO-DR. It also gave us the opportunity to point out practical implementation details mainly to make the computation of harmonic Ritz vector more stable.

We demonstrate the benefit of these techniques by computing the coupled derivatives for an aeroelastic configuration of the ONERA-M6 fixed wing in transonic flow representative of the numerical complexity of realistic industrial applications. For this exercise the fluid grid was coupled to a structural model specifically designed to exhibit a high flexibility. All computations were performed using RANS flow modeling and a fully linearized one-equation Spalart-Allmaras turbulence model, which typically results in stiff linear systems. 

For the non flexible GCRO-DR solver, subspace recycling between fluid-structure cycles achieved a reduction of 39\% in terms of matrix-vector products with respect to the legacy solver. Even recycling early from the second cycle led to a reduction in computational cost of 37\% showing the robustness of the proposed strategy. The conclusions with respect to the Flexible GCRO-DR are more mitigated, even if recycling almost systematically ended up with a reduction of the total number of matrix-vector products. Indeed, gains were achieved when recycling was activated from the fourth fluid-structure cycle only, regardless of the deflation strategy A, B or C. Deflation strategy A seems to be more consistent when comparing convergence profiles for various starting recycling cycles. For our application, a reduction of approximately 16\% to 19\% in terms of matrix-vector products was still obtained. However, the convergence seems less smooth exhibiting some local stagnations. This suggests that the recycled subspace may not always be suitable and a dynamic recycling approach would certainly be beneficial in this case. This will be part of a future work.

\newpage

\bibliography{Paper_for_arXiv_Recycling_Subspaces_for_coupledAdjoint_v2}

\end{document}

%% file: macros.tex
\newcommand{\ve}[1]{\mathbf{#1}}					
\newcommand{\te}[1]{\boldsymbol{\operatorname{#1}}}			

\newcommand{\matI}{\te{I}}					

\DeclarePairedDelimiter\abs{\lvert}{\rvert}%
\DeclarePairedDelimiter\norm{\lVert}{\rVert}%
\makeatletter
\let\oldabs\abs
\def\abs{\@ifstar{\oldabs}{\oldabs*}}
\let\oldnorm\norm
\def\norm{\@ifstar{\oldnorm}{\oldnorm*}}
\makeatother